\documentclass[a4paper,11pt]{amsart}
\usepackage[margin=2.5cm,bottom=2.2cm]{geometry}
\usepackage{hyperref}
\hypersetup{
    colorlinks=false,
    linkbordercolor={1 0 0},
    pdfborder={0 0 1}
}
\usepackage{amssymb}
\usepackage{mathrsfs}
\usepackage{needspace}
\usepackage{tikz-cd}
\tikzcdset{column sep=small,row sep=small}
\usepackage{amsmath}

\newcommand{\R}{\mathbb{R}}
\newcommand{\C}{\mathbb{C}}
\newcommand{\Z}{\mathbb{Z}}
\newcommand{\del}{\partial}
\newcommand{\dbar}{\overline{\partial}}
\newcommand{\calh}{\mathcal{H}}

\newcommand{\tr}{\operatorname{tr}}

\theoremstyle{definition}
\newtheorem{rmk}{Remark}
\newtheorem{exa}{Example}
\theoremstyle{plain}
\newtheorem{lemma}{Lemma}
\newtheorem{cor}{Corollary}
\newtheorem{prop}{Proposition}
\newtheorem{thm}{Theorem}
\newtheorem{introthm}{Theorem}

\author{Tommaso Sferruzza}
\address[Tommaso Sferruzza]{Instituto de Matemática, Estatística e Computação Científica \\
Universidade Estadual de Campinas\\
Rua Sérgio Buarque de Holanda, 651, 13083-859, Campinas, SP, Brasil.}
\email{tommaso.sferruzza@gmail.com}
\title[Formal Hermitian metrics on surfaces]{Geometrically formal Hermitian metrics on compact complex surfaces}
\begin{document}
\raggedbottom
\thanks{The author is partially supported by GNSAGA of INdAM}

\begin{abstract}
We classify compact complex surfaces admitting either a geometrically Dolbeault
formal or an ABC-geometrically formal Hermitian metric. 
Furthermore, we prove that on a compact K\"ahler surface,  every ABC-geometrically formal metric is K\"ahler, whereas, for geometrically Dolbeault formal Hermitian metrics,  Gauduchon, LCK, and
K\"ahler conditions are equivalent. On non-K\"ahler
surfaces, an ABC-geometrically formal Hermitian metric is LCK
and, except on Inoue surfaces, Vaisman.
\end{abstract}
\maketitle

\tableofcontents

\section{Introduction}

A Riemannian metric on a closed oriented manifold is called
\emph{geometrically formal} if the wedge product of harmonic forms is still harmonic \cite{Kot}. Then, the space of harmonic forms provides a chain of quasi-isomorphisms from the de Rham complex to its cohomology algebra, so that the underlying smooth manifold is formal in the sense of Sullivan \cite{Sul77,Kot}.
However, geometric formality is considerably more rigid than its 
homotopy-theoretic counterpart. For instance, in his seminal paper, Kotschick \cite{Kot} obtained bounds on the
Betti numbers and intersection numbers and, in dimensions at most four, proved that a geometrically formal manifold has the real
cohomology algebra of a compact symmetric space. Furthermore, in the complex setting, although every compact K\"ahler manifold is formal \cite{DGMS75},
geometric formality of a K\"ahler metric imposes further restrictions, see, e.g.,
\cite{Huybrechts,GN}. Subsequent work by many authors
has investigated geometric formality under curvature assumptions
\cite{AZ,Bar,Kot2}, on (locally) homogeneous spaces
\cite{GNO,Kas4,KT1,KT2}, and on Vaisman manifolds \cite{OP}.

On a compact complex manifold, Hodge theory associated to the Dolbeault, Bott-Chern, and Aeppli
Laplacians \cite{Sch01} allows to define corresponding notions of geometric formality
\cite{TomTor,AngTom,TaTom17,MS,STprelim}. We focus on two of them.
A Hermitian metric is \emph{geometrically Dolbeault formal} if its
Dolbeault harmonic forms $\mathcal{H}_{\dbar}^{\bullet,\bullet}(X)$ are closed under wedge product. It is
\emph{ABC-geometrically formal} if
$$
\mathcal H^{\bullet,\bullet}(X)
:=\mathcal H_{BC}^{\bullet,\bullet}(X)
+\mathcal H_A^{\bullet,\bullet}(X)
$$
is closed under $\partial$, $\bar\partial$, and wedge product, where $\mathcal H_{BC}^{\bullet,\bullet}(X)
$ and $\mathcal H_A^{\bullet,\bullet}(X)$ are, respectively, the spaces of Bott-Chern and Aeppli harmonic forms.

These conditions represent an important bridge between Hermitian geometry and Dolbeault and
pluripotential homotopy theories \cite{NT78,Ste}. Indeed, they imply the corresponding formality properties and vanishing results for Massey
products \cite{TomTor,AngTom,MS}, providing a way to study the complex homotopy type of a
compact complex manifold through geometric conditions on its Hermitian metrics.
As in the Riemannian case, they also force constant pointwise norms and inner products of the
relevant harmonic forms, yielding cohomological and topological bounds \cite{STprelim}. For a K\"ahler metric,
both conditions agree with Riemannian geometric formality, by the K\"ahler identities.

The main goal of this paper is a complete classification of compact complex surfaces admitting geometrically Dolbeault
formal or ABC-geometrically formal Hermitian metrics. A feature of our results is that, compared to classifications of geometrically formal four-folds \cite{Bar,Kot2}, we are not bound to special
curvature assumptions on possibly different metrics on the same surface. As a second goal, we determine how these geometrically formal Hermitian notions constrain individual metrics, particularly in relation to
the K\"ahler, locally conformally K\"ahler, Vaisman, and Gauduchon
conditions.

Taking a dimensional step back, the situation is completely understood for compact complex curves. Indeed, here every Hermitian metric is K\"ahler and the geometric formality conditions coincide. Moreover, geometric formality forces the compact complex curve to be either  $\mathbb{CP}^1$ with any K\"ahler metric,
or $\mathbb T^1_{\C}$ with a K\"ahler flat metric.

In complex dimension two, we combine the properties descending from geometric Dolbeault or ABC-geometric formality with the special features of compact complex surfaces to determine their cohomology algebra (Propositions \ref{prop:geom-dolb-non-kahler-diamond}, \ref{prop:abc-non-kahler-b1-one-diamond}, \ref{prop:abc-non-kahler-b1-three-diamond}, and \ref{prop:geom-dolb-kahler-b_1=0}) and, by
the  Enriques--Kodaira classification \cite{BPV}, we directly inspect a finite number of families. Nonetheless, it is worth mentioning that several families require additional geometric
arguments, since the determination of the cohomology algebra alone is not sufficient to prove our classification. We present below the results, starting from the non-K\"ahler setting.

\subsection*{Non-K\"ahler surfaces}
A compact complex surface $X=(M,J)$ is non-K\"ahler if and only if its first
Betti number is odd \cite{Buc,Lam}. We call a surface
\emph{cohomologically Hopf} if $b_1(M)=1$ and $b_2(M)=0$.
A Hermitian metric is \emph{locally conformally K\"ahler} (shortly, {\em LCK}) if
$d\omega=\theta\wedge\omega$ for a closed real one-form $\theta$,
called its {\em Lee form}. It is \emph{Vaisman} if $\theta$ is parallel for
the Levi--Civita connection. We prove the following.
\begin{introthm}[Theorems \ref{thm:classification-dolbeault-nonkahler} and \ref{thm:classification-abc-nonkahler}]\label{intro:nonkahler}
Let $X=(M,J)$ be a compact non-K\"ahler complex surface. Then $X$ admits a
geometrically Dolbeault formal Hermitian metric if and only if $X$ is either
one of the following:
\begin{enumerate}
\item a Hopf surface,
\item an Inoue surface,
\item a secondary Kodaira surface,
\item a cohomologically Hopf properly elliptic surface.
\end{enumerate}
$X$ admits an ABC-geometrically formal Hermitian metric if and only if $X$ is either
\begin{enumerate}
\item a finite unramified quotient of a primary diagonal Hopf surface,
\item an LCK Inoue surface,
\item a primary or secondary Kodaira surface,
\item a cohomologically Hopf properly elliptic surface.
\end{enumerate}
Every ABC-geometrically formal metric on $X$ is LCK and is Vaisman
unless $X$ is an Inoue surface. If $X$ is cohomologically Hopf and is
not an Inoue surface, a Hermitian metric on $X$ is
ABC-geometrically formal if and only if it is Vaisman.
\end{introthm}
The existence of a geometrically Dolbeault formal metric on Inoue,
secondary Kodaira, and cohomologically Hopf properly elliptic surfaces
is either directly proved or descends immediately from \cite{AS20,IO}. For Hopf surfaces, we prove that existence
is stable under small deformations within smooth families of Hopf
surfaces (Lemma~\ref{lem:Hopf-openness}), in contrast with the general situation \cite[Theorem 4.3]{TomTor}. Then we use Hasegawa's linearization and
diagonalization families \cite{Hasegawa93,Hasegawa23} to extend
existence from (finite unramified quotients of) primary diagonal Hopf surfaces to
all Hopf surfaces. From the point of view of rigidity, geometric Dolbeault formality alone does not force
the metric to be LCK, as an explicit family on a homogeneous Hopf
surface shows (Example~\ref{ex:hopf-gdf-not-lck}).

For ABC-geometric formality, the algebra
$\mathcal H^{\bullet,\bullet}(X)$ determines the Lee form $\theta$ of the ABC-geometrically formal metric $g$ and forces  $g$ to be LCK. When $b_1(M)=1$, its two possible multiplicative
structures distinguish the LCK Inoue case from the Vaisman case.
In the latter, the differential identities imply the pluricanonical
condition $(\nabla^g\theta)^{(1,1)}\equiv 0$, which forces an LCK metric on a compact manifold to be Vaisman
\cite{MoroianuMoroianu2017} (Lemma \ref{lemma:abc-b1-one-lck}). When $b_1(M)=3$, the surface is primary
Kodaira, and an ABC-geometrically formal metric is Vaisman
and has constant scalar curvature (Proposition \ref{prop:abc-non-kahler-b1-three-diamond}). The nonconstant-scalar-curvature Hermitian metrics of \cite{FGP24} hence provide examples of Vaisman non ABC-geometrically formal metrics on
primary Kodaira surfaces.

\subsection*{K\"ahler surfaces}
We will recall that a  $\mathbb{CP}^1$-bundle over an elliptic curve $E$ is
\emph{projectively unitary} if it is defined by a representation
$\pi_1(E)\to PU(2)$. A \emph{fake quadric} is a minimal surface $X$ of
general type with $h^{1,0}(X)=h^{2,0}(X)=0$ and $c_1^2(X)=8$. We prove the following.

\begin{introthm}[Theorem \ref{thm:classification-kahler}]\label{intro:kahler}
Let $X=(M,J)$ be a compact K\"ahler surface. Then $X$ is geometrically Dolbeault formal if and only if $X$ is either
\begin{enumerate}
\item a complex two-dimensional torus,
\item $\mathbb{CP}^2$ or $\mathbb{CP}^1\times\mathbb{CP}^1$,
\item a fake projective plane,
\item a fake quadric $\Gamma\backslash(\mathbb D\times\mathbb D)$, where
$\Gamma\subset\operatorname{Aut}(\mathbb D)\times
\operatorname{Aut}(\mathbb D)$ is a torsion-free cocompact lattice,
\item a projectively unitary $\mathbb{CP}^1$-bundle over an elliptic curve,
\item a hyperelliptic surface,
\item a quotient $(E\times C)/G$, where $E$ is an elliptic curve,
$g(C)\geq2$, and a finite group $G$ acts freely and diagonally,
by translations on $E$, with $C/G\simeq\mathbb{CP}^1$.
\end{enumerate}
In particular, the following are equivalent:
\begin{enumerate}
\item $X$ admits a geometrically Dolbeault formal Hermitian metric,
\item $X$ admits an ABC-geometrically formal Hermitian metric,
\item $X$ admits a geometrically formal K\"ahler metric.
\end{enumerate}
\end{introthm}
We first show that the Dolbeault
cohomology algebra of $X$ is that of either $\mathbb{T}^2_{\C}$, $\mathbb{CP}^2$,
$\mathbb{CP}^1\times\mathbb{CP}^1$, or
$\mathbb{T}^1_{\C}\times\mathbb{CP}^1$, and  that $X$ is minimal
 (Corollary~\ref{cor:dolb-cohom-ring-kahler-1}). {\em En passant}, combining this minimality result with Theorem \ref{intro:nonkahler}, we have proved that no blow-up of a compact complex surface at a finite number of points is either geometrically Dolbeault or ABC-geometrically formal.

As we mentioned earlier, the determination of the cohomology algebras does not suffice by itself to yield a complete classification. Certain families of compact K\"ahler surfaces from the Enriques-Kodaira classification need new {\em ad hoc} arguments. 

For \emph{Hirzebruch surfaces} $\Sigma^n$, $n\geq 0$, we combine intersection theory with the
holomorphic vector field generating the fiberwise $\mathbb C^*$-action
to prove that $\Sigma^n$ admits a geometrically Dolbeault formal metric
if and only if $n=0$, i.e., $\Sigma^0=\mathbb{CP}^1\times \mathbb{CP}^1$ (Proposition~\ref{prop:no-hirzebruch}). This encompasses the nonexistence result of geometrically formal K\"ahler metrics on odd type Hirzebruch surfaces in
\cite[Theorem~7]{STprelim}.

For \emph{fake quadrics}, symplectic and Seiberg–Witten theories allow us to apply Kotschick’s refinement \cite{Kot97} of Leung’s uniformization theorem \cite{Leung96}, with the small-perturbation chamber made explicit in Lemma \ref{lemma:positive-canonical-chamber}. As a result, $X$ must be  
a quotient $\Gamma\backslash (\mathbb D\times \mathbb D)$ with $\Gamma\subset \text{Aut}(\mathbb D)\times \text{Aut}(\mathbb D)$ a discrete lattice preserving each factor
(Proposition ~\ref{prop:fake-quadric-uniformization}).

For \emph{elliptic ruled surfaces},
we
make use of
nonnegative intersection pairing  
(Lemma~\ref{lem:effective-euler-intersection}), to show that any geometrically Dolbeault formal Hermitian metric is
K\"ahler (Propositions \ref{prop:elliptic-vertical-rigidity} and \ref{thm:elliptic-ruled-metric-rigidity}). Then, the parallelism of
holomorphic one-forms with respect to a geometrically formal K\"ahler metric \cite{GN} makes the bundle projectively
unitary (Corollary \ref{cor:elliptic-ruled-classification}).

For \emph{properly elliptic K\"ahler surfaces}, the study of the Albanese map yields that $X$ is a holomorphic fiber bundle over an elliptic curve, with fibers curves of genus at least two. Isotriviality \cite{VZ} and finite monodromy
give the quotient description in Theorem~\ref{intro:kahler}
(Proposition~\ref{prop:kappa-one-dolbeault}).

\subsection*{Metric rigidity} The proof of the second part of Theorem \ref{intro:kahler} follows by the construction of a geometrically formal K\"ahler metric for every element of the classification. We can also characterize ABC-geometrically formal Hermitian metrics on Kähler surfaces and establish rigidity for geometrically Dolbeault formal metrics under additional assumptions.
Recall that a Hermitian metric on a complex surface is
\emph{Gauduchon} if $\partial\bar\partial\omega=0$.

\begin{introthm}[Propositions \ref{lemma:aeppli-forces-kahler} and \ref{prop:conformal-gdf-uniqueness}, Theorem \ref{thm:kahler-gdf-metric-rigidity}]\label{intro:metric-rigidity}
Let $X$ be a compact K\"ahler surface. Every ABC-geometrically formal
Hermitian metric on $X$ is K\"ahler. For a geometrically Dolbeault
formal Hermitian metric $g$, the following conditions are equivalent:
\begin{enumerate}
\item $g$ is locally conformally K\"ahler,
\item $g$ is Gauduchon,
\item $g$ is K\"ahler.
\end{enumerate}
Moreover, each Hermitian conformal class on $X$ contains at most one
geometrically Dolbeault formal metric up to constant scaling.
\end{introthm}
Thus, ABC-geometric formality is a substantially stronger property than geometric Dolbeault formality on K\"ahler surfaces. To stress this concept, we observe that on $\mathbb{CP}^2$ and fake projective planes, every Hermitian conformal
class contains a unique geometrically Dolbeault formal metric up to
constant scaling (Proposition~\ref{thm:conformal-gdf-b2-one}).
We can use this fact to show that these surfaces admit geometrically Dolbeault formal Hermitian metrics
which are neither Gauduchon nor locally conformally K\"ahler
(Corollary~\ref{cor:nonkahler-gdf-b2-one}), thus showing that the additional
metric assumptions in Theorem~\ref{intro:metric-rigidity} are sharp.

The paper is organized as follows. We set the notation in Section \ref{sec:notation} and, after the Hodge-theoretic and cohomological preliminaries (Section \ref{sec:preliminaries}-\ref{section:geom-dolb-surfaces}),
Sections~\ref{sec:nonkahler-classification}
and~\ref{sec:nonkahler-classification-abc} establish the non-K\"ahler
classifications. Section~\ref{sec:kahler-surfaces} proves the
K\"ahler classification and Section \ref{sec:rigidity} the metric rigidity results.

\section{Notation}\label{sec:notation}
Let $X=(M,J)$ be a compact complex manifold of complex dimension $n$, endowed with a Hermitian structure $(g,\omega)$. We write $\bigwedge^kM$ for $\Lambda^kT^*M$ and $\bigwedge^{p,q}X$ for the bundle of $(p,q)$-forms. We use the complex orientation and the convention
$\omega(U,V)=g(JU,V)$, so $dV_g=\omega^n/n!$.
On real one-forms, $J\alpha=-\alpha\circ J$.
Hermitian inner products are $\C$-linear in the first argument and $\C$-antilinear in the second. We also denote by $g$ the induced Hermitian inner product on forms.
For a $(1,1)$-form $\alpha$, its trace with respect to $\omega$ is
$
\tr_\omega\alpha:=g(\alpha,\omega)=\Lambda_\omega\alpha,
$
where $\Lambda_\omega$ is the adjoint of $L=\omega\wedge\cdot$. In particular, $\text{tr}_{\omega}\omega=n$.
Throughout the paper, $X$ is assumed to be connected. We collect the topological and complex invariants used below, omitting $M$ or $X$ from the notation when no ambiguity arises.
The \emph{Betti numbers} of $M$ are defined as
\begin{align*}
    &b_k(M):=\dim_\R H^{k}(M;\R)=\dim_\C H^{k}(M;\C).
\end{align*}
The \emph{Euler characteristic} of $M$ is given by
\begin{align*}
e(M):=\sum_{i\geq 0} (-1)^ib_i(M).
\end{align*}
The Riemannian metric $g$ induces a Hodge $\star$-operator on the closed oriented manifold $M$. If $\dim_\R M=4m$, the bundle $\bigwedge^{2m}M$ decomposes into the $\pm1$-eigenbundles of $\star$ as
\begin{equation}\label{eq:decomp}
\textstyle\bigwedge^{2m}M=\bigwedge_g^{+
}\oplus\bigwedge^-_g.
\end{equation}
Define $\mathcal{A}_g^{\pm}:=\Gamma(M,\bigwedge^{\pm}_g)$.
If $\dim_{\R}M=4$, we have the explicit expressions
\begin{equation*}
\textstyle\bigwedge_g^+=\R\langle\omega\rangle\oplus (\bigwedge^{2,0}X\oplus\bigwedge^{0,2}X)_{\R}
\end{equation*}
and
\begin{equation}\label{eq:decom-ASD}
\textstyle\qquad \bigwedge^-_g=(\bigwedge^{1,1}_0X)_{\R}:=\{\alpha\in\bigwedge^{1,1}_{\R}X:\tr_{\omega}\alpha=0\}.
\end{equation}
Since $\star$ commutes with the Hodge--de Rham Laplacian $\Delta=d^\star d+dd^{\star}$, decomposition \eqref{eq:decomp} holds at the level of harmonic forms
\[
\mathcal{H}^{2m}(M)=\mathcal{H}_{g}^+\oplus \mathcal{H}_g^-.
\]
We therefore set
\begin{align*}
b^{\pm}(M)=\dim \mathcal{H}^{\pm}(M).
\end{align*}
By Hodge theory, $b^+$ and $b^-$ are the positive and negative indices of the intersection form. In particular, they depend only on the topology and orientation of $M$.
The \emph{Hodge numbers}, \emph{Bott--Chern numbers}, and \emph{Aeppli numbers} of $X$ are
\begin{align*}
h_{\dbar}^{p,q}(X):=\dim_\C H^{p,q}_{\dbar}(X), &\quad h_{BC}^{p,q}(X):=\dim_\C H^{p,q}_{BC}(X), \\
h_{A}^{p,q}(X):=&\dim_\C H^{p,q}_{A}(X).\nonumber
\end{align*}
The \emph{holomorphic Euler characteristic} of $X$ is defined as
\begin{align*}
\chi(X):=\sum_{i\geq 0} (-1)^i h^{0,i}(X).
\end{align*}
We denote by $c_i(X)$ the \emph{$i$-th Chern class} of $T^{1,0}X$.
For a surface, we also write $c_1^2(X)$ and $c_2(X)$ for the
corresponding Chern numbers when they occur in numerical identities.
Let $K_X:=\bigwedge^n(T^{1,0}X)^*$ be the canonical bundle of $X$ and let $P_m(X):=\dim_\C H^{0}(X,K_X^{\otimes m})$ denote its $m$-th plurigenus. Then the \emph{Kodaira dimension} of $X$ is
\[
k(X):=\begin{cases}-\infty, \qquad \qquad \qquad \qquad H^{0}(X,K_X^{\otimes m})=0, \,\,\forall m\geq 1\\
\limsup_{m\to \infty}{\frac{\log P_m}{\log m}}, \qquad \text{otherwise}.
\end{cases}
\]
One has $k(X)\in\{-\infty,0,1,\dots,n\}$. If $\dim_\C X=2$, the topology and complex structure are related by the following standard identities \cite[Chapter~IV]{BPV}. {\em Noether's formula} is
\begin{align}
c_1^2(X)=12\chi(X)-e(M)\label{eq:c_1^2},
\end{align}
and the signature formula is
\begin{align*}
c_1^2(X)=2e(M)+3\sigma(M), \qquad \sigma(M):=b^+(M)-b^-(M).
\end{align*}
Moreover, the \emph{Fr\"olicher spectral sequence} degenerates at the first page, so
\begin{equation*}
b_k(M)=\sum_{p+q=k}h_{\dbar}^{p,q}(X).
\end{equation*}
The top Chern number is the Euler characteristic, so
\begin{align*}
c_2(X)=e(M).
\end{align*}
It will be useful to rewrite Noether's formula as
\begin{align*}
    8h^{0,1}(X)+2(2h^{0,1}(X)-b_1(M))+b_2(M)=10+12h^{0,2}(X)-c_1^2(X).
\end{align*}
Finally, Noether's formula gives $c_1^2(X)+e(M)\equiv0\pmod{12}$. Combining this with the signature formula yields
\begin{equation*}
1-b_1(M)+b^+(M)\equiv0\pmod2.
\end{equation*}
\section{Preliminaries}\label{sec:preliminaries}
Let $X$ be a compact complex manifold endowed with a Hermitian metric
$g$. We recall the \emph{Dolbeault}, \emph{Bott--Chern}, and \emph{Aeppli cohomologies}
\[
H_{\bar\partial}^{\bullet,\bullet}(X)
=\frac{\ker\bar\partial}{\operatorname{im}\bar\partial},
\qquad
H_{BC}^{\bullet,\bullet}(X)
=\frac{\ker\partial\cap\ker\bar\partial}
       {\operatorname{im}\partial\bar\partial},
\qquad
H_A^{\bullet,\bullet}(X)
=\frac{\ker\partial\bar\partial}
       {\operatorname{im}\partial+\operatorname{im}\bar\partial}.
\]
Dolbeault and Bott--Chern cohomologies are bigraded algebras with respect to the cup product, whereas
Aeppli cohomology is a module over $H_{BC}^{\bullet,\bullet}(X)$.
The corresponding spaces of harmonic forms are
\begin{align*}
\mathcal H_{\bar\partial}^{p,q}(X,g)
&=\{\alpha\in\mathcal A^{p,q}(X):
       \bar\partial\alpha=\bar\partial^*\alpha=0\},\\
\mathcal H_{BC}^{p,q}(X,g)
&=\{\alpha\in\mathcal A^{p,q}(X):
       \partial\alpha=\bar\partial\alpha
       =(\partial\bar\partial)^*\alpha=0\},\\
\mathcal H_A^{p,q}(X,g)
&=\{\alpha\in\mathcal A^{p,q}(X):
       \partial\bar\partial\alpha=\partial^*\alpha
       =\bar\partial^*\alpha=0\}.
\end{align*}
They are the kernels of the standard self-adjoint elliptic Dolbeault,
Bott--Chern, and Aeppli Laplacians, respectively, see, e.g., \cite{Sch01}. Hodge theory gives
vector-space isomorphisms
\[
\mathcal H_\sharp^{p,q}(X,g)\cong H_\sharp^{p,q}(X),
\qquad \sharp\in\{\bar\partial,BC,A\}.
\]
When the metric is understood, we omit it from the notation. We now recall the geometric formality notions that we will use in this paper. The Hermitian metric $g$ is said to be
\begin{itemize}
\item \emph{geometrically Dolbeault formal} if
\(\mathcal H_{\bar\partial}^{\bullet,\bullet}(X,g)\) is closed under
wedge product \cite{TomTor},
\item \emph{geometrically Bott--Chern formal} if
\(\mathcal H_{BC}^{\bullet,\bullet}(X,g)\) is closed under wedge
product \cite{AngTom},
\item \emph{ABC-geometrically formal} if
\(\mathcal H_{BC}^{\bullet,\bullet}(X,g)
+\mathcal H_A^{\bullet,\bullet}(X,g)\) is a bidifferential bigraded
subalgebra of \((\mathcal A^{\bullet,\bullet}(X),
\partial,\bar\partial,\wedge)\) \cite{MS},
\item \emph{geometrically Aeppli formal} if
\(\mathcal H_A^{\bullet,\bullet}(X,g)\) is a module over
\(\mathcal H_{BC}^{\bullet,\bullet}(X,g)\) under wedge product
\cite{STprelim}.
\end{itemize}
For K\"ahler metrics, the four formality conditions above agree with each other and with
Riemannian geometric formality, by the K\"ahler identities. We then use the expression
\emph{geometrically formal K\"ahler metric}.

Throughout this paper, we use two distinct extensions of the real Hodge operator $\star$ associated to the Hermitian metric $g$: the \emph{$\C$-antilinear Hodge $\ast$-operator} and the \emph{$\C$-linear Hodge $\star$-operator}, related by $\star=\overline{\cdot}\circ\ast$. Clearly,  $\ast=\star$ on real forms.

The $\C$-antilinear Hodge $\ast$-operator is the isomorphism
\[
\ast\colon \textstyle\bigwedge^{p,q}X\rightarrow \bigwedge^{n-p,n-q}X
\]
defined implicitly by
\[
\alpha\wedge \ast \beta=g(\alpha,\beta)\frac{\omega^{n}}{n!}, \qquad \alpha,\beta\in\textstyle\bigwedge^{p,q}X,
\]
where $g(\alpha,\beta)$ is the extension of the Hermitian metric $g$ to $\bigwedge^{p,q}X$.

The \emph{$\C$-linear $\star$-operator} is defined as the isomorphism
\[
\textstyle\star\colon \bigwedge^{p,q}X\rightarrow \bigwedge^{n-q,n-p}X
\]
defined implicitly by
\[
\alpha\wedge \star\beta=g(\alpha,\overline{\beta})\frac{\omega^n}{n!}, \qquad \alpha\in\textstyle\bigwedge^{q,p}X, \quad \beta\in\bigwedge^{p,q}X.
\]

With these conventions, the formal adjoint of $\dbar$ with respect to the $L^2$-inner product induced by $g$ can be expressed using either Hodge operator:
$$\bar\partial^*=-\ast\bar\partial\ast=-\star\partial\star.$$
Both $\ast$ and $\star$ act on the Dolbeault, Bott--Chern, and Aeppli harmonic forms and induce the following identities:
\begin{align*}
    \ast \mathcal{H}_{\dbar}^{p,q}=\mathcal{H}_{\dbar}^{n-p,n-q}, \qquad \ast\mathcal{H}_{\del}^{p,q}=\mathcal{H}_{\del}^{n-p,n-q}\\
    \ast \mathcal{H}_{BC}^{p,q}=\mathcal{H}_{A}^{n-p,n-q}, \qquad \ast\mathcal{H}_{A}^{p,q}=\mathcal{H}_{BC}^{n-p,n-q}
\end{align*}
and
\begin{align*}
    \star \mathcal{H}_{\dbar}^{p,q}=\mathcal{H}_{\del}^{n-q,n-p}, \qquad \star\mathcal{H}_{\del}^{p,q}=\mathcal{H}_{\dbar}^{n-q,n-p}\\
    \star \mathcal{H}_{BC}^{p,q}=\mathcal{H}_{A}^{n-q,n-p}, \qquad \star\mathcal{H}_{A}^{p,q}=\mathcal{H}_{BC}^{n-q,n-p}.
\end{align*}

\begin{lemma}\label{lemma:n,0-inv}
The bundle $\bigwedge^{n,0}X$ is invariant under $\star$. Moreover,
\[
\mathcal{H}_{\dbar}^{n,0}(X)=\mathcal{H}_{\del}^{n,0}(X)=\mathcal{H}_{BC}^{n,0}(X)=\mathcal{H}_{A}^{n,0}(X)=\mathcal{H}_{dR}^n(M;\C)\cap \mathcal{A}^{n,0}(X),
\]
and similarly
\[
\mathcal{H}_{\dbar}^{0,n}(X)=\mathcal{H}_{\del}^{0,n}(X)=\mathcal{H}_{BC}^{0,n}(X)=\mathcal{H}_{A}^{0,n}(X)=\mathcal{H}_{dR}^n(M;\C)\cap \mathcal{A}^{0,n}(X).
\]
Moreover, these spaces are independent of the metric.
\begin{proof}
The pointwise Hodge-star identity on $(n,0)$-forms is
$$
\star\alpha=\lambda_n\alpha,\qquad
\lambda_n=\begin{cases}1,&n\text{ even},\\-i,&n\text{ odd}.
\end{cases}
$$
For such a form, $\partial\alpha=0$ and $\bar\partial^*\alpha=0$
for bidegree reasons. The identity
$\partial^*=-\star\bar\partial\star$ then shows that
$\partial^*\alpha=0$ if and only if $\bar\partial\alpha=0$.
Thus Dolbeault, conjugate Dolbeault, and de Rham harmonicity all
reduce to holomorphicity. Bott--Chern harmonicity has the same
characterization, since $(\partial\bar\partial)^*\alpha=0$ for
bidegree reasons. For Aeppli harmonicity,
$\partial\bar\partial\alpha=0$ and $\bar\partial^*\alpha=0$
are automatic, while $\partial^*\alpha=0$ again means that
$\alpha$ is holomorphic. This also proves metric independence.
Complex conjugation gives the assertions in bidegree $(0,n)$.
\end{proof}
\end{lemma}
For the remainder of this section, let $X=(M,J)$ be a compact
connected complex surface endowed with a Hermitian metric $g$.
Recall that $g$ is \emph{Gauduchon} if
$\partial\bar\partial\omega=0$. Every Hermitian conformal class
contains a Gauduchon metric, unique up to constant scaling
\cite{Gauduchon1977}.
\begin{lemma}[{\cite[Lemma 3.2]{TarTom23}}]\label{lemma:conform-inv}
    The spaces $\mathcal{H}^{2}_{dR}(X;\R)$ and $\mathcal H_{\bar\partial}^{1,1}(X)$ are conformal invariants of $g$.

\end{lemma}
\begin{lemma}[{\cite[Proposition 3.4]{TarTom23}}]\label{lemma:dolbeault-middle-trace}
Let $\alpha\in\mathcal{A}^{1,1}(X)$. If $\ast\alpha=\alpha$, then
\[
\alpha=f\omega+i\phi, \qquad f\in\mathcal{C}^{\infty}(X)_{\R}, \qquad \phi\in \mathcal{A}_g^-.
\]
If $g$ is Gauduchon, every Dolbeault harmonic $(1,1)$-form
has constant trace; in particular, $f$ is constant whenever
$\alpha$ above is harmonic.

\end{lemma}
\begin{lemma}[{cf. \cite[Proposition 3.8]{TarTom23}}]\label{lemma:asd-dolb}
\[
\mathcal{H}_{g}^-\subset \mathcal{H}_{\sharp}^{1,1}(X), \qquad \sharp \in\{\dbar,BC,A\}.
\]
\begin{proof}
Let $\alpha\in\mathcal{H}_g^-$. By \eqref{eq:decom-ASD}, $\alpha$ is a real anti-self-dual $(1,1)$-form, so $\ast\alpha=\star\alpha=-\alpha$. Since $d\alpha=0$, its pure type implies $\partial\alpha=\bar\partial\alpha=0$. The Hodge-star identities then give $\partial^*\alpha=\bar\partial^*\alpha=0$. Thus $\alpha$ is Dolbeault, Bott--Chern, and Aeppli harmonic.

\end{proof}
\end{lemma}

Let $\bigwedge_g^-$ denote the real rank-three bundle of primitive real
$(1,1)$-forms. The next lemma
relates their existence to the Chern classes of $X$.

\begin{lemma}\label{lem:primitive-complement-chern}
Let $(X,g)$ be a compact Hermitian surface and let $\beta$ be a
nowhere-vanishing real primitive $(1,1)$-form. Set
$\mathscr E:=\langle\beta\rangle^\perp\subset\bigwedge_g^-$,
an oriented real rank-two bundle. Then there is a smooth
splitting into complex line bundles
$$
T^{1,0}X=L_+\oplus L_-.
$$
Writing $x:=c_1(L_+)$ and $y:=c_1(L_-)$, we have
$$
c_1(X)=x+y,\qquad c_2(X)=xy,\qquad
c_1(X)^2-4c_2(X)=(x-y)^2.
$$
In particular,
$$
e(\mathscr E)^2=c_1(X)^2-4c_2(X).
$$
If $\mathscr E$ admits a nowhere-vanishing section, then
$L_+\simeq L_-$ and $c_1(X)^2=4c_2(X)$.
\end{lemma}
\begin{proof}
To every real $(1,1)$-form $\gamma$, we associate the endomorphism $A_\gamma$
defined by
$$
\gamma(u,v)=g(JA_\gamma u,v).
$$
It is self-adjoint and commutes with $J$, so it acts as a Hermitian
endomorphism of $T^{1,0}X$. Primitivity of $\gamma$ is equivalent to
$\operatorname{tr}_{\C}A_\gamma=0$.
Since $\beta$ never vanishes, $A_\beta$ has distinct eigenvalues
$\lambda$ and $-\lambda$, with $\lambda>0$. Its eigenbundles yield the smooth
complex line decomposition  $T^{1,0}(X)=L_+\oplus L_-$. The Whitney formula gives
$c_1(X)=x+y$ and $c_2(X)=xy$.

In a unitary frame adapted to this splitting, a trace-free Hermitian
endomorphism orthogonal to $A_\beta$ has the form
$$
\begin{pmatrix}0&\bar z\\ z&0\end{pmatrix}.
$$
Its lower-left entry represents a complex-linear map $L_+\to L_-$.
This identifies $\mathscr E=\langle\beta\rangle^{\perp}$ with the underlying real bundle of
$L_+^*\otimes L_-$. Hence $e(\mathscr E)=\pm(y-x)$, depending on the
choice of orientation, and $e(\mathscr E)^2=c_1(X)^2-4c_2(X)$. Finally, a nowhere-vanishing section gives a fiberwise nonzero map between
the two complex line bundles, hence an isomorphism. Then $x=y$.
\end{proof}

\section{Geometrically Dolbeault and ABC-geometrically formal manifolds}\label{sec:dolbeault-surfaces}
Let $X=(M,J)$ be a compact complex manifold of complex dimension $n$. Assume that there exists a geometrically Dolbeault formal metric on $X$. We recall the following.
\begin{thm}[{\cite[Proposition~5 and Theorem~1]{STprelim}}]\label{thm:st25}
Dolbeault harmonic forms have constant pointwise norms and pairwise inner products, and the same holds for conjugate Dolbeault harmonic forms. In particular, any basis of $\mathcal{H}_{\dbar}^{\bullet,\bullet}(X)$ (and of $\mathcal{H}_{\del}^{\bullet,\bullet}(X)$) can be orthogonalized with a Gram--Schmidt procedure with constant coefficients. One has
\[
h_{\dbar}^{p,0}(X)\cdot h_{\dbar}^{0,q}(X)\leq h_{\dbar}^{p,q}(X)\leq h_{\dbar}^{p,q}(\mathbb{T}_{\C}^n)
\]
and the analogous inequalities hold for $h_{\del}^{p,q}(X)$.
Moreover,  $b_k(X)\leq b_k(\mathbb{T}^{2n})$.
\end{thm}
\begin{proof}
We recall the argument for completeness. The antilinear Hodge $\ast$ operator preserves the space of Dolbeault harmonic forms. Hence,  $\alpha\wedge\ast\beta=g (\alpha,\beta)\omega^n/n!$ is Dolbeault harmonic
for every $\alpha,\beta\in\mathcal{H}_{\dbar}^{p,q}(X)$. Since the metric is geometrically Dolbeault formal,  $\ast(g(\alpha,\beta)\omega^n/n!)=\overline{g(\alpha,\beta)}$ is holomorphic, hence constant on the compact manifold $X$. Thus, evaluation at any point is
injective, giving $h^{p,q}_{\bar\partial}(X)\leq h^{p,q}(\mathbb{T}_{\C}^n)$.
At a fixed point, products of independent $(p,0)$-forms with independent
$(0,q)$-forms are independent, under the canonical tensor-product
identification with $(p,q)$-forms. Restricting this identification to Dolbeault harmonic forms yields
the lower bound $h_{\dbar}^{p,0}(X)\cdot h_{\dbar}^{0,q}(X)\leq h^{p,q}(X)$. Conjugation gives the assertions for $\partial$.
Finally, from the upper bound and the Fr\"olicher inequality, we have
$
b_k(X)\leq b_k(\mathbb{T}^{2n}).
$
\end{proof}

We observe the following facts.
\begin{lemma}\label{lemma:geom-dolb-1}
$h_{\dbar}^{1,0}(X)+h_{\dbar}^{0,1}(X)\neq 2n-1$.
\begin{proof}
Suppose, for a contradiction, that $h_{\dbar}^{1,0}(X)+h_{\dbar}^{0,1}(X)=2n-1$. By Theorem~\ref{thm:st25}, both $h_{\dbar}^{1,0}(X)$ and $h_{\dbar}^{0,1}(X)$ are at most $n$. Consider first the case $h_{\dbar}^{1,0}(X)=n$ and $h_{\dbar}^{0,1}(X)=n-1$. Suppose that $\{\alpha_1,\dots,\alpha_n\}$ and $\{\beta_1,\dots, \beta_{n-1}\}$ are orthonormal bases, respectively, of $\mathcal{H}_{\dbar}^{1,0}(X)$ and $\mathcal{H}_{\dbar}^{0,1}(X)$. Since the metric is geometrically Dolbeault formal,
\[
0\neq \alpha_1\wedge\dots\wedge \alpha_n\wedge\beta_1\wedge\dots\wedge\beta_{n-1}\in\mathcal{H}_{\dbar}^{n,n-1}(X).
\]
Moreover, the $\C$-antilinear $\ast$ is an automorphism of $\mathcal{H}_{\dbar}^{\bullet,\bullet}(X)$, hence the form
\[
\ast(\alpha_1\wedge\dots\wedge \alpha_n\wedge\beta_1\wedge\dots\wedge\beta_{n-1})\in\mathcal{H}_{\dbar}^{0,1}(X)
\]
is pointwise orthogonal to $\beta_1,\dots,\beta_{n-1}$ and nonzero. This contradicts $h_{\dbar}^{0,1}(X)=n-1$. The case $h_{\dbar}^{1,0}(X)=n-1$ and $h_{\dbar}^{0,1}(X)=n$ is analogous.
\end{proof}
\end{lemma}
\begin{lemma}\label{lemma:geom-dolb-2}
If $h_{\dbar}^{1,0}(X)=\dim_\C X$, then $X$ is complex parallelisable. If, moreover, $X$ is a compact complex surface, then $X$ is a complex torus and the metric is flat and K\"ahler.
\begin{proof}
Since $h_{\dbar}^{1,0}(X)=n$, Theorem~\ref{thm:st25} gives a global orthonormal coframe of holomorphic $(1,0)$-forms. Thus the holomorphic cotangent bundle is trivial. By Wang's theorem \cite{Wan56}, $X=G/\Gamma$, where $G$ is a simply connected complex Lie group and $\Gamma$ is a cocompact lattice.

Suppose now that $\dim_{\C}X=2$. Then $G$ has complex dimension two. Since $G$ must be unimodular, this forces the Lie algebra of $G$ to be abelian, i.e., $G\simeq\C^2$ and $X$ is a complex torus. By Theorem \ref{thm:st25}, we may choose a global orthogonal coframe $\{\alpha_1,\alpha_2\}$ of holomorphic closed forms with $|\alpha_i|_g^2=2$, so that
\[
g=\alpha_1\odot  \overline{\alpha_1}+\alpha_2\odot  \overline{\alpha_2},
\qquad
\omega=\frac{i}{2}(\alpha_1\wedge\overline{\alpha_1}+\alpha_2\wedge\overline{\alpha_2}).
\]
Since $d\alpha_1=d\alpha_2=0$, we have $d\omega=0$. Thus, the geometrically Dolbeault formal metric $g$ is K\"ahler.
Let $\pi\colon\C^2\to X$ be the universal covering map. In the standard coordinates on $\C^2$,
\[
\pi^*\alpha_i=\sum_{j=1}^2a_{ij}dz_j, \qquad a_{ij}\in\C,
\]
since holomorphic $1$-forms are invariant under translations on complex tori.
In particular, the coefficients of $\pi^*g$ on $\C^2$ are constant. Hence, $\mathcal{R}^g\equiv 0$ and $g$ is flat.
\end{proof}
\end{lemma}
\begin{lemma}\label{lemma:geom-dolb-3}
If either $h^{1,0}_{\dbar}(X)\neq 0$, or $h_{\dbar}^{0,1}(X)\neq 0$, then $e(M)= 0$.
\begin{proof}
Choose a nonzero Dolbeault harmonic $(1,0)$-form $\alpha=u+iJu$. Its norm
is a positive constant, so $u=\mathfrak{Re}(\alpha)$ is nowhere vanishing, since a zero of
$u$ would also be a zero of $\alpha$. Thus $u^\sharp$ is a
nowhere-vanishing real vector field, and the classical Poincar\'e--Hopf
theorem gives $e(M)=0$. The same argument applies to a nonzero Dolbeault 
harmonic $(0,1)$-form.
\end{proof}
\end{lemma}
Assume now that $X$ admits an ABC-geometrically formal metric. We recall the following.
\begin{thm}[{\cite[Proposition~6, Theorem~2 and Theorem~6]{STprelim}}]\label{thm:geom-BC-A}
Aeppli harmonic forms have constant pointwise norms and pairwise inner products, and the same holds for Bott--Chern harmonic forms. In particular, any basis of $\mathcal{H}_{A}^{\bullet,\bullet}(X)$ and $\mathcal{H}_{BC}^{\bullet,\bullet}(X)$ can be orthogonalized with a Gram--Schmidt procedure with constant coefficients. One has
\[
h_{BC}^{p,0}(X)\cdot h_{BC}^{0,q}(X)\leq h_{BC}^{p,q}(X)\leq h_{BC}^{p,q}(\mathbb{T}_{\C}^n), \quad h_{A}^{p,q}(X)\leq h_{A}^{p,q}(\mathbb{T}_\C^n), \quad b_k(M)\leq b_k(\mathbb{T}^{2n}).
\]
Moreover, if $h_{BC}^{1,0}(X)=n$, then $X$ is biholomorphic to a complex torus of dimension $n$ and the ABC-geometrically formal metric is K\"ahler and flat.
\end{thm}
\begin{proof}
Set $\mathcal H=\mathcal H_{BC}+\mathcal H_A$.
Both $\ast$ and wedge product preserve $\mathcal H$.
In bidegree $(0,0)$, both harmonic spaces consist of constants:
for Aeppli harmonicity this follows from compactness and
$\partial\bar\partial f=0$. By Hodge duality,
$\mathcal H^{n,n}=\mathbb C\langle\omega^n\rangle$. Therefore
$$
\alpha\wedge\ast\beta=\langle\alpha,\beta\rangle_g\,\omega^n/n!
\in\mathcal H^{n,n}
$$
forces $\langle\alpha,\beta\rangle_g$ to be constant for
$\alpha,\beta\in\mathcal H^{p,q}$. Evaluation at each point is
injective, so, more generally,
$
\dim_\C\mathcal{H}^{p,q}\leq h^{p,q}(\mathbb{T}^n_{\C}).
$
This proves the upper bounds for both summands.
Products of Bott--Chern harmonic forms are Bott--Chern harmonic. Indeed, if $\alpha,\beta\in\mathcal{H}^{\bullet,\bullet}_{BC}(X)$, then $\alpha\wedge \beta=\gamma+\delta$, with
$\gamma\in\mathcal H_{BC}$ and $\delta\in\mathcal H_A$. Then
$d\delta=0$, and since $d^*\delta=0$, we have $(\del\dbar)^*(\alpha\wedge \beta)=0$, hence $\alpha\wedge \beta$ is
Bott--Chern harmonic \cite{MS}. The same argument as in
Theorem~\ref{thm:st25} gives the lower bound. Finally, combining the upper bound and the  inequality in \cite[Theorem A]{AngTom13} gives the upper bound on Betti numbers.

If $h_{BC}^{1,0}(X)=n$, constant orthonormalization gives a global
coframe of closed holomorphic one-forms. Its real and imaginary
parts are a closed orthonormal coframe, so the metric is flat. Since $X$ is compact, integrating the coframe on the universal cover identifies it with
$\mathbb C^n$, and deck transformations act by translations.
Thus $X$ is a complex torus and $g$ is flat and K\"ahler.
\end{proof}

As in Lemma~\ref{lemma:geom-dolb-3}, we obtain the following consequence.

\begin{lemma}\label{lemma:abc-geom-1}
    If either $h_{BC}^{1,0}(X)\neq 0$ or $h_A^{1,0}(X)\neq 0$, then $e(M)=0.$
\end{lemma}

We shall also use the following comparison.
\begin{lemma}[{\cite[Lemma~1]{STprelim}}]\label{lemma:abc-harmonic-comparison}
Let $X=(M,J)$ be a compact complex manifold with a Hermitian metric $g$.
Then the following are equivalent:
\begin{enumerate}
    \item $(M,J,g)$ is geometrically Aeppli formal,
    \item $(M,J,g)$ is ABC-geometrically formal and satisfies the $\del\dbar$-lemma,
    \item the Dolbeault, Bott--Chern, and Aeppli harmonic spaces
    coincide. Writing their common value as
    $\mathcal H^{\bullet,\bullet}(M)$, one has
    $\mathcal H_{dR}^{k}(M;\C)=\bigoplus_{p+q=k}\mathcal H^{p,q}(M)$,
    and this space is closed under wedge product.
\end{enumerate}
In particular, if any of the above conditions holds, the metric $g$ is geometrically Dolbeault formal and
geometrically formal in the Riemannian sense.
\end{lemma}
\begin{proof}
Under (1), multiplication by $1\in\mathcal H_A^{0,0}$ gives
$\mathcal H_{BC}^{p,q}\subseteq\mathcal H_A^{p,q}$ for every $p,q$. Since $\ast$ is an isomorphism, we have  $\mathcal{H}_A^{p,q}=\ast\mathcal{H}_{BC}^{n-p,n-q}\subseteq \ast\mathcal{H}_{A}^{n-p,n-q}=\mathcal{H}_{BC}^{p,q}$, so $\mathcal{H}_{BC}^{p,q}=\mathcal{H}_A^{p,q}$.
The map induced by the identity $H_{BC}\to H_A$ is then an isomorphism represented by the identity map between $\mathcal{H}_{BC}$ and $\mathcal{H}_A$, hence $X$ satisfies the
$\partial\bar\partial$-lemma. The space $\mathcal{H}$ collapses to $\mathcal{H}_{BC}$, which satisfies $\del|_{\mathcal{H}}\equiv \dbar|_{\mathcal{H}}\equiv 0$ and has a structure of algebra induced by the module structure, proving (2).
\smallskip

Assume (2).
For $\alpha\in\mathcal H$, one has $\partial\bar\partial\alpha=0$.
The space $\mathcal{H}$ is closed under both $\del$ and $\dbar$, hence $\partial\alpha\in\mathcal H$. This form is
$d$-closed and $\partial$-exact, so the $\partial\bar\partial$-lemma
makes it $\partial\bar\partial$-exact. But both summands $\mathcal{H}_{BC}$ and $\mathcal{H}_A$ of
$\mathcal H$ are orthogonal to $\operatorname{im}\partial\bar\partial$.
Hence $\partial\alpha=0$, and similarly $\bar\partial\alpha=0$.
Every Aeppli harmonic form is therefore Bott--Chern harmonic.
From the $\del\dbar$-lemma, $H_{BC}$ and $H_A$ have the same dimension, thus
$\mathcal H_A=\mathcal H_{BC}$, and
$\mathcal H=\mathcal H_{\bar\partial}=\mathcal H_{BC}=\mathcal H_A$.
Indeed, the common Bott--Chern and Aeppli harmonic forms are
Dolbeault harmonic, and the $\partial\bar\partial$-lemma gives
equality of dimensions. Each is also de Rham harmonic. The Hodge decomposition supplied by
the $\partial\bar\partial$-lemma yields that $\mathcal{H}_{dR}^{k}(X;\C)=\bigoplus_{p+q=k}\mathcal{H}^{p,q}$, and hence (3).
\smallskip

The implication (3)$\Rightarrow$(1) is immediate.
\end{proof}

\section{Cohomological preliminaries for compact complex surfaces}\label{section:geom-dolb-surfaces}
We collect the standard numerical identities used in the surface classifications \cite[Chapter~IV]{BPV}. Let $X=(M,J)$ be a compact complex surface. The Fr\"olicher spectral sequence degenerates at the first page, so
\begin{align}
b_1(M)&=h_{\dbar}^{1,0}(X)+h_{\dbar}^{0,1}(X), \label{eq:1}\\
b_2(M)&=h_{\dbar}^{2,0}(X)+h_{\dbar}^{1,1}(X)+h_{\dbar}^{0,2}(X),\notag
\end{align}
Moreover,
\begin{equation}
2h_{\dbar}^{0,1}(X)\geq b_1(M). \label{eq:3}
\end{equation}
\begin{lemma}\label{lemma:geom-dolb-4}
If $X=(M,J)$ is a compact complex surface admitting a geometrically Dolbeault formal metric, then
\[
(h_{\dbar}^{1,0}(X),h_{\dbar}^{0,1}(X))\in\{(0,0),(0,1),(1,1),(2,2)\},
\]
\begin{proof}
The bounds in Theorem~\ref{thm:st25} and Lemma~\ref{lemma:geom-dolb-1} leave the possibilities
\[
(h_{\dbar}^{1,0}(X),h_{\dbar}^{0,1}(X))\in\{(0,0),(1,0),(0,1),(2,0),(1,1),(0,2),(2,2)\}.
\]
If $(h_{\dbar}^{1,0}(X),h_{\dbar}^{0,1}(X))$ were $(0,2)$ or $(2,0)$, then $b_1(M)=2$ and $X$ would be K\"ahler, contradicting Hodge symmetry.
Moreover, by \eqref{eq:3}, we cannot have $(h_{\dbar}^{1,0}(X),h_{\dbar}^{0,1}(X))=(1,0)$, otherwise $b_1(M)=1>2h_{\dbar}^{0,1}(X)=0$.
\end{proof}
\end{lemma}
So, we are left with the cases
\[
(h_{\dbar}^{1,0}(X),h_{\dbar}^{0,1}(X))\in\{(0,0),(0,1),(1,1),(2,2)\}.
\]
If we assume that $(h_{\dbar}^{1,0}(X),h_{\dbar}^{0,1}(X))=(2,2)$, then by Lemma \ref{lemma:geom-dolb-2}, $X$ is a complex  torus and the metric $g$ is flat and K\"ahler. We consider the remaining cases, i.e.,
\[
(h_{\dbar}^{1,0}(X),h_{\dbar}^{0,1}(X))\in
\begin{cases}
\{(0,1)\}, \qquad\qquad \text{$X$ is non-K\"ahler,}\\
\{(0,0),(1,1)\},\,\,\quad \text{$X$ is K\"ahler}.
\end{cases}
\]

On the other hand, the Bott--Chern diamond of a compact complex surface can be described in terms of the Betti and intersection numbers of the underlying topological manifold, depending on the parity of the first Betti number, see, e.g., \cite{SteWil}. If $b_1$ is even, i.e., $X$ is K\"ahler, then

{\small
\[
\begin{tikzcd}
  & & 1 & &\\
  &\frac{b_1}{2} & &  \frac{b_1}{2}&\\
  \frac{b^+-1}{2} && b^-+1 && \frac{b^+-1}{2}\\
  &\frac{b_1}{2} & &  \frac{b_1}{2}&\\
  &&1&&
\end{tikzcd}
\]}
Conversely, if $b_1$ is odd, i.e., $X$ is not K\"ahler, we have
{\small
\begin{equation}\label{eq:b+-hbc}
\begin{tikzcd}
  & & 1 & &\\
  &\frac{b_1-1}{2} & &  \frac{b_1-1}{2}&\\
  \frac{b^+}{2} && b^-+1 && \frac{b^+}{2}\\
  &\frac{b_1+1}{2} & &  \frac{b_1+1}{2}&\\
  &&1&&
\end{tikzcd}
\end{equation}}

\begin{lemma}\label{lemma:h(1,0)-abc}
If $X=(M,J)$ is a compact complex surface admitting an ABC-geometrically formal metric, then
\[
(h_{BC}^{1,0}(X),h_{A}^{1,0}(X))\in\{(0,0),(0,1),(1,1),(1,2),(2,2)\}.
\]
\begin{proof}
By Theorem~\ref{thm:geom-BC-A}, both $h_{BC}^{1,0}$ and
$h_A^{1,0}$ are at most two. The displayed Bott--Chern diamonds
and Bott--Chern--Aeppli duality therefore give $b_1\leq4$ in
the K\"ahler case and $b_1\leq3$ in the non-K\"ahler case.
The five possible values of $b_1$ give the stated pairs.
\end{proof}
\end{lemma}
We observe that
\[
(h^{1,0}_{BC}(X),h_{A}^{1,0}(X))\in\begin{cases}
    \{(0,0),(1,1),(2,2)\},\quad  X\,\,\text{is K\"ahler}\\
    \{(0,1),(1,2)\},\qquad \,\,\,\,\quad  X\,\,\text{is non-K\"ahler}.
\end{cases}
\]
If $h_{BC}^{1,0}(X)=2$, Theorem~\ref{thm:geom-BC-A} shows that $X$ is a complex torus. We therefore focus on the remaining cases.

\section{Geometrically Dolbeault formal non-K\"ahler surfaces: a classification}\label{sec:nonkahler-classification}
In this section,  we classify compact non-K\"ahler surfaces admitting geometrically Dolbeault formal metrics. We start by characterizing the Dolbeault cohomology algebra of such surfaces.

\subsection{Dolbeault cohomology algebra}
Let $X$ be a compact non-K\"ahler surface admitting a geometrically Dolbeault formal metric. By Lemma~\ref{lemma:geom-dolb-4}, $h_{\dbar}^{1,0}=0$ and $h_{\dbar}^{0,1}=1$. By Lemma \ref{lemma:geom-dolb-3} and \eqref{eq:1},
$$
0=e(M)=b_0-b_1+b_2-b_3+b_4=b_2.
$$
so that the only nonvanishing Hodge numbers of $X$ are
\[
h^{0,0}_{\dbar}=h_{\dbar}^{0,1}=h_{\dbar}^{2,1}=h_{\dbar}^{2,2}=1.
\]
In the paper, we will refer to a compact complex surface with this Hodge diamond as \emph{cohomologically Hopf}.

Since a non-K\"ahler surface has Kodaira dimension at most one, we obtain the following.
\begin{prop}\label{prop:geom-dolb-non-kahler-diamond}
Let $X$ be a compact complex non-K\"ahler surface admitting a geometrically Dolbeault formal metric. Then, the Hodge diamond of $X$ must be

{\small
\begin{equation}\label{dolb-alg-non-kahl-gd}
\begin{tikzcd}
  & & 1 & &\\
  & 0 &  & 1 &\\
  0 & & 0 & & 0 \\
  & 1 &  & 0 & \\
  & & 1 & &
\end{tikzcd}
\end{equation}}

The Dolbeault cohomology algebra is represented by $\C\langle1,\alpha,\ast\alpha,\alpha\wedge\ast\alpha\rangle$, where $0\neq\alpha\in\mathcal{H}_{\dbar}^{0,1}(X)$. Furthermore, $k(X)\in\{-\infty,0,1\}$.
\end{prop}

\subsection{Geometrically Dolbeault formal surfaces}
We determine which compact non-K\"ahler surfaces admit a geometrically Dolbeault formal metric, via their Kodaira dimension, see, e.g., \cite{BPV}.

\subsubsection{$k(X)=-\infty$}
Any non-K\"ahler surface with $k(X)=-\infty$ is a class VII surface, and those with $b_2=0$ are Hopf or Inoue surfaces \cite{Bog76,LYZ90}. Every invariant Hermitian metric on an Inoue surface is geometrically Dolbeault formal \cite[\S5]{AS20}, whereas it is proved in \cite[Theorem 6.6]{IO} that diagonal Hopf surfaces and their finite quotients are geometrically Dolbeault formal, since they are Vaisman and cohomologically Hopf.
We extend this result to all Hopf surfaces.

We first prove a normalization lemma.
\begin{lemma}\label{lem:Hopf-normalization}
Let $X=(M,J)$ be a cohomologically Hopf compact complex surface. If $X$ admits a Hermitian metric $g$ for which $\eta\in\calh^{0,1}_{\dbar}(X)$ is nowhere-vanishing, then $X$
admits a geometrically Dolbeault formal metric.
\end{lemma}

\begin{proof}
Let $\omega$ be the fundamental form of $g$ and set
$\ell=|\eta|_g^2>0$. Define
\[
\tilde{\omega}:=\omega-i\left(1-\frac{1}{\ell}\right)
\eta\wedge\overline{\eta}.
\]
In a local $g$-unitary $(1,0)$-coframe with
$\eta=\sqrt{\ell}\,\overline{\theta^1}$,
\[
\tilde\omega
=i\ell\,\theta^1\wedge\overline{\theta^1}
+i\theta^2\wedge\overline{\theta^2}>0.
\]
Hence $\tilde\omega$ is the fundamental form of a Hermitian metric $\tilde{g}$.
We next show that $\eta$ remains Dolbeault harmonic for $\tilde g$. Since $\eta$ is a primitive $(0,1)$-form, the Hodge-star identity gives $\ast_g\eta=-i\bar\eta\wedge\omega$. Thus, since
\[
\ast_{\tilde{g}}\eta=-i \bar\eta\wedge\tilde\omega
 =-i\bar\eta\wedge\omega=\ast_g\eta,
\]
we have $\dbar_g^*\eta=0$ if and only if $\dbar^{*}_{\tilde g}\eta=0$, hence
$\eta\in\calh^{0,1}_{\dbar}(X,\tilde g)$ and $\mathcal{H}_{\dbar}^{2,1}(X,\tilde{g})=\C\langle\ast_{\tilde{g}}\eta\rangle$. Furthermore, since $|\eta|_{\tilde{g}}^2\equiv 1$, we have $2\eta\wedge \ast_{\tilde{g}}\eta=\tilde{\omega}^2$, which is harmonic because it is twice the volume form of $\tilde g$. Hence, $\tilde{g}$ is geometrically Dolbeault formal.
\end{proof}

We next prove a stability result for Hopf surfaces. Recall that a Hopf surface is defined as $X:=(\C^2\setminus\{0\})/G$, where $G=H\rtimes \langle f\rangle$, with $f$ a contraction, $H$ a finite group on which $\langle f\rangle$ acts by conjugation. If $H$ is trivial, then the Hopf surface $X$ is {\em primary}, otherwise, it is {\em secondary}. Every element of $G$ extends
holomorphically across the origin. Hasegawa's linearization construction associates to $G$ its
derivative group, which again acts freely and properly discontinuously
on $\C^2\setminus\{0\}$
\cite[Theorem~3.3]{Hasegawa93}; see also \cite[Section~3]{Hasegawa23}. Denote by $L(G)=H_{\text{lin}}\rtimes \langle f\rangle_{\text{lin}}$ the linear part of $G$. The Hopf manifold $X$ is said to be {\em linear} if $G=L(G)$.
\begin{lemma}\label{lem:Hopf-openness}
Let $\pi\colon\mathcal X\to\Delta$ be a smooth proper holomorphic family all
of whose fibers are Hopf surfaces.  If $X_0$ admits a geometrically Dolbeault
formal metric $g_0$, then every sufficiently small deformation $X_t$ is still geometrically Dolbeault formal.
\end{lemma}

\begin{proof}
By Ehresmann's theorem, up to shrinking $\Delta$, we may identify the fibers
differentiably and choose a smooth family of Hermitian metrics $g_t$ whose
value at $t=0$ is $g_0$.  Every Hopf
surface has $h^{0,1}_{\dbar}=1$, see, e.g., \cite[Chapter V]{BPV}.
Therefore, since the kernels of the self-adjoint elliptic operators
$\Delta_{\dbar,t}\colon\mathcal A^{0,1}(X_t)\to
\mathcal A^{0,1}(X_t)$ have constant dimension, the
projection map $P_t\colon \mathcal{A}^{0,1}(X_t)\rightarrow \mathcal{H}_{\dbar}^{0,1}(X_t,g_t)$ depends smoothly on $t$ \cite[Theorem 4.5]{KodMor}. Choose $0\ne\eta_0\in\calh^{0,1}_{\dbar}(X_0,g_0)$ and a smooth extension $\tilde\eta_t$ to the family. Set $\eta_t:=P_t(\tilde\eta_t)$. Then $\eta_t$ is a smooth family of Dolbeault harmonic forms with value $\eta_0$ at $t=0$.

Since $g_0$ is geometrically Dolbeault formal,
$|\eta_0|_{g_0}$ is a positive constant.  Compactness and smoothness
then imply that $\eta_t$ is nowhere-vanishing for all sufficiently small
$t$.  Lemma \ref{lem:Hopf-normalization} then applies on $X_t$, completing the
proof.
\end{proof}
We make use of the deformation stability of Lemma \ref{lem:Hopf-openness} to obtain the following general result.

\begin{prop}\label{prop:all-Hopf-GDF}
Every Hopf surface admits a geometrically Dolbeault formal
Hermitian metric.
\end{prop}

\begin{proof}
Let $X=(\mathbb C^2\setminus\{0\})/G$ be any Hopf surface. Hasegawa constructs two proper
holomorphic families related to $X$, the first by the linearization  $G\mapsto L(G)$
\cite[Theorem~3.3]{Hasegawa93}, and the second by diagonalizing its
cyclic part \cite[proof of Theorem~3.7]{Hasegawa93}. In particular, there exist
\[
\pi_1\colon \mathcal X_1\longrightarrow\C,
\qquad
\pi_2\colon \mathcal X_2\longrightarrow\C,
\]
such that
\[
\pi_1^{-1}(0)=X_{\mathrm{lin}},
\qquad
\pi_1^{-1}(t)\simeq X\quad(t\neq0),
\]
and
\[
\pi_2^{-1}(0)=X_{\mathrm{diag}},
\qquad
\pi_2^{-1}(t)\simeq X_{\mathrm{lin}}\quad(t\neq0),
\]
where $X_{\mathrm{lin}}$ is linear, and $X_{\mathrm{diag}}$ is a finite unramified quotient of a
primary diagonal Hopf surface. In the first construction, for $t\ne0$, set
$T_t(z)=tz$ and $G_t=T_t^{-1}GT_t$.
For each $g\in G$, the map $T_t^{-1}gT_t$ extends to $(Dg)_0$
at $t=0$, giving $G_0=L(G)$. Hasegawa's construction gives a proper holomorphic submersion $\mathcal{X}_1\rightarrow \C$, and for
$t\ne0$, the action of $G_t$ on $\C^2\setminus\{0\}$ is conjugate to that of $G$, i.e., $X_t$ is biholomorphic to  $X$.

For clarity, the second family can be described directly in dimension two.
Put $A=(Df)_0$. If $A$ is diagonalizable, choose diagonal coordinates and
define $\mathcal{X}_2\rightarrow \C$ as the constant family. Otherwise, set $A=\lambda I+N$ with $N^2=0\ne N$.
A power $A^m$ centralizes the finite group $H_{\mathrm{lin}}$ because
conjugation by $A$ permutes this group. Thus every element of
$H_{\mathrm{lin}}$ commutes with $N$ and has the form $aI+bN$.
Finite order forces $b=0$, so $H_{\mathrm{lin}}$ acts by scalar multiplication.
In a Jordan basis, put $T_t=\operatorname{diag}(1,t)$. Then
\[
A_t=T_t^{-1}AT_t=\begin{pmatrix}\lambda&t\\0&\lambda\end{pmatrix},
\]
while leaving $H_{\mathrm{lin}}$ fixed. At $t=0$ the contraction is
diagonal, and for $t\ne0$ the whole group is conjugate to $L(G)$.
This is the second proper family above, including the secondary case.

Now, by Belgun's classification \cite{Bel},
$X_{\mathrm{diag}}$ admits a Vaisman metric, which is
geometrically Dolbeault formal by \cite[Theorem 6.6]{IO}.

Up to restricting to a small enough neighborhood $\Delta_2\subset \C$ of $0$, we can apply Lemma \ref{lem:Hopf-openness} to the family
$\pi_2\colon \pi_2^{-1}(\Delta_2)\to\Delta_2$, so that
$X_{\mathrm{lin}}$ admits a geometrically Dolbeault formal
metric. Restricting to a small neighborhood $\Delta_1\subset \C$ of $0$, we apply once again Lemma \ref{lem:Hopf-openness} to
$\pi_1\colon\pi_1^{-1}(\Delta_1)\to\Delta_1$, which finally shows that $X$ is geometrically Dolbeault formal.
\end{proof}

\subsubsection{$k(X)=0$}
Cohomologically  Hopf surfaces of Kodaira dimension zero are precisely the secondary Kodaira surfaces, for which every invariant Hermitian metric is geometrically Dolbeault formal \cite[\S5]{AS20}.

\begin{rmk}
Lemma \ref{lem:Hopf-normalization} recovers the case of Inoue surfaces and secondary Kodaira surfaces. Indeed, they all admit a structure of (finite unramified quotients of) solvmanifolds endowed with an invariant complex structure and  the space of Dolbeault harmonic forms with respect to any invariant metric is determined by the complex of invariant forms, which have constant norm \cite[\S\S4--5]{AS20}. It is worth mentioning that a stability result was obtained in \cite[Theorem 4.2]{TaTom17} for geometrically Bott-Chern formal metrics on Hopf surfaces.
\end{rmk}

\subsubsection{ $k(X)=1$}\label{subsub:k(x)=1-gd}
$X$ is a properly elliptic surface. Since $X$ realizes the Hodge diamond \eqref{dolb-alg-non-kahl-gd}, it is cohomologically Hopf. Then, by Belgun's classification \cite[Theorem 1]{Bel},  $X$ admits a Vaisman metric $g$ which is geometrically Dolbeault formal by \cite[Theorem 6.6]{IO}, since $X$ is cohomologically Hopf. \smallskip

To summarize, geometrically Dolbeault formal non-K\"ahler surfaces are classified as follows.
\begin{thm}\label{thm:classification-dolbeault-nonkahler}
Let $X=(M,J)$ be a compact non-K\"ahler surface. Then $X$ is geometrically Dolbeault formal if and only if $X$ is either
\begin{enumerate}
    \item a Hopf surface,
    \item an Inoue surface,
    \item a secondary Kodaira surface,
    \item a cohomologically Hopf properly elliptic surface.
\end{enumerate}
\end{thm}

The following example provides a smooth family of geometrically Dolbeault
formal metrics on a homogeneous Hopf surface whose only locally
conformally K\"ahler member is the initial metric. Thus, geometric
Dolbeault formality does not imply LCK on non-K\"ahler surfaces.

\begin{exa}
\label{ex:hopf-gdf-not-lck}
On the homogeneous Hopf surface $S^3\times S^1$, choose an invariant
real coframe with $de^1=-2e^2\wedge e^3$, $de^2=2e^1\wedge e^3$,
$de^3=-2e^1\wedge e^2$, $de^4=0$, and $(1,0)$-forms $e^1+ie^2,e^3+ie^4$.
For $t\in(-1,1)$, the forms
\[
 \omega_t=e^1\wedge e^2+e^3\wedge e^4+t(e^1\wedge e^3+e^2\wedge e^4)
\]
are positive, with Hermitian eigenvalues $1\pm t$.
The form $\eta=e^3-ie^4$ is $\bar\partial$-closed, and
$\bar\partial_t^*\eta$ is an invariant function with integral zero,
hence vanishes. A basis of Dolbeault harmonic forms with respect to $\omega_t$ is given by  $1,\eta,\ast_t\eta,\omega_t^2$. Since $|\eta|_{g_t}$ is constant,
its only nonzero positive-degree product,
$2\eta\wedge\ast_t\eta=|\eta|_{g_t}^2\omega_{t}^2$, is harmonic.
Thus $g_t$ is geometrically Dolbeault formal. Its Lee form satisfies
\[
 \theta_t=\frac{4t}{1-t^2}e^1-\frac{2(1+t^2)}{1-t^2}e^4,
 \qquad d\theta_t=-\frac{8t}{1-t^2}e^2\wedge e^3,
\]
so $g_t$ is not locally conformally K\"ahler when $t\ne0$.
\end{exa}

\section{ABC-geometrically formal non-K\"ahler surfaces: a classification}\label{sec:nonkahler-classification-abc}

We now classify the non-K\"ahler surfaces that admit ABC-geometrically formal metrics.
\subsection{Aeppli--Bott--Chern cohomology algebra}
We first determine the structure of $\mathcal H=\mathcal H_{BC}+\mathcal H_A$.

For an ABC-geometrically formal non-K\"ahler surface, Lemma \ref{lemma:h(1,0)-abc}
gives $b_1\in\{1,3\}$, corresponding respectively to
$(h^{1,0}_{BC},h^{1,0}_A)=(0,1)$ and $(1,2)$.
We first treat $b_1=1$, when the two possible harmonic algebras lead
to the distinction between Inoue surfaces and Vaisman surfaces.

\begin{prop}\label{prop:abc-non-kahler-b1-one-diamond}
Let $X=(M,J)$ be a compact non-K\"ahler complex surface endowed
with an ABC-geometrically formal Hermitian metric $g$. Assume that
$b_1(M)=1$. Then $b_2(M)=0$, and the Bott--Chern and Aeppli
diamonds of $X$ are respectively
{\small \[
\begin{tikzcd}[column sep=small,row sep=small]
 &&1&&\\
 &0&&0&\\
 0&&1&&0\\
 &1&&1&\\
 &&1&&
\end{tikzcd},
\qquad \quad 
\begin{tikzcd}[column sep=small,row sep=small]
 &&1&&\\
 &1&&1&\\
 0&&1&&0\\
 &0&&0&\\
 &&1&& .
\end{tikzcd}
\]}
Moreover
\[\mathcal H^{\bullet,\bullet}(X)=\C\langle 1,\eta,\overline\eta,\eta\wedge\overline{\eta},f,\eta\wedge f,\overline{\eta}\wedge f, \eta\wedge \overline{\eta}\wedge f\rangle,
\]
where
\[
\eta\in\mathcal{H}_A^{1,0}(X),\qquad
f\in\mathcal{H}^{1,1}(X),
\]
and exactly one of the following two alternatives occurs:
\begin{enumerate}
    \item $\eta\wedge\overline\eta\in\mathcal{H}_{BC}^{1,1}(X)$ and $f\in\mathcal{H}_A^{1,1}(X)$,
    \item $\eta\wedge\overline\eta\in\mathcal{H}_{A}^{1,1}(X)$  and
    $f\in\mathcal{H}_{BC}^{1,1}(X)$.
\end{enumerate}
\end{prop}

\begin{proof}
The case $b_1(M)=1$ corresponds to
\[
h_{BC}^{1,0}(X)=0, \qquad h_A^{1,0}(X)=1.
\]
Lemma~\ref{lemma:abc-geom-1} gives $e(M)=0$. Since $b_1(M)=1$, it follows that $b_2(M)=0$. In particular, $b^+=b^-=0$, so that from \eqref{eq:b+-hbc} and duality, we obtain
\[
h^{1,1}_{BC}(X)=h^{2,1}_{BC}(X)=h^{1,2}_{BC}(X)=h_A^{1,0}(X)=h^{0,1}_A(X)=h_A^{1,1}(X)=1,
\]
Apart from bidegrees $(0,0)$ and $(2,2)$, all other Bott--Chern and Aeppli numbers vanish.

We now determine the multiplicative structure of
$\mathcal H^{\bullet,\bullet}(X)$. Fix a real nonzero generator
$
\alpha\in\mathcal H_{BC}^{1,1}(X)$.
Then,
\[
\alpha^2=0.
\]
Indeed, $\alpha^2\in\mathcal{H}^{2,2}(X)=\C\langle\omega^2\rangle$, where $\omega$ is the fundamental form of $g$. Hence, $\alpha^2=c\omega^2$. Since $b_2(M)=0$ and $d\alpha=0$, there exists
$\xi\in\mathcal{A}_{\C}^{1}(X)$ such that $\alpha=d\xi$. Thus, by Stokes' theorem, we obtain that $c=0$.
Set
$\gamma:=\ast\alpha\in\mathcal H_A^{1,1}(X)$.
Since $\alpha$ is real, so is
$\gamma$. Moreover,
\[
\gamma^2=(\ast\alpha)^2=\alpha^2=0,
\]
whereas
\[
\alpha\wedge\gamma
 =
\alpha\wedge\ast\alpha
 =
|\alpha|_g^2\,\frac{\omega^2}{2}
\neq0.
\]
We also observe that
$\mathcal H_{BC}^{1,1}(X)
 \cap
\mathcal H_A^{1,1}(X)
 =
\{0\}$. Indeed, a form in this intersection is both $d$-closed and
$d^*$-closed, hence de Rham harmonic. Since $b_2(M)=0$, it must
vanish. Consequently,
\[
\mathcal H^{1,1}(X)
 =
\C\langle\alpha\rangle\oplus\C\langle\gamma\rangle.
\]
Choose $0\neq\eta\in\mathcal H^{1,0}(X)$. By Theorem~\ref{thm:geom-BC-A}, $\eta$ is nowhere-vanishing. Since $\mathcal H^{\bullet,\bullet}(X)$ is closed under wedge product,
\[
e:=i\eta\wedge\overline\eta\in\mathcal H^{1,1}(X,g).
\]
We claim that $e$ lies in exactly one of $\mathcal H_A^{1,1}(X,g)$ and $\mathcal H_{BC}^{1,1}(X,g)$. Indeed, $e$ is real and nowhere-vanishing, so
$e=a\alpha+b\gamma$
for constants $(a,b)\in\R^2\setminus\{(0,0)\}$.
Since $e$ is decomposable, $e^2=0$. Therefore,
\[
0=e^2
 =2ab\,\alpha\wedge\gamma.
\]
As $\alpha\wedge\gamma\neq0$, we obtain $ab=0$. However, since
$e\neq0$, exactly one of
$a$ and $b$ is nonzero.
Therefore, if
$
f:=\ast e,
$
we have
\[
e\in\mathcal H_{BC}^{1,1}(X),
\qquad
f\in\mathcal H_A^{1,1}(X),
\]
or
\[
e\in\mathcal H_A^{1,1}(X),
\qquad
f\in\mathcal H_{BC}^{1,1}(X).
\]
We observe that
\[
\eta\wedge f\neq 0, \qquad \overline{\eta}\wedge f\neq 0,
\]
since, up to a nonzero constant, $\eta\wedge f=\iota_{\overline{\eta}^{\sharp}}\omega^2$, and $\overline{\eta}^{\sharp}$ is nowhere-vanishing. The argument for $\overline{\eta}\wedge f$ is analogous.

By ABC-geometric formality,
\[
\mathcal{H}^{2,1}(X)=\mathcal{H}_{BC}^{2,1}(X)=\C\langle\eta\wedge f\rangle, \qquad \mathcal{H}^{1,2}(X)=\mathcal{H}^{1,2}_{BC}(X)=\C\langle\overline{\eta}\wedge f\rangle,
\]
since
$h_A^{2,1}(X)
 =
h_A^{1,2}(X)
 =
0$.
Finally,
\[
\mathcal H^{2,2}(X)=\mathcal H_{BC}^{2,2}(X)
 =
\mathcal H_A^{2,2}(X)=\C\langle e\wedge f\rangle
 =\C\langle
i\eta\wedge\overline\eta\wedge f\rangle.
\]
\end{proof}

The differential structure of this algebra gives the following
restriction on the metric.

\begin{lemma}\label{lemma:abc-b1-one-lck}
Let $g$ be an ABC-geometrically formal Hermitian metric on a cohomologically Hopf compact
complex surface $X$. Then $g$ is locally conformally
K\"ahler and its Lee form $\theta$ has constant positive squared norm
$q=|\theta|_g^2$. Furthermore, for $0\ne\eta\in\mathcal H^{1,0}(X)$,
the alternatives of Proposition~\ref{prop:abc-non-kahler-b1-one-diamond}
give respectively
\[
\begin{aligned}
\eta\wedge\bar\eta\in\mathcal H_{BC}^{1,1}
&\quad\Longrightarrow\quad d(J\theta)=-\theta\wedge J\theta,\\
\eta\wedge\bar\eta\in\mathcal H_A^{1,1}
&\quad\Longrightarrow\quad d(J\theta)=\theta\wedge J\theta-q\omega.
\end{aligned}
\]
In the second case, $g$ is Vaisman.

\end{lemma}

\begin{proof}

Let $g$ be such an ABC-geometrically formal metric on $X$, with fundamental form $\omega$, and put
$\mathcal H^k:=\text{Tot}^k(\mathcal{H}^{\bullet,\bullet}(X))$.
By Proposition~\ref{prop:abc-non-kahler-b1-one-diamond} and its proof, a nonzero real generator
$\alpha\in\mathcal H_{BC}^{1,1}$ has constant positive norm and
satisfies $\alpha^2=0$. Thus $\alpha$ has rank one. After changing
its sign and rescaling by a constant, its unique nonzero eigenvalue
with respect to $\omega$ is one. Setting $\beta=\ast\alpha$, we obtain
\[
\omega=\alpha+\beta,\qquad
\mathcal H^2=\mathcal H^{1,1}=\C\langle\alpha,\beta\rangle,
\qquad \beta\in\mathcal H_A^{1,1}.
\]
In particular, $\omega\in\mathcal H^2$. Moreover,
$\dim_\C\mathcal H^1=\dim_\C\mathcal H^3=2$.
Since $\mathcal H^\bullet$ is closed under wedge product, the map
$L_\omega\colon\mathcal H^1\to\mathcal H^3$,
$\xi\mapsto\xi\wedge\omega$, is well defined. It is injective,
because wedge multiplication by $\omega$ is pointwise an isomorphism
$\bigwedge^1X\to\bigwedge^3X$. Thus $L_\omega$ is an isomorphism
between these two-dimensional spaces.
\smallskip

\noindent Let $\theta$ be the Lee form of $g$, so that
$d\omega=\theta\wedge\omega$. Since $g$ is ABC-geometrically formal, $d\omega\in\mathcal H^3$.
The isomorphism $L_\omega$ therefore gives $\theta\in\mathcal H^1$. Therefore $d\theta=a\alpha+b\beta$ for
constant real coefficients $a,b$. Since $\beta$ is Aeppli harmonic,
$d^*\beta=0$, while $\alpha$ and $\beta$ are pointwise orthogonal
because $\alpha\wedge\ast\beta=\alpha^2=0$. It follows that
$0=(d\theta,\beta)_{L^2}=b\|\beta\|_{L^2}^2$, so $b=0$.
Finally,
\[
0=d^2\omega=d\theta\wedge\omega=a\alpha\wedge\beta.
\]
Since $\alpha\wedge\beta$ never vanishes, $a=0$. Thus $d\theta=0$,
and $g$ is locally conformally K\"ahler.

We now determine $d(J\theta)$. The normalization above gives
$\alpha^2=\beta^2=0$ and $\alpha\wedge\beta=\omega^2/2$.
Since $\mathcal H^1=\mathcal H_A^{1,0}\oplus\mathcal H_A^{0,1}$,
we also have $J\theta\in\mathcal H^1$, and
$q:=|\theta|_g^2$ is constant by Theorem~\ref{thm:geom-BC-A}.
Moreover, $q>0$, since $X$ is non-K\"ahler.

Since $\mathcal H$ is closed under $d$, we have
$d(J\theta)=a\alpha+b\beta$ for constants $a,b\in\R$.
 Since $d^*\beta=0$ and $\alpha$ and $\beta$ are
orthogonal, the identity $0=\langle dJ\theta,\beta\rangle_{L^2}$ gives $b=0$.  Stokes' theorem and
$\theta\wedge J\theta\wedge\omega=q\omega^2/2$ now yield
\[
0=\int_X d(J\theta\wedge\omega)
 =\int_X\bigl(a\alpha\wedge\omega
       +\theta\wedge J\theta\wedge\omega\bigr)
 =(a+q)\operatorname{Vol}_g(X).
\]
Consequently, $d(J\theta)=-q\alpha$.

The decomposable semipositive form $\theta\wedge J\theta$ belongs
to $\mathcal H^2=\C\langle\alpha,\beta\rangle$, i.e., $\theta\wedge J\theta=a'\alpha+b'\beta$, for $a',b'\in \R$. Since $(\theta\wedge J\theta)^2=\alpha^2=\beta^2=0$, we obtain  $a'b'=0$. From $\text{tr}_{\omega}(\theta\wedge J\theta)=|\theta|_g^2=q$, we must have that
\[
\theta\wedge J\theta=q\alpha
\qquad\text{or}\qquad
\theta\wedge J\theta=q\beta.
\]

Since $\theta\in\mathcal H^1$ is real and nowhere-vanishing,
$\theta^{1,0}$ spans $\mathcal H_A^{1,0}$. Thus
$i\eta\wedge\bar\eta$ is a positive constant multiple of
$\theta\wedge J\theta$. The two alternatives above correspond
exactly to those of Proposition~\ref{prop:abc-non-kahler-b1-one-diamond}.
In the first case, $d(J\theta)=-q\alpha=-\theta\wedge J\theta$.

In the second case, $\theta\wedge J\theta=q\beta$, and therefore
\begin{equation*}
d(J\theta)=\theta\wedge J\theta-|\theta|_g^2\omega.
\end{equation*}
Let $\nabla$ denote the Levi--Civita connection of $g$.
For any locally conformally K\"ahler metric, one has
\[
\bigl(d(J\theta)-\theta\wedge J\theta
      +|\theta|_g^2\omega\bigr)(U,JV)
=
2(\nabla\theta)^{1,1}(U,V).
\]
Here $(\nabla\theta)^{1,1}(U,V)
=\tfrac12((\nabla\theta)(U,V)+(\nabla\theta)(JU,JV))$.
The identity follows by antisymmetrizing the formula for
$\nabla(J\theta)$ in \cite[equations~(3)--(6)]{MoroianuMoroianu2017}.
Hence the identity above implies
$(\nabla\theta)^{1,1}=0$, i.e., $g$ is {\em pluricanonical} \cite{Kok}.
By \cite[Theorem~1]{MoroianuMoroianu2017}, pluricanonical
metrics on compact complex surfaces are Vaisman. Thus $g$ is Vaisman.

\end{proof}

The Vaisman implication has a converse on every cohomologically Hopf surface.

\begin{lemma}\label{lemma:vaisman-coho-hopf-abc}
Every Vaisman metric on a cohomologically Hopf compact complex surface
is ABC-geometrically formal.
\end{lemma}

\begin{proof}
Let $g$ be any Vaisman metric on a cohomologically Hopf surface $X$.
After rescaling $g$ by a positive constant, which does not change its
harmonic spaces, assume that its Lee form $\theta$ has unit norm.
Since $b_1(X)=1$ and $b_2(X)=0$, the description in
\cite[Theorem~4.2]{IO}, together with Hodge duality, gives
\[
\mathcal{H}_{BC}(X)=\C\langle1, \omega_T,\omega_T\wedge \theta^{1,0},\omega_T\wedge\theta^{0,1}, \omega_T\wedge \theta^{1,0}\wedge \theta^{0,1}\rangle,
\]
and
\[
\mathcal{H}_{A}(X)=\C\langle1,\theta^{1,0},\theta^{0,1}, \theta^{1,0}\wedge \theta^{0,1}, \omega_T\wedge \theta^{1,0}\wedge \theta^{0,1}\rangle,
\]
where $\omega_T=-d(J\theta)=-2i\partial\theta^{0,1}$ is the
transverse K\"ahler form. Their sum is the algebra generated by
$\theta^{1,0}$, $\theta^{0,1}$, and $\omega_T$, with $\omega_T^2=0$.
Moreover,
\[
\partial\theta^{1,0}=\bar\partial\theta^{0,1}=0,\qquad
\bar\partial\theta^{1,0}=-\frac i2\omega_T,\qquad
\partial\theta^{0,1}=\frac i2\omega_T,\qquad d\omega_T=0.
\]
Hence the sum is also closed under $\partial$ and $\bar\partial$,
so $g$ is ABC-geometrically formal.

\end{proof}

\subsection{Vaisman ABC-geometrically formal surfaces}
We now prove the converse to Lemma~\ref{lemma:vaisman-coho-hopf-abc}
for cohomologically Hopf surfaces other than Inoue surfaces.

\begin{prop}\label{prop:abc-vaisman-coho-hopf}\label{prop:abc-hopf}\label{prop:abc-elliptic-vaisman}

Let $X$ be a cohomologically Hopf compact complex surface which is not
an Inoue surface. A Hermitian metric $g$ on $X$ is ABC-geometrically
formal if and only if it is Vaisman.
Moreover, $X$ admits such metrics if and only if it is either a finite
unramified quotient of a primary diagonal Hopf surface, a secondary
Kodaira surface, or a cohomologically Hopf properly elliptic surface.

\end{prop}

\begin{proof}

The Vaisman implication is Lemma~\ref{lemma:vaisman-coho-hopf-abc}.
Conversely, let $g$ be ABC-geometrically formal. By
Lemma~\ref{lemma:abc-b1-one-lck}, it suffices to exclude
$d(J\theta)=-\theta\wedge J\theta$. Suppose, for a contradiction,
that this identity holds, and use the normalization of that lemma.
The Enriques--Kodaira classification leaves Hopf surfaces,
secondary Kodaira surfaces, and cohomologically Hopf properly elliptic
surfaces. We treat the first case and the two elliptic cases separately.

Suppose first that $X$ is a Hopf surface.

On the simply
connected universal cover $\widetilde X=\C^2\setminus\{0\}$,
write the lifted Lee form as $\theta=df$ and put $u=e^f>0$.
Then
\[
d(uJ\theta)=e^f\theta\wedge J\theta+e^fdJ\theta=0,
\]
since $-dJ\theta=\theta\wedge J\theta=q\alpha$. Hence,  $uJ\theta=dv$ for a real function $v$.
The function $w=u+iv$ is holomorphic, because
$dw=u(\theta+iJ\theta)$ has type $(1,0)$, and $\operatorname{Re}w>0$.
Its Cayley transform $F(w):=(w-1)/(w+1)$ is bounded on $\C^2\setminus\{(0,0)\}$, since $|F|\leq 1$. Hence, it extends across $0\in\C^2$ by Hartogs' theorem to a bounded holomorphic function on $\C^2$, which is constant by Liouville's theorem.
Thus $u$ is constant, contradicting $du=u\theta$ and $|\theta|_g^2=q>0$.

Suppose now that $X$ is either a secondary Kodaira surface or a
cohomologically Hopf properly elliptic surface. In both cases, the elliptic
fibration $\pi:X\to B$ has smooth reduced fibers, possibly with
multiplicities \cite[Chapter~V]{BPV}, and $B$ has genus zero because
$2g(B)\leq b_1(X)=1$. Over every regular fiber $F$, Stokes' theorem gives
$\int_F\theta\wedge J\theta=0$. The form $\theta\wedge J\theta$ is
semipositive, hence $\theta|_{TF}=0$, i.e., $\theta$ is horizontal.
Moreover, for every vertical vector field $Y$, since $d\theta=0$, we have
\begin{equation}\label{eq:lee-basic}
\mathcal L_Y\theta=\iota_Yd\theta+d(\iota_Y\theta)=0,
\end{equation}
so $\theta$ is also basic. As a consequence, it descends to $B$ over
the regular values.

We now take care of the multiple fibers. Let $F=mF_0$ be a multiple
fiber, and choose a disk $\Delta\subset B$ with coordinate $z$ centered
at $\pi(F_0)$. Near each point of $F_0$, there are coordinates $(x,w)$
on $X$ such that $z(\pi(x,w))=x^m$. Let $\widehat\Delta$ be a new disk
with coordinate $t$ and define $p:\widehat\Delta\to\Delta$ by $p(t)=t^m$.
Let $\widehat X$ be the normalization of the fiber product
$$
\widehat\Delta\times_\Delta\pi^{-1}(\Delta)
=\{(t,P)\in\widehat\Delta\times\pi^{-1}(\Delta):p(t)=\pi(P)\}.
$$
Locally, this fiber product is given by $t^m=x^m$ and the normalization
separates the branches $x=\zeta^kt$, where $k=0,\ldots,m-1$ and
$\zeta=e^{2\pi i/m}$. Each normalized branch has coordinates $(t,w)$,
and the maps are
$$
\rho(t,w)=(\zeta^kt,w),\qquad \widehat\pi(t,w)=t.
$$
Thus $\rho$ is locally biholomorphic and $\widehat\pi$ is a submersion,
including along the central fiber. Moreover, $\widehat\pi$ is proper
and its noncentral fibers are connected. Since a proper submersion is
locally differentiably trivial, its central fiber is connected as well.

On this local base change, $\rho^*(\theta|_{\pi^{-1}(\Delta)})$ is smooth
and closed. It is horizontal away from the central fiber, and hence
everywhere by continuity. The same Lie derivative identity
\eqref{eq:lee-basic} shows that it is basic with respect to
$\widehat\pi:\widehat X\to\widehat\Delta$. Since the fibers are connected,
there is a smooth closed one-form $\lambda$ on $\widehat\Delta$ such that
$\rho^*\theta=\widehat\pi^*\lambda$.
The cyclic action $(t,P)\mapsto(\zeta t,P)$ on the fiber product lifts to
$\widehat X$, preserves $\rho$, and therefore preserves $\rho^*\theta$.
Consequently, $\lambda$ is invariant under $t\mapsto\zeta t$.
It represents an orbifold one-form on the quotient disk $\Delta$, with
an orbifold point of order $m$ at $\pi(F_0)$. Together with the descent
over the regular values, this shows that $\theta$ is the pullback of a
closed orbifold one-form $\widehat\theta$ on $B$.

For each multiple fiber, choose a smooth primitive of $\lambda$ on
$\widehat\Delta$. Averaging it over the cyclic group makes it invariant,
so it descends to $\Delta$ as an orbifold-smooth function. Thus $\widehat\theta$ has local primitives
also near the multiple fibers. On connected overlaps, these primitives
differ by constants. Since $B$ is a sphere and $H^1(B;\R)=0$, after
refining the cover we can modify the primitives by constants so that
they glue to a global orbifold-smooth function $f$. Its pullback $\pi^*f$ is smooth on
$X$ and satisfies $d(\pi^*f)=\theta$. Compactness of $X$ forces $\pi^*f$
to have a critical point, contradicting the fact that $\theta$ is
nowhere-vanishing.

Thus the first alternative is impossible in every case, and
Lemma~\ref{lemma:abc-b1-one-lck} implies that $g$ is Vaisman.
Finally, by \cite[Theorem~1]{Bel}, a cohomologically Hopf compact complex
surface $X$ admits a Vaisman metric if and only if it is either a finite
unramified quotient of a primary diagonal Hopf surface, a secondary
Kodaira surface, or a cohomologically Hopf properly elliptic surface.
Together with Lemma~\ref{lemma:vaisman-coho-hopf-abc}, this proves the
existence statement.

\end{proof}

\begin{rmk}
The final surface-classification assertion in \cite[Theorem~6.6]{IO}
must also include cohomologically Hopf properly elliptic surfaces.
Indeed, they admit Vaisman metrics by \cite[Theorem~1]{Bel}, and
Lemma~\ref{lemma:vaisman-coho-hopf-abc} makes these metrics
ABC-geometrically formal, hence Bott--Chern geometrically formal.

\end{rmk}

\subsection{\texorpdfstring{Inoue surfaces}{Inoue surfaces}}\label{subsub:k(x)=-infty}
It remains to treat Inoue surfaces, which admit no Vaisman metrics \cite[Theorem~1]{Bel}.

Recall that Inoue surfaces are compact complex surfaces of class VII
with $b_2=0$ and no compact complex curves. They are obtained as
compact quotients of $\mathbb H\times\mathbb C$, where $\mathbb H\subset\C$
denotes the upper half-plane, by groups of affine holomorphic
transformations acting freely and properly discontinuously
\cite{Ino}. They fall into three families, denoted by
\[
S_M,\qquad S^+_{N,p,q,r;t},\qquad S^-_{N,p,q,r}.
\]
Here $M\in\mathrm{SL}(3,\mathbb Z)$ has one real eigenvalue
$\alpha>1$ and two non-real conjugate eigenvalues $\beta,\bar\beta$.
For the $S^+$ and $S^-$ families, respectively, the matrix
$N\in\mathrm{GL}(2,\mathbb Z)$ has eigenvalues
$\alpha>1$ and $\pm\alpha^{-1}$, and hence determinant $\pm1$.
In both cases $p,q\in\mathbb Z$ and
$r\in\mathbb Z\setminus\{0\}$, while  $t\in\mathbb C$ is an additional parameter of the $S^+$ family.

\begin{prop}\label{prop:abc-inoue}
An Inoue surface admits an ABC-geometrically formal Hermitian
metric if and only if it admits a locally conformally K\"ahler
metric. Equivalently, it is of type $S_M$, $S^-_{N,p,q,r}$, or
$S^+_{N,p,q,r;t}$ with $t\in\R$.
\end{prop}
\begin{proof}
Let $X$ be an Inoue surface.
From Lemma~\ref{lemma:abc-b1-one-lck}, if $X$ admits an ABC-geometrically formal metric, it must be a locally conformally K\"ahler Inoue surface.  In particular,
$S^+_{N,p,q,r;t}$ with $t\notin\R$ are excluded by
\cite[Theorem~7]{Bel}.

Conversely, consider first $S_M$ or $S^+_{N,p,q,r;t}$ with
$t\in\R$, endowed with its standard diagonal invariant metric.
In the coframe $\{\varphi^1,\varphi^2\}$ used in
\cite[Sections~3 and~4]{AS20}, the harmonic representatives give
\begin{align*}
\mathcal H_{BC}^{\bullet,\bullet}
=\C\langle1,\varphi^{2\overline{2}},\varphi^{12\overline{1}},\varphi^{1\overline{12}},\varphi^{12\overline{12}}\rangle, \qquad
\mathcal H_A^{\bullet,\bullet}
=\C\langle1,\varphi^{2},\varphi^{\overline{2}},\varphi^{1\overline{1}},\varphi^{12\overline{12}}\rangle,
\end{align*}
where, e.g., $\phi^{ij\overline{k}}:=\phi^{i}\wedge \phi^j\wedge\overline{\phi^k}$.
Their sum $\mathcal{H}^{\bullet,\bullet}=\mathcal{H}^{\bullet,\bullet}_{BC}+\mathcal{H}_A^{\bullet,\bullet}$ is the bigraded algebra generated by $\varphi^2$, $\varphi^{\overline{2}}$,
and $\varphi^{1\overline{1}}$. The structure equations in \cite[Section~3]{AS20} yield that $\mathcal{H}^{\bullet,\bullet}$
is closed under $\partial$ and $\bar\partial$. Hence, the metric is
ABC-geometrically formal.

Every $S^-_{N,p,q,r}$ has an unramified double cover of type
$S^+_{N^2,p',q',r;0}$. The standard diagonal metric above is invariant
under the deck transformation and descends to $S^-_{N,p,q,r}$.
Harmonic forms on the quotient identify with the invariant harmonic
forms on the cover. Let $\pi:\widetilde X\to X$ denote the double cover, with deck
group $\Gamma$. Since the deck transformations are holomorphic
isometries, the pullback $\pi^*$ identifies harmonic forms on $X$ with
$\Gamma$-invariant harmonic forms on $\widetilde X$. Moreover,
\[
\pi^*\bigl(\mathcal H_{BC}(X)+\mathcal H_A(X)\bigr)
=
\bigl(\mathcal H_{BC}(\widetilde X)
      +\mathcal H_A(\widetilde X)\bigr)^\Gamma.
\]
Indeed, any $\Gamma$-invariant form in the sum on the right
admits a decomposition into $\Gamma$-invariant harmonic
summands, obtained by averaging each summand over $\Gamma$.
Since the deck transformations preserve wedge products and
commute with $\partial$ and $\bar\partial$, the invariant
part of the bidifferential algebra on $\widetilde X$ is
itself a bidifferential algebra. Thus the descended metric
is ABC-geometrically formal.
\end{proof}

\begin{rmk}\label{rmk:abc-strict-bc}
ABC-geometrically formal compact complex manifolds form a proper
subclass of geometrically Bott--Chern formal compact complex
manifolds. As recalled in Theorem \ref{thm:geom-BC-A}, an ABC-geometrically formal metric is geometrically Bott-Chern formal \cite{MS}.
The inclusion is strict already among compact complex surfaces, as Inoue surfaces $S^+_{N,p,q,r;t}$ with $t\notin\R$ admit
geometrically Bott--Chern formal metrics \cite[Proposition~3.2]{AngTom},
but admit no ABC-geometrically formal metric by
Proposition~\ref{prop:abc-inoue}.

\end{rmk}

\subsection{\texorpdfstring{The case $b_1=3$: primary Kodaira surfaces}{The case b1=3: primary Kodaira surfaces}}
We complete the picture by determining the ABC-geometrically formal non-K\"ahler surfaces with $b_1(M)=3$.

\begin{prop}\label{prop:abc-non-kahler-b1-three-diamond}
Let $X=(M,J)$ be a compact non-K\"ahler surface endowed with an
ABC-geometrically formal Hermitian metric $g$. If $b_1(M)=3$,
then the Bott--Chern diamond of $X$ must be
{\small
\[
\begin{tikzcd}
  & & 1 & &\\
  &1 & &  1&\\
  1 && 3 && 1\\
  &2 & &  2&\\
  &&1&&
\end{tikzcd}
\]}
and $k(X)=0$.
In particular, $X$ is a primary Kodaira surface and $g$ is Vaisman.
\end{prop}
\begin{proof}
From \eqref{eq:b+-hbc} we have
\[
h_{BC}^{1,0}(X)=1, \qquad h^{1,0}_{A}(X)=2.
\]
Lemma~\ref{lemma:abc-geom-1} gives $e(M)=0$, so $b_2(M)=4$. Since $b^+/2=h_{BC}^{2,0}(X)\leq1$, we have $b^+\in\{0,2\}$. If $b^+=0$, then $b^-=4$ and $h_{BC}^{1,1}(X)=5$, contradicting the bound $h_{BC}^{1,1}(X)\leq4$ in Theorem~\ref{thm:geom-BC-A}. Thus $b^+=b^-=2$, and \eqref{eq:b+-hbc} gives
\[
h^{2,0}_{BC}(X)=1, \quad h^{1,1}_{BC}(X)=3.
\]
Choose $0\neq\sigma\in\mathcal H_{BC}^{2,0}(X)$. Its constant
positive norm implies that it is a nowhere-vanishing holomorphic
section of $K_X$. Thus $K_X$ is holomorphically trivial and $k(X)=0$.
Since $b_1(M)=3$, the Enriques--Kodaira classification shows that
$X$ is a primary Kodaira surface.

We now prove that every ABC-geometrically formal metric is Vaisman.
Constant Gram--Schmidt orthogonalization and rescaling provide a global coframe
$\alpha,\beta\in\mathcal H^{1,0}$, with
$\alpha\in\mathcal H_{BC}^{1,0}$, normalized by
\[
\omega=\frac i2\bigl(\alpha\wedge\bar\alpha
                         +\beta\wedge\bar\beta\bigr),
\qquad d\alpha=0.
\]
By Lemma~\ref{lemma:n,0-inv},
\[
\Omega:=\alpha\wedge\beta
\in\mathcal H^{2,0}
=\mathcal H_{BC}^{2,0}(X)=\mathcal H_A^{2,0}(X).
\]
Thus $d\Omega=0$ and $\partial^*\Omega=0$. Writing
$\partial\beta=A\Omega$, with $A\in\C$, we obtain
\[
A\|\Omega\|_{L^2}^2
=(\partial\beta,\Omega)_{L^2}
=(\beta,\partial^*\Omega)_{L^2}=0.
\]
Hence $\partial\beta=0$. Since $\mathcal{H}^{1,1}=\C\langle\alpha\wedge\overline{\alpha},\alpha\wedge\overline{\beta},\beta\wedge\overline{\alpha},\beta\wedge\overline{\beta}\rangle$, the identity $-\alpha\wedge d\beta=d\Omega=0$ gives
\[
d\beta=B\,\alpha\wedge\bar\alpha
       +C\,\alpha\wedge\bar\beta,
\qquad B,C\in\C.
\]
Applying $d$ once more yields
$0=d^2\beta=-|C|^2\alpha\wedge\bar\alpha\wedge\beta$, so $C=0$.
Consequently,
\[
d\alpha=0,\qquad d\beta=B\,\alpha\wedge\bar\alpha,
\qquad
\theta=\bar B\beta+B\bar\beta
\]
satisfy $d\omega=\theta\wedge\omega$ and $d\theta=0$.
Here $B\neq0$, since otherwise $g$ would be K\"ahler. After rotating
$\beta$ by a constant phase, assume $B>0$. Writing
$\alpha=e^1+ie^2$ and $\beta=e^3+ie^4$ gives a real orthonormal
coframe with
\[
de^1=de^2=de^3=0,\qquad
 de^4=-2B\,e^1\wedge e^2,\qquad \theta=2B\,e^3.
\]
The dual vector field $E_3$ commutes with every frame vector and is
orthogonal to every bracket. The Koszul formula therefore gives
$\nabla E_3=0$, and hence $\nabla\theta=0$. Thus $g$ is Vaisman.

Conversely, primary Kodaira surfaces admit ABC-geometrically formal
metrics \cite{STprelim}. Indeed, for the standard diagonal invariant
metric, the sum of the Bott--Chern and Aeppli harmonic spaces is the
full double complex of invariant forms \cite[Section~4]{AS20},
which is a bidifferential bigraded algebra.

\end{proof}

\begin{rmk}
On primary Kodaira surfaces, the Vaisman condition does not imply
ABC-geometric formality. Indeed, the orthonormal coframe constructed
in the proof of Proposition~\ref{prop:abc-non-kahler-b1-three-diamond}
has constant structure coefficients. Hence every ABC-geometrically
formal metric on a primary Kodaira surface has constant scalar
curvature. On the other hand, Fino--Grantcharov--Perez~\cite{FGP24}
construct Vaisman metrics on the Kodaira--Thurston surface with
nonconstant scalar curvature. These metrics are therefore not
ABC-geometrically formal.
\end{rmk}

\subsection{\texorpdfstring{Classification}{Classification}}

To summarize, we have obtained the following.
\begin{thm}\label{thm:classification-abc-nonkahler}
Let $X=(M,J)$ be a compact complex non-K\"ahler surface. Then $X$
admits an ABC-geometrically formal Hermitian metric if and only if
it is one of the following:
\begin{enumerate}
\item a finite unramified quotient of a primary diagonal Hopf surface;
\item an Inoue surface of type $S_M$, $S^-_{N,p,q,r}$, or
$S^+_{N,p,q,r;t}$ with $t\in\R$;
\item a primary or secondary Kodaira surface;
\item a cohomologically Hopf properly elliptic surface.
\end{enumerate}
Every ABC-geometrically formal Hermitian metric on $X$ is locally
conformally K\"ahler and is Vaisman whenever $X$ is not an Inoue surface.
 Among the surfaces in this list, precisely those which
are not Inoue surfaces admit Vaisman metrics.
For $b_1(X)=1$ and $X$ not Inoue, a Hermitian metric on $X$ is
ABC-geometrically formal if and only if it is Vaisman.

\end{thm}

\begin{proof}
From Lemma \ref{lemma:h(1,0)-abc}, it suffices to study the cases  $b_1(M)=1$ and $b_1(M)=3$, which follow, respectively, from Propositions \ref{prop:abc-non-kahler-b1-one-diamond},
\ref{prop:abc-vaisman-coho-hopf}, and~\ref{prop:abc-inoue},
and Proposition~\ref{prop:abc-non-kahler-b1-three-diamond}.
The statements on metric rigidity follow from these propositions and
Lemma~\ref{lemma:abc-b1-one-lck}. The assertion about the existence
of Vaisman metrics is precisely \cite[Theorem~1]{Bel}.
\end{proof}

\section{Geometrically formal Hermitian metrics on K\"ahler surfaces}\label{sec:kahler-surfaces}\label{sec:examples}
We characterize the compact K\"ahler surfaces admitting either
geometrically Dolbeault formal or ABC-geometrically formal Hermitian
metrics. We first establish the cohomological restrictions and
minimality. The case-by-case
arguments then lead to the common classification in
Theorem~\ref{thm:classification-kahler}.

\subsection{Cohomological restrictions and minimality}

Let $X=(M,J)$ be a compact K\"ahler surface with a geometrically Dolbeault formal Hermitian metric $g$.
For the complex orientation, the Hodge and intersection numbers satisfy \cite{Draghici-li-Zhang-on-J-antiinvariant}
\[
b^+=2h^{0,2}(X)+1, \qquad b^-=h^{1,1}(X)-1.
\]
By Theorem \ref{thm:st25}, $h^{0,2}(X)\leq1$ and $h^{1,1}(X)\leq 4$. Hence,
\begin{equation*}
b^+,b^-\leq 3.
\end{equation*}
Recall that if $X$ is compact K\"ahler, Lemma \ref{lemma:geom-dolb-4} gives
\[
(h^{1,0}(X),h^{0,1}(X))\in\{(0,0),(1,1),(2,2)\},
\]
corresponding to $b_1(M)\in\{0,2,4\}$. By Lemma \ref{lemma:geom-dolb-2}, the case $h^{1,0}(X)=2$ means that $X$ is a complex torus with a flat K\"ahler metric. We therefore assume in the next two lemmas that $X$ is not a complex torus, i.e., $b_1(M)\leq 2$.

\begin{lemma}\label{lemma:geom-dolb-5}
\[
b^+=1.
\]
\begin{proof}
From $b^+=2h^{0,2}(X)+1$, we have $b^+\in \{1,3\}$.

Suppose, for a contradiction, that $b^+=3$, so that
$h^{2,0}(X)=h^{0,2}(X)=1$. By Theorem \ref{thm:st25} and Lemma \ref{lemma:n,0-inv}, the Dolbeault harmonic $(2,0)$-form $\Omega$ is a nowhere-vanishing holomorphic section of the canonical bundle $\bigwedge^{2,0}X$. Hence $c_1(X)=0$.

Combining this with Noether's formula \eqref{eq:c_1^2}, we obtain that
\begin{align*}
0=c_1^2(X)&=12\chi(X)-e(M)\\
&=12(2-h^{0,1}(X))-(2-2b_1(M)+b_2(M)).
\end{align*}
Since $X$ is not a torus, $b_1(M)\in\{0,2\}$. These two cases give $b_2(M)=22$ and $b_2(M)=14$, respectively, both contradicting Theorem~\ref{thm:st25}. Therefore $b^+=1$ and $h^{2,0}(X)=h^{0,2}(X)=0$.
\end{proof}
\end{lemma}
\begin{rmk}\label{rmk:proj}
In particular, unless $X$ is a complex torus, $X$ is projective.
Indeed, Lemma \ref{lemma:geom-dolb-5} gives $h^{2,0}(X)=0$, so
$H^2(X;\mathbb R)=H^{1,1}(X;\mathbb R)$.
Hence the K\"ahler cone contains a rational class.
The conclusion follows from the Kodaira embedding theorem
\cite[Chapter~IV]{BPV}.
\end{rmk}

We can also restrict the value of $b^-$.

\begin{lemma}\label{lemma:geom-dolb-7}
\[
b^-\in\{0,1\}.
\]
\end{lemma}
\begin{proof}
Suppose $b^-\geq2$ and choose two linearly independent real anti-self-dual
de Rham harmonic forms. By Lemma~\ref{lemma:asd-dolb}, they are
Dolbeault harmonic, so their pointwise inner products are constant by Theorem \ref{thm:st25}.
Constant Gram--Schmidt orthogonalization therefore gives two
everywhere independent real primitive $(1,1)$-forms.
Lemma~\ref{lem:primitive-complement-chern} implies $c_1(X)^2=4c_2(X)$.
On the other hand, $b^+=1$ by Lemma \ref{lemma:geom-dolb-5} and $b_1\in\{0,2\}$ give
$$
c_2(X)=3-2b_1+b^-,\qquad c_1(X)^2=9-4b_1-b^-.
$$
The equality of Chern numbers would give
$$
5b^-=4b_1-3\leq5,
$$
contradicting $b^-\geq2$.
\end{proof}

The preceding lemmas give the following restrictions.

\begin{prop}\label{prop:geom-dolb-kahler-b_1=0}
Let $X=(M,J)$ be a compact K\"ahler surface admitting a geometrically Dolbeault formal metric $g$.\vspace{0.1cm}\\
(1) If $b_1(M)=0$, then the Hodge diamond of $X$ is one of the following
{\small
\begin{equation*}
\begin{tikzcd}
  & & 1 & &\\
  & 0 &  & 0 &\\
  0 & & 1 & & 0\\
  & 0 &  & 0 & \\
  & & 1 & &
\end{tikzcd} \qquad\qquad \begin{tikzcd}
  & & 1 & &\\
  & 0 &  & 0 &\\
  0 & & 2 & & 0 \\
  & 0 &  & 0 & \\
  & & 1 & &
\end{tikzcd}
\end{equation*}}
with Dolbeault cohomology algebras, respectively, given by
$$
H_{\dbar}^{\bullet,\bullet}(X)\cong \C[t]/(t^3), \quad |t|=(1,1),
$$
and
$$
H_{\dbar}^{\bullet,\bullet}(X)\cong \C[s]/(s^2)\otimes_\C \C[t]/(t^2), \quad |s|=|t|=(1,1).
$$
Furthermore, $k(X)\in \{-\infty,2\}$.\vspace{0.1cm}\\
(2) If $b_1(M)=2$, then the Hodge diamond of $X$ must be
{\small
\[
\begin{tikzcd}
  & & 1 & &\\
  & 1 &  & 1 &\\
  0 & & 2 & & 0 \\
  & 1 &  & 1 & \\
  & & 1 & &
\end{tikzcd}
\]}
with Dolbeault cohomology algebra
\[
\textstyle H_{\dbar}^{\bullet,\bullet}(X)
\simeq
\bigwedge_{\C}\langle a,b\rangle
\otimes_{\C}\C[t]/(t^2),
\qquad
|a|=(1,0),\quad |b|=(0,1),\quad |t|=(1,1).
\]
Furthermore, $k(X)\in\{-\infty,0,1\}$.
\end{prop}

\begin{proof}
(1) We assume that $b_1(M)=h^{1,0}(X)=h^{0,1}(X)=0$ and $b_2(M)=h^{1,1}(X)\in\{1,2\}$. For either of the two Hodge diamonds above we have
\[
h^{1,0}(X)=h^{0,1}(X)=0,\qquad \chi(X)=1,\qquad e(M)=2+b_2(M),
\]
with $b_2(M)\in\{1,2\}$. Hence Noether's formula gives
\[
c_1^2(X)=12\chi(X)-e(M)=10-b_2(M)\in\{9,8\},
\]
so in particular $c_1^2(X)>0$.

Let us now consider the complex bilinear form $B(\alpha,\beta):=g(\alpha,\ast\beta)$, for $\alpha,\beta\in\mathcal{H}_{\dbar}^{1,1}(X)$. Since the antilinear Hodge operator $\ast$
preserves $\calh_{\dbar}^{1,1}(X)$, the pairing $B$ is nondegenerate with constant coefficients. 
If $h^{1,1}(X)=1$, any generator $\tau$ satisfies $\tau^2\ne0$,
giving the first algebra. If $h^{1,1}(X)=2$, choose an isotropic
basis $\sigma,\tau$ for $B$, normalized so that
$$
\sigma^2=\tau^2=0,\qquad \sigma\wedge\tau=\omega^2,
$$
and  $\mathcal{H}_{\dbar}^{\bullet,\bullet}(X)=\C\langle 1,\sigma,\tau,\sigma\wedge\tau\rangle$. Then, passing to Dolbeault cohomology classes
gives the second algebra.

It remains to exclude $k(X)=0$ and $k(X)=1$.
Suppose that $k(X)\in\{0,1\}$, and let $X\to X_{\min}$ be the blow-down map of $X$ to its minimal model. Minimal surfaces of Kodaira dimension $0$ or $1$ satisfy $c_1^2(X_{\min})=0$, hence if $X$ is minimal, we have a contradiction. However, if $X$ is obtained from $X_{\min}$ by $r$ blow-ups, then
\[
c_1^2(X)=c_1^2(X_{\min})-r=-r\le0,
\]
again contradicting $c_1^2(X)>0$. Therefore $k(X)\notin\{0,1\}$, and hence $k(X)\in\{-\infty,2\}$.\vspace{0.1cm}\\
(2) From $b_1(M)=2$, we have $h^{1,0}(X)=h^{0,1}(X)=1$. Then $e(M)=0$ by Lemma \ref{lemma:geom-dolb-3}, hence
\[
0=b_0-b_1+b_2-b_3+b_4=2-4+b_2,
\]
i.e., $b_2(M)=2$. In particular, $b^+=1$ and $b^-=1$, while
$h^{2,0}(X)=h^{0,2}(X)=0$, and $ h^{1,1}(X)=2$.

Let us choose nonzero forms
$
\alpha\in\calh_{\dbar}^{1,0}(X)$ and $
\beta\in\calh_{\dbar}^{0,1}(X),
$
and set $0\neq u=\alpha\wedge\beta$ and $\tau=\ast u$.
Both $u$ and $\tau$ are Dolbeault harmonic. Moreover, $\alpha$ and
$\beta$ have constant positive norms by Theorem~\ref{thm:st25}.
Since $u$ is decomposable, we have
$u^2=\tau^2=0$, and 
$2u\wedge\tau=|u|_g^2\,\omega^2\ne0.
$
Thus $u,\tau$ are linearly independent, and
$\alpha\wedge\tau$ and $\beta\wedge\tau$ are nonzero.
By geometric Dolbeault formality and the Hodge numbers, we have that
$$
\calh_{\dbar}^{\bullet,\bullet}(X)
=\C\langle 1,\alpha,\beta,u,\tau,
\alpha\wedge\tau,\beta\wedge\tau,u\wedge\tau\rangle.
$$
Passing to Dolbeault classes and setting
$a=[\alpha]$, $b=[\beta]$, and $t=[\tau]$, gives the asserted algebra.

It remains only to exclude $k(X)=2$. The Hodge diamond gives
$\chi(X)=0.$ If $k(X)=2$, let $X_{\min}$ be the minimal model of $X$. The holomorphic Euler characteristic is birationally invariant, so $\chi(X_{\min})=0$. Since $X_{\min}$ is a minimal surface of general type, both $c_1^2(X_{\min})>0$ and $c_2(X_{\min})>0$. Noether's formula then gives
\[
0=12\chi(X_{\min})=c_1^2(X_{\min})+c_2(X_{\min})>0,
\]
a contradiction. Thus $k(X)\neq2$.
\end{proof}
We immediately obtain the following corollary.
\begin{cor}\label{cor:dolb-cohom-ring-kahler-1}
    Let $X$ be a compact K\"ahler surface. Then, if $X$ admits a geometrically Dolbeault formal metric, its Dolbeault cohomology ring must be that of
    \[
    (1)\,\,\mathbb{T}^2_\C, \qquad (2)\,\, \mathbb{T}^1_{\C}\times\mathbb{CP}^1, \qquad (3)\,\, \mathbb{CP}^1\times \mathbb{CP}^1, \qquad (4)\,\, \mathbb{CP}^2.
    \]
    Moreover, $X$ is minimal.
\begin{proof}
The assertion about the Dolbeault cohomology algebra follows from
Proposition~\ref{prop:geom-dolb-kahler-b_1=0}. We prove minimality in the
four cases listed above.\vspace{0.1cm}

\noindent (1) If the algebra is that of $\mathbb T^2_\C$, then $X$ is a
complex torus by Lemma~\ref{lemma:geom-dolb-2}, and hence is minimal.\vspace{0.1cm}

\noindent (2) If the algebra is that of
$\mathbb{CP}^1\times\mathbb T^1_\C$, then $e(X)=0$ and $h^{1,0}(X)=1$.
If $X$ were obtained from its minimal model $X_{\min}$ by $r\geq1$
point blow-ups, then $e(X_{\min})=-r<0$.
The classification of compact K\"ahler surfaces implies that
$X_{\min}$ is ruled over a curve of genus at least two
\cite[Chapter VI]{BPV}. This contradicts the birational invariance of
$h^{1,0}$.\vspace{0.1cm}

\noindent (3) Suppose the algebra is that of
$\mathbb{CP}^1\times\mathbb{CP}^1$, so that $b_1(X)=0$ and $b_2(X)=2$.
If $X$ is not minimal, positivity of $b_2(X_{\min})$ implies that exactly
one point is blown up and $b_2(X_{\min})=1$.
The classification of compact K\"ahler surfaces then gives
$X_{\min}\cong\mathbb{CP}^2$ or a fake projective plane, i.e., a compact quotient of the unit ball by a torsion-free cocompact discrete subgroup of $PU(2,1)$
with $e(X_{\min})=3$.
The first possibility would give
$X\cong Bl_p\mathbb{CP}^2\cong\Sigma^1$, i.e., the first Hirzebruch surface, which is excluded by
Proposition~\ref{prop:no-hirzebruch} below.

Suppose instead that $X=Bl_pY$, where $Y$ is a fake projective plane.
Let $g$ be a geometrically Dolbeault formal Hermitian metric on $X=(M,J)$.
Since $b^-(X)=1$, choose a nonzero real anti-self-dual $d$-harmonic
$2$-form $\beta$, which is also a Dolbeault harmonic $(1,1)$-form by Lemma \ref{lemma:asd-dolb}. Geometric Dolbeault formality of $g$ implies that $\beta$ has constant positive
pointwise norm and hence it is a symplectic form on the underlying manifold
$\overline M$ with the opposite orientation.
The exceptional curve is a smoothly embedded sphere of
self-intersection $+1$ in $\overline M$.
By \cite[Corollary 2]{Li99}, a closed symplectic four-manifold containing
such a sphere is rational or ruled. Since $b_1(\overline M)=0$,
$\overline M$ must be rational, and in particular simply connected.
However,
$\pi_1(\overline M)\cong\pi_1(X)\cong\pi_1(Y)$
is infinite, since a fake projective plane is a compact quotient of the
complex two-ball \cite{PY}. This contradiction proves minimality.\vspace{0.1cm}

\noindent (4) If the Dolbeault algebra of $X$ is that of $\mathbb{CP}^2$, then $b_2(X)=1$.
Any point blow-down would give a compact K\"ahler surface with $b_2(X)=0$,
which is impossible. Thus $X$ is minimal.
\end{proof}
\end{cor}

\Needspace{14\baselineskip}
\subsection{Surfaces with \(b_1=0\)}\label{sec:dolb-b1-0}

\subsubsection{$k(X)=-\infty$}\label{sec:dolb-hirzebruch}
Among rational surfaces, $\mathbb{CP}^1\times\mathbb{CP}^1$ and $\mathbb{CP}^2$ admit geometrically formal K\"ahler metrics. The remaining surfaces with $b_1=0$, $h^{1,1}=2$, and $k(X)=-\infty$ are Hirzebruch surfaces.
Recall that, for $n\geq 0$, the $n$-th Hirzebruch surface is the ruled surface
\[
\Sigma^n:=\mathbb{P}(H^n\oplus \mathbb{C}_{\mathbb{CP}^1})\longrightarrow \mathbb{CP}^1,
\]
where $\C_{\mathbb{CP}^1}$ is the trivial bundle over $\mathbb{CP}^1$ and $H$ denotes the hyperplane line bundle on $\mathbb{CP}^1$. One has $\Sigma^0=\mathbb{CP}^1\times \mathbb{CP}^1$, and $\Sigma^1\cong Bl_1\mathbb{CP}^2$, so that every $\Sigma^n$ is a minimal rational surface for $n\neq 1$. Denote by $D_0,D_\infty\subset\Sigma^n$ the two distinguished constant sections of $\Sigma^n$, labeled so that $D_0^2=-n$ and $D_\infty^2=n$, and by $F$ the fiber class. For $n\geq 1$,
\begin{align*}
D_0\cdot F=D_\infty \cdot F=1, \qquad D_0\cdot D_\infty=0, \qquad D_\infty^2=n, \qquad D_0^2=-n,
\end{align*}
and $D_0,D_\infty$ span $H^{1,1}(\Sigma^n;\R)$.

We first recall how to integrate Dolbeault classes over complex curves.

\begin{lemma}\label{lem:restriction-lemma}
Let $X$ be a compact complex surface satisfying the
$\partial\bar\partial$-lemma. For $\bar\partial$-closed $(1,1)$-forms
$\alpha,\gamma$ and a smooth compact complex curve $C\subset X$,
\[
\int_C\alpha=[\alpha]\cdot[C],
\qquad
\int_X\alpha\wedge\gamma=[\alpha]\cdot[\gamma],
\]
where $[\alpha]$ and $[\gamma]$ are identified with their corresponding
de Rham classes of type $(1,1)$.
\end{lemma}
\begin{proof}
By the $\del\dbar$-lemma, we can write $\alpha=\rho+\bar\partial\xi$, with $\rho$ a $d$-closed
$(1,1)$-form and $\xi\in\mathcal A^{1,0}(X)$.
If $i\colon C\hookrightarrow X$ denotes the inclusion, the pullback $i^*\partial\xi$ vanishes on $C$ for bidegree reasons, so that $i^*(\dbar\xi)=d(i^*\xi)$. On $X$, the same bidegree argument and $\bar\partial\gamma=0$ give $\bar\partial\xi\wedge\gamma=d(\xi\wedge\gamma)$.
Stokes' theorem, therefore, allows us to replace $\alpha$, and
similarly $\gamma$, by $d$-closed representatives, in the integration of $\alpha$ over $C$ and $\alpha\wedge \gamma$ over $X$.
\end{proof}

\begin{prop}\label{prop:no-hirzebruch}\label{prop:no-hirzebruch-hermitian}
The Hirzebruch surface $\Sigma^n$ admits a geometrically Dolbeault formal
Hermitian metric if and only if $n=0$.
\end{prop}
\begin{proof}
For $n=0$, the product Fubini--Study metric on
$\Sigma^0=\mathbb{CP}^1\times\mathbb{CP}^1$ is geometrically Dolbeault formal.
Suppose now that $n\geq1$ and that $g$ is a geometrically Dolbeault formal
Hermitian metric on $\Sigma^n$. Since $\Sigma^n$ is K\"ahler, the $\partial\bar\partial$-lemma identifies its Dolbeault cohomology algebra with the pure-type decomposition of its complex de Rham cohomology algebra. Set
\[
a=2[D_0]+n[F].
\]
Since $F^2=0$ and $[D_\infty]=[D_0]+n[F]$, we have
\begin{equation}\label{eq:intersect-a-hirz}
a^2=0,\qquad a\cdot[D_0]=-n,\qquad
 a\cdot[D_\infty]=n,\qquad a\cdot[F]=2.
\end{equation}
Let $\sigma$ be the Dolbeault-harmonic representative of $a$.
Geometric Dolbeault formality implies that $\sigma\wedge\sigma$ is
harmonic. Since $a^2=0$, it follows that $\sigma\wedge\sigma=0$.

Let $V=w\partial_w$ be the holomorphic vector field generating the standard fiberwise
$\mathbb C^*$-action, with $w$ the fiber coordinate
so that $D_0=\{w=0\}$ and $D_\infty=\{w=\infty\}$.
This vector field is globally defined and vanishes along both sections.
Since $V$ is holomorphic,
$\bar\partial\iota_V+\iota_V\bar\partial=0$ and,
thus, $\iota_V\sigma$ is a $\bar\partial$-closed $(0,1)$-form.
As $H^{0,1}_{\bar\partial}(\Sigma^n)=0$, there exists a smooth
complex-valued function $f$ such that
$
\bar\partial f=\iota_V\sigma.
$
Consequently,
\[
\bar\partial(f\sigma)
=(\iota_V\sigma)\wedge\sigma
=\tfrac12\iota_V(\sigma\wedge\sigma)=0.
\]
Moreover, $(f\sigma)\wedge\sigma=0$, and hence
$[f\sigma]\cdot a=0$.
The intersection pairing on $H^{1,1}_{\bar\partial}(\Sigma^n)$ is
nondegenerate and this space has complex dimension two. Since $a$ is a
nonzero isotropic class, $a^\perp=\mathbb Ca$. Therefore
\begin{equation}\label{eq:fsigma-ca}
[f\sigma]=c\,a
\end{equation}
for some $c\in\mathbb C$.
The restriction of $f$ to each of $D_0$ and $D_\infty$ is holomorphic,
since $V$ vanishes there, and is therefore constant. Denote $f_0:=f|_{D_0}$
and $f_\infty:=f|_{D_{\infty}}$. By Lemma~\ref{lem:restriction-lemma}, \eqref{eq:intersect-a-hirz}, and \eqref{eq:fsigma-ca}, we have
\[
-nf_0=\int_{D_0}f\sigma=-nc,\qquad
nf_\infty=\int_{D_\infty}f\sigma=nc.
\]
Since $n>0$, we obtain $f_0=f_\infty=c$.

Finally, we restrict to a fiber $F\cong\mathbb{CP}^1$ of $\Sigma^n$. On
$F\setminus\{0,\infty\}$ write
$\sigma|_F=s(w,\bar w)\,dw\wedge d\bar w$. The potential equation $\iota_V\sigma=\dbar f$ gives
$\bar\partial(f|_F)=ws(w,\bar w)\,d\bar w$, and hence
\[
d\left(f\frac{dw}{w}\right)=-\sigma|_F.
\]
We then consider the annulus $A_{\varepsilon,R}:=\{\varepsilon\leq|w|\leq R\}\subset F$ and integrate $\sigma$ over $A_{\varepsilon,R}$. By Stokes' theorem,
\[
\int_{A_{\varepsilon,R}} \sigma=\int_{|w|=\varepsilon}f\frac{dw}{w}-\int_{|w|=R}f\frac{dw}{w}.
\]
By letting
$\varepsilon\to0$ and $R\to\infty$, smoothness of $f$ at the two endpoints
yields
\[
\int_F\sigma=2\pi i(f_0-f_\infty)=0.
\]
This contradicts $\int_F\sigma=a\cdot[F]=2$. Hence, $\Sigma^n$ does not admit any geometrically Dolbeault formal metric.
\end{proof}

\subsubsection{$k(X)=2$}\label{subsection:k=2b1=0}
Let $X$ be a surface of general type. By Corollary~\ref{cor:dolb-cohom-ring-kahler-1}, $X$ is minimal. Noether's formula
$
c_1^2(X)=10-b_2(M)
$
gives
\[
c_1^2(X)=\begin{cases}9, \qquad \text{if}\,\,\,h^{1,1}(X)=1,\\
8, \qquad \text{if}\,\,\,h^{1,1}(X)=2.
\end{cases}
\] We discuss the two cases separately.

On the one hand, the case $c_1^2(X)=9$, $h^{1,0}(X)=h^{0,1}(X)=0$, and $b_2(M)=h^{1,1}(X)=1$ corresponds to the \emph{fake projective planes}, i.e., quotients of the (open) unit ball $B\cong PU(2,1)/U(2)\subset \C^2$  by a torsion-free cocompact discrete subgroup of $PU(2,1)$. The classification was initiated by Prasad--Yeung and completed with the work of Cartwright--Steger \cite{PY,CS}. Since $b_2(M)=1$ and $h^{2,0}(X)=h^{0,2}(X)=0$, for any K\"ahler metric $g$, the space of harmonic forms is $\mathcal{H}(X)=\C\langle1,\omega,\omega^2\rangle$, which is clearly an algebra with respect to $\wedge$. Hence,
\begin{prop}\label{prop:fake-pr-dolb}
On a fake projective plane, every K\"ahler metric is geometrically formal.
\end{prop}

On the other hand, a minimal surface of general type with $b_1(X)=0$ and
$b_2(X)=2$ is called a \emph{fake quadric}: it has the same Betti numbers as
$\mathbb{P}^1\times\mathbb{P}^1$, and equivalently
\[
h^{1,0}(X)=0,\qquad c_2(X)=4,\qquad c_1^2(X)=8.
\]
We now characterize all fake quadrics admitting either geometrically Dolbeault formal or ABC-geometrically formal metrics.

We first make explicit the Seiberg--Witten chamber used in Lemma \ref{lemma:positive-canonical-chamber}. To this aim, we recall that on a closed symplectic four-manifold $(M,\Omega)$ endowed with any compatible almost complex structure $J$, the {\em symplectic canonical class} is defined as $\kappa:=-c_1(TM,J)$.

\begin{lemma}\label{lemma:positive-canonical-chamber}
Let \((M,\Omega)\) be a closed symplectic four-manifold with
\(b_1(M)=0\), \(b^+(M)=1\), and symplectic canonical class \(\kappa\).
If \(\kappa^2>0\) and \(\kappa\cdot[\Omega]>0\), the canonical
\(\operatorname{Spin}^c\) structure has nonzero Seiberg--Witten
invariant in the small-perturbation chamber for every Riemannian metric.
In particular, \(M\) admits no metric of positive scalar curvature.
\end{lemma}
\begin{proof}
We fix the component $\mathcal{C}^+$ of the positive cone containing $[\Omega]$.
Since $\kappa^2>0$ and $\kappa \cdot [\Omega]>0$, we have $\kappa\in\mathcal{C}^+$. Let $h$ be any Riemannian metric on $M$ and denote by $\alpha_h$ the unique real self-dual harmonic form satisfying the normalizations $\int_M\alpha_h\wedge\alpha_h=1$ and $\int_M\alpha_h\wedge \Omega>0$.  Set $u_h:=[\alpha_h]$.

Since $b^+(M)=1$, for every $v\in H^2(M;\R)$ such that $v^2>0$,  the intersection form restricted to $v^{\perp}$ is negative definite. If $\kappa\cdot v= 0$, then $\kappa\in v^{\perp}$ and $\kappa^2<0$, which is a contradiction with our assumptions. Hence, $\kappa\cdot v\neq 0$, for every $v\in H^2(M;\R)$ with $v^2>0$. In particular, the continuous function $v\mapsto \kappa\cdot v$ from $\mathcal{C}^+$ to $\R$ is never vanishing. Since $\kappa\cdot [\Omega]>0$ by assumption and $\mathcal{C}^+$ is connected, also $\kappa\cdot v>0$ for every $v\in\mathcal{C}^+$. In particular, $\kappa\cdot u_h>0$.

In the
conventions of \cite[Section~2.1]{LL01}, the canonical
$\operatorname{Spin}^c$ structure $\mathcal{L}$ of $(M,J)$ has determinant class $-\kappa$. Indeed, it is equal to $c_1(\det W^+)$, where $W^+=\bigwedge^{0,0}M\oplus \bigwedge^{0,2}M\cong \C\oplus K^{-1}$.
The discriminant of the $\text{Spin}^c$ structure $\mathcal{L}$ with respect to $h$ and the zero perturbation $\mu\equiv 0$ is, by definition,
\begin{equation}\label{eq:discriminant}
\Delta(h,0)=c_1(\det W^+)\cdot u_h=-\kappa\cdot u_h<0.
\end{equation}
Thus every zero-perturbation pair $(h,0)$ lies in the negative
chamber. With the homology orientation induced by $\Omega$, from Taubes's nonvanishing theorem \cite[Section~2.1]{LL01}, we have $SW_-(M,\mathcal{L})=1$.
From \eqref{eq:discriminant}, all sufficiently small perturbations remain in the negative chamber and admit no reducible solutions \cite[§3.3]{HT}. If $h$ had positive scalar curvature, the Seiberg--Witten Weitzenböck estimate would also exclude irreducible solutions for sufficiently small perturbations \cite[Lemma 3.4]{HT}. Thus, for a sufficiently small generic perturbation, the moduli space would be empty, contradicting the nonvanishing of the invariant in this chamber.

\end{proof}

\begin{prop}\label{prop:fake-quadric-uniformization}
For a fake quadric $X=(M,J)$, the following conditions are equivalent:
\begin{enumerate}
\item $X$ admits a geometrically Dolbeault formal Hermitian metric,
\item the oppositely oriented smooth manifold $\overline{M}$ is symplectic,
\item $X\cong\Gamma\backslash(\mathbb D\times\mathbb D)$, where
$\Gamma\subset\operatorname{Aut}(\mathbb D)\times
\operatorname{Aut}(\mathbb D)$ is a torsion-free cocompact lattice,
\item $X$ admits a geometrically formal K\"ahler metric;
\item $X$ admits an ABC-geometrically formal Hermitian metric.
\end{enumerate}
\end{prop}
\begin{proof}
$(1)\Rightarrow(2)$ Let $g$ be a geometrically Dolbeault formal metric on $X$. Since $b^-=1$, choose $0\ne\beta\in\mathcal H_g^-$.
By Lemma~\ref{lemma:asd-dolb}, this form is Dolbeault harmonic,
so geometric Dolbeault formality gives it constant positive norm.
It is closed and satisfies $\beta^2=-|\beta|_g^2dV_g$, hence
defines a symplectic form on $\overline M$.\smallskip

\noindent $(2)\Rightarrow (3)$ We proceed by steps.\\
\noindent \textit{Step 1.} We prove that $M$ does not admit metrics of positive scalar curvature.

Let $\Omega_-$ be a symplectic form on
$\overline M$. First, we observe that $X$ admits no positive-scalar-curvature
metric. Indeed, choose a background K\"ahler metric with fundamental form $\omega$ on $X$. Since $X$ is a minimal surface of general type, its canonical class $\kappa=-c_1(X)$ is nef, i.e., $\kappa$ lies in the closure of the K\"ahler cone, and $\kappa\cdot [\omega]\geq 0$. If $\kappa\cdot [\omega]=0$, then $\kappa\in [\omega]^{\perp}$, so $\kappa^2<0$ by the Hodge index theorem. This contradicts $\kappa^2>0$, hence $\kappa\cdot[\omega]>0$. Then, Lemma \ref{lemma:positive-canonical-chamber} applied to $(X,\omega)$ yields that $X$ does not admit any metric of positive scalar curvature. The existence of a positive-scalar-curvature metric is independent of orientation, so the same obstruction holds on \(\overline M\).\vspace{0.1cm}

\noindent \textit{Step 2.} $(\overline{M},\Omega_-)$ admits a nonzero Seiberg--Witten invariant in the small-perturbation chamber.

Let $\kappa_-$ be the symplectic canonical class of
$(\overline M,\Omega_-)$. With the intersection form of
$\overline M$,
\[
\kappa_-^2=2e(\overline M)+3\sigma(\overline M)=8.
\]
Now, on $\overline{M}$, we have $[\Omega_-]^2>0$ and
$b^+(\overline{M})=b^-(M)=1$, so if
$\kappa_-\cdot [\Omega_-]=0$, then $8=\kappa_-^2<0$,
a contradiction. Hence, $\kappa_-\cdot[\Omega_-]\ne0$.
If $\kappa_-\cdot [\Omega_-]<0$, then by
\cite[Theorem B]{Liu96}, $\overline{M}$ would be either a
rational or ruled surface. Since $b_1=0$ and $b_2=2$, its
diffeomorphism type would be either $S^2\times S^2$ or
$\mathbb{CP}^2\#\overline{\mathbb{CP}}{}^2$, both of which
admit metrics of positive scalar curvature. This contradicts
Step 1. Thus $\kappa_-\cdot [\Omega_-]>0$, and
Lemma~\ref{lemma:positive-canonical-chamber} supplies a nonzero
small-perturbation Seiberg--Witten invariant in the opposite
orientation.\vspace{0.1cm}

\noindent {\it Step 3.} We prove that
$X\cong \Gamma\backslash (\mathbb{D}\times \mathbb{D})$,
with $\Gamma\subset\text{Aut}(\mathbb{D})\times
\text{Aut}(\mathbb{D})$.

Since $X$ is a surface of general type and the Seiberg--Witten
invariant in the opposite orientation is nonzero,
\cite[Theorem~1]{Kot97} implies that $K_X$ is ample and,
since $\sigma(X)=0$, that $X$ is uniformized by the product
of two unit discs $\mathbb{D}\times\mathbb{D}$.
Furthermore, from the equality case in the proof of
\cite[Theorem~1]{Kot97}, the K\"ahler--Einstein metric $g$
induces a globally defined nonzero $g$-parallel anti-self-dual
two-form $\eta$, rescaled so that $|\eta|_g^2=2$.
Let $\omega$ be the fundamental form of $g$.
Setting $\omega_1:=1/2(\omega+\eta)$ and
$\omega_2:=1/2(\omega-\eta)$ defines two global
nowhere-vanishing forms. The distributions
$E_1:=\ker\omega_2$ and $E_2:=\ker\omega_1$ are $J$-invariant
and parallel and give the splitting $TX=E_1\oplus E_2$.
Let $\pi\colon \tilde{X}\to X$ be the universal cover of $X$.
By the de Rham decomposition,
$\tilde{X}\cong \mathbb{D}\times \mathbb{D}$ isometrically,
with $E_1$ tangent to the first factor and $E_2$ tangent to
the second. If $\gamma$ is a deck transformation of
$\tilde{X}$, then
$\gamma^*\pi^*\omega_j=\pi^*\omega_j$, i.e.,
$\gamma$ preserves $\pi^*E_1$, respectively, $\pi^*E_2$,
and hence each factor $\mathbb{D}$, i.e., the deck
transformation group $\Gamma$ is a subgroup of
$\text{Aut}(\mathbb{D})\times \text{Aut}(\mathbb{D})$.
This proves (3).\smallskip

\noindent $(3)\Rightarrow(4)$ The product hyperbolic metric descends to a K\"ahler metric $g$ on
$X$, and the two volume forms on each factor descend to parallel forms
$\Omega_1,\Omega_2$. Since $b_1=0$ and $b_2=2$, all
harmonic forms with respect to $g$ are generated by
$1,\Omega_1,\Omega_2,\Omega_1\wedge\Omega_2$, and their products are
parallel. \smallskip

\noindent $(4)\Rightarrow(5)$ For a K\"ahler metric, the three formality conditions agree, giving $(5)$ and $(1)$.\smallskip

\noindent $(5)\Rightarrow(1)$
It follows from Lemma~\ref{lemma:abc-harmonic-comparison}, since $X$
satisfies the $\partial\bar\partial$-lemma.
\end{proof}

Recall that a surface is \emph{isogenous to a product of curves} if it
can be written as $(C_1\times C_2)/G$, where the $C_i$ are compact
curves of genus at least two and $G$ is a finite group acting freely
by holomorphic automorphisms. The action is \emph{unmixed} if it
preserves both factors, and \emph{mixed} otherwise
\cite[Section~2]{FL}. Examples of unmixed fake quadrics isogenous to a product of curves are provided by {\em unmixed Beauville surfaces}
with $\chi(X)=1$
\cite{Cat2}.
\begin{cor}\label{prop:fake-quadric-mixed-unmixed}
A fake quadric isogenous to a product admits a geometrically
Dolbeault formal Hermitian metric if and only if it is of unmixed
type.
\end{cor}
\begin{proof}
Write $X=(C_1\times C_2)/G$, and let $\Gamma$ be the deck group
of the universal covering of $X$ given by the composition
$\mathbb D^2 \longrightarrow C_1\times C_2 \longrightarrow X$.
The elements of $\Gamma$ are precisely the lifts of elements of
$G$. Since the first covering is a product covering, a lift
preserves the disk factors if and only if the corresponding
element of $G$ preserves the curve factors. Thus the action is
unmixed if and only if
$\Gamma\subset H:=
\operatorname{Aut}(\mathbb D)\times\operatorname{Aut}(\mathbb D)$.

Any other biholomorphic identification of the universal cover
with $\mathbb D^2$ replaces $\Gamma$ by $F\Gamma F^{-1}$ for some
$F\in\operatorname{Aut}(\mathbb D^2)$. But the subgroup $H$ is normal in $\text{Aut}(\mathbb{D}^2)$,
being the kernel of the homomorphism
$\rho\colon \operatorname{Aut}(\mathbb D^2)\to\Z_2$ given by $\rho(g)=1$ if $g$ interchanges factors, and $\rho(g)=0$, otherwise. Consequently, $F\Gamma F^{-1}\subset H\iff \Gamma\subset H$.
Hence condition~(3) of
Proposition~\ref{prop:fake-quadric-uniformization} holds exactly in the
unmixed case.
\end{proof}

\subsubsection{The metric geometry of fake quadrics}\label{sec:fake-quadrics-isogenous}
We now characterize geometric formality for a fixed K\"ahler
metric on a fake quadric.

\begin{prop}\label{prop:fake-quadric-local-product}
Let $X$ be a fake quadric and $g$ a K\"ahler metric on $X$.
Then $g$ is geometrically formal if and only if it admits a
$g$-orthogonal almost-K\"ahler structure $I$ of the opposite
orientation. In this case, $(\overline{M},I,g)$ is a K\"ahler surface  and $(X,g)$ is locally isometric to a
product of two Riemann surfaces.
\end{prop}
\begin{proof}
If $g$ is geometrically formal, we consider $0\neq \beta\in \mathcal{H}_g^-$ with constant normalized norm
$|\beta|^2=2$. Then $\beta=g(I\cdot,\cdot)$ defines the required
opposite almost-K\"ahler structure.

Conversely, suppose that $I$ is a $g$-orthogonal almost-K\"ahler
structure of the opposite orientation. Its fundamental form
$\beta:=g(I\cdot,\cdot)$ is closed and anti-self-dual with respect
to the complex orientation of $X$, hence harmonic.
Let $\omega:=g(J\cdot,\cdot)$ be the fundamental form of the original
K\"ahler metric on $(M,J)$. Since $b_1=0$ and $b^+=b^-=1$, all
real harmonic forms with respect to the K\"ahler metric $g$ on $X=(M,J)$ are spanned by
$\{1, \omega, \beta, \omega^2/2\}$.  The identities
$\omega\wedge\beta=0$ and $\beta^2=-\omega^2$ show that this space is closed under wedge products.
Thus $g$ is geometrically formal.

It remains to prove $\nabla^g I=0$. By
Proposition~\ref{prop:fake-quadric-uniformization}, we have $X=\Gamma\backslash (\mathbb{D}\times \mathbb{D})$, with $\Gamma\subset\text{Aut}(\mathbb{D})\times \text{Aut}(\mathbb{D}) $.
Let $\Omega_1,\Omega_2$ be the forms on $X$ descended from
the hyperbolic volume forms of the disk factors. Let  $\omega=g(J\cdot,\cdot)$ and  $\beta=g(I\cdot,\cdot)$ with norm $|\beta|_g^2=2$. The forms
$\alpha:=(\omega+\beta)/2$ and $\gamma:=(\omega-\beta)/2$
are closed semipositive $(1,1)$-forms of complex rank one,
satisfying
$\alpha^2=\gamma^2=0$ and $
\alpha\wedge\gamma=\frac{\omega^2}{2}>0$.

The classes $[\Omega_1]$ and $[\Omega_2]$ form a basis of
$H^{1,1}(X;\mathbb R)$, with
$[\Omega_1]^2=[\Omega_2]^2=0$ and $
[\Omega_1]\cdot[\Omega_2]>0$.
Thus they span the two isotropic lines of the intersection
form. Since $[\alpha]\cdot[\gamma]>0$, the classes
$[\alpha],[\gamma]$ are nonzero and lie on distinct isotropic
lines. Their positive pairing with the K\"ahler class
$[\Omega_1+\Omega_2]$ therefore gives, after relabeling,
\[
[\alpha]=r[\Omega_1],
\qquad
[\gamma]=s[\Omega_2],
\qquad r,s>0.
\]

It follows that $\int_X\alpha\wedge\Omega_1=0$.
Since $\alpha\wedge\Omega_1$ is pointwise nonnegative,
it vanishes identically. Two nonzero semipositive
$(1,1)$-forms of complex rank one whose wedge product
vanishes are positively proportional, hence $\alpha=u\Omega_1$, with $u\in\mathcal{C}^{\infty}(X;\R)$ positive. The same argument yields $ \gamma=v\Omega_2$, with $v\in\mathcal{C}^{\infty}(X;\R)$ positive.

On the universal cover $\pi\colon \mathbb D^2\rightarrow X$, the equations
$d\alpha=d\gamma=0$ give
$du\wedge\Omega_1=0$ and $ dv\wedge\Omega_2=0.$
Thus $u$ depends only on the first disk coordinate and $v$
only on the second. Consequently, the lifted metric $\pi^*g$ is the product
$\tilde g
=u(z_1)g_{\mathrm{hyp},1}
+v(z_2)g_{\mathrm{hyp},2}.$
Since the bidisk identification is holomorphic, this gives
a local holomorphic and isometric product decomposition of $(X,g)$.
Moreover, the lifted form $\beta=u\Omega_1-v\Omega_2$ is the
difference of the area forms of the two factor metrics, hence
is parallel. Since $\beta=g(I\cdot,\cdot)$, we obtain
$\nabla^g I=0$. In particular, $(\overline M,I,g)$ is a
K\"ahler surface.

\end{proof}
Recall that a lattice
$\Gamma\subset\operatorname{Aut}(\mathbb D)\times
\operatorname{Aut}(\mathbb D)$ is called \emph{irreducible}
if its projection onto each factor is dense in
$\operatorname{Aut}(\mathbb D)$ with respect to the usual
Lie group topology. We call the corresponding bidisk quotient
\emph{irreducible} when its defining lattice is irreducible.
\begin{cor}\label{prop:fake-gauduchon-rigidity}
\label{cor:fake-conformal-uniqueness}
\label{cor:fake-irreducible-kahler-rigidity}
Let $X$ be a fake quadric. If $X$ is an irreducible factor-preserving bidisk quotient, every
geometrically formal K\"ahler metric lifts to a product of constant multiples of the hyperbolic metrics.
\end{cor}
\begin{proof}
Proposition~\ref{prop:fake-quadric-local-product} gives
$
\widetilde g=u(z_1)g_{\mathrm{hyp},1}+v(z_2)g_{\mathrm{hyp},2}.
$
For an irreducible lattice, both projections are dense in
$\operatorname{Aut}(\mathbb D)$. Deck invariance therefore makes
$u$ and $v$ invariant under dense subgroups of the isometry groups
of the respective disks. Continuity and transitivity imply that
both functions are constant.
\end{proof}

\subsection{Surfaces with \(b_1=2\)}\label{sec:dolb-b1-2}
In this case \(b_2(X)=h^{1,1}(X)=2\). We treat the three possible
Kodaira dimensions separately.

\subsubsection{$k(X)=-\infty$}
By minimality and the Enriques--Kodaira classification, these surfaces are geometrically ruled over an elliptic curve. We call such a $\mathbb{CP}^1$-bundle
\emph{projectively unitary} if it can be written as
\[
X\cong (\mathbb{C}\times\mathbb{CP}^1)/\Lambda,
\qquad
\lambda\cdot(z,p)=(z+\lambda,\rho(\lambda)p),
\]
where $\Lambda\subset\mathbb{C}$ is a lattice and
$\rho\colon\Lambda\rightarrow PU(2)$ is a group homomorphism.

We set out to prove that if an elliptic ruled surface is geometrically Dolbeault formal, then it must be projectively unitary.
We start with a useful normalization lemma. In the ruled case, $\theta$ below
is the pullback of a holomorphic one-form on the elliptic base.

\begin{lemma}\label{lem:ruled-harmonic-normalization}
Let $X$ be a compact K\"ahler surface with $h^{1,0}(X)=1$ and let
$g$ be a geometrically Dolbeault formal Hermitian metric. There are forms
\[
0\ne\theta\in H^{1,0}(X),\qquad
\eta\in\mathcal H^{0,1}_{\bar\partial}(X)
\]
such that
\[
|\theta|_g=|\eta|_g=1,\qquad
[\eta]_{\bar\partial}=c[\bar\theta]_{\bar\partial},\qquad
a:=\langle\eta,\bar\theta\rangle_g\in(0,1],
\]
where
\begin{align}\label{eq:proj-lemma-norm-c}
c=\frac{\operatorname{Vol}_g(X)}{\int_Xa\,dV_g}\geq1.
\end{align}
\end{lemma}
\begin{proof}
From Theorem \ref{thm:st25}, we can normalize $\theta$ and $\eta$ so that $|\theta|_g=|\eta|_g=1$. Since $g$ is geometrically Dolbeault formal, the $(1,1)$-form $\Phi=i\theta\wedge\eta$ is Dolbeault harmonic on $X$ with respect to $g$.  We fix the complex-valued function
$a:=\text{tr}_\omega\Phi=g(\eta,\overline{\theta})$. Since $H^{0,1}_{\bar\partial}(X)=\mathbb C[\bar\theta]$, there exists $c\in\C\setminus\{0\}$ such that $[\eta]=c[\overline{\theta}]$, i.e., $\eta-c\overline{\theta}$ is $\dbar$-exact. Since $\dbar$-exact forms are $L^2$-orthogonal to Dolbeault harmonic forms, we have
\begin{align*}
0=\langle \eta,\eta-c\overline{\theta}\rangle_{L^2}=\text{Vol}_g(X)-\overline{c}\int_X a\frac{\omega^2}{2},
\end{align*}
which shows that $a\not\equiv0$. Let $\hat g=f g$, with $f>0$, be the Gauduchon
representative in the conformal class of $g$. By Lemma~\ref{lemma:conform-inv},
$\Phi$ is harmonic also with respect to $\hat g$, and by
Lemma~\ref{lemma:dolbeault-middle-trace} its trace with respect to $\hat g$
is a constant $C$. Since $\text{tr}_{\hat\omega}\Phi=f^{-1}a$, we have
$a=Cf$. Thus $C\neq0$, and $a$ is nowhere zero.
Write $C=|C|e^{it}$. Replacing $\eta$ by $e^{-it}\eta$ and $c$ by
$e^{-it}c$, and retaining the same notation, gives
$a=|C|f>0$ and preserves $[\eta]_{\bar\partial}=c[\bar\theta]_{\bar\partial}$.
The norms are unchanged  by the change of phase, so the Cauchy--Schwarz inequality yields
$0<a\leq |\eta|_g|\theta|_g=1$.
The orthogonality identity now gives
\[
\overline c=\frac{\operatorname{Vol}_g(X)}{\int_Xa\,dV_g}>0.
\]
Hence $c$ is positive real and satisfies \eqref{eq:proj-lemma-norm-c}. Moreover,
$a\leq1$ implies $c\geq1$.
\end{proof}
We next establish a nonnegativity criterion for an intersection pairing, which will force $c=1$. This will be crucial for the proof of Proposition \ref{prop:elliptic-vertical-rigidity}.

\begin{lemma}
\label{lem:effective-euler-intersection}
Let $X$ be a compact connected complex surface, let $\theta$ be a nowhere-vanishing holomorphic one-form, and set
$\psi:=i\theta\wedge\bar\theta$.
Let $\beta$ be a smooth real semipositive $(1,1)$-form of rank
one with $\psi\wedge\beta>0$, and let $\rho\geq0$ be a smooth function
such that $\alpha:=\beta-\rho\psi$ is closed.
Denote by $Z$ the section of $\ker\beta\cap T^{1,0}X$
normalized by $\theta(Z)=1$.
Suppose a holomorphic line bundle $L$ admits a nonzero
holomorphic section $V$ and a smooth section $s$ satisfying,
for a smooth Hermitian metric $h$ and a constant $A>0$,
\[
\rho\leq A|s|_h^2,\qquad
|(\bar\partial_Ls)(\bar Z)|_h^2\leq A|s|_h^2.
\]
Then $c_1(L)\cdot[\alpha]\geq0$.
\end{lemma}
\begin{proof}
Let $D$ be the Chern connection of the Hermitian metric $h$ in the hypothesis. Set
\[
q=|s|_h^2,\qquad p=|V|_h^2,\qquad
\kappa=\frac{i}{2\pi}F_D.
\]
Let $e$ be a local unitary frame of $L$ and write $s=fe$, $V=ve$,
and $De=i\gamma\otimes e$, where $\gamma$ is a real one-form.
We use the abbreviations $Df:=df+if\gamma$ and $Dv:=dv+iv\gamma$,
so that $Ds=Df\otimes e$ and $DV=Dv\otimes e$.
In this frame, $q=|f|^2$, $p=|v|^2$, and $\kappa=-d\gamma/(2\pi)$.
For a section $\sigma=ue$, set $Du:=du+iu\gamma$ and
$j_\sigma:=\operatorname{Im}(\overline uDu)$.
This real one-form is globally defined, since intrinsically
$j_\sigma(Y)=\operatorname{Im}h(D_Y\sigma,\sigma)$ for every real tangent vector $Y$,
with $h$ taken to be linear in its first argument.
For $0<\varepsilon\leq1$, set
\[
Q_\varepsilon:=q+\varepsilon^2p,\qquad
R_\varepsilon:=Q_\varepsilon+\varepsilon^4,\qquad
j_\varepsilon:=j_s+\varepsilon^2j_V.
\]
The smooth real two-form
\[
T_\varepsilon:=\kappa+
\frac1{2\pi}d\left(\frac{j_\varepsilon}{R_\varepsilon}\right)
\]
is closed and satisfies $[T_\varepsilon]=c_1(L)$.
We will estimate $T_\varepsilon\wedge\alpha$ from below.

We first express $T_\varepsilon$ in terms of the sections $s$ and $V$.
Since $j_s=\operatorname{Im}(\overline fDf)
=\operatorname{Im}(\overline fdf)+q\gamma$, we have
\[
dj_s=i\,df\wedge d\overline f+dq\wedge\gamma+q\,d\gamma
=i\,Df\wedge\overline{Df}+q\,d\gamma.
\]
Applying the same calculation to $V$ gives
\[
dj_\varepsilon
=i\bigl(Df\wedge\overline{Df}+\varepsilon^2Dv\wedge\overline{Dv}\bigr)
-2\pi Q_\varepsilon\kappa.
\]
Set $B:=\overline fDf+\varepsilon^2\overline vDv$. Metric compatibility yields
\[
dR_\varepsilon=dQ_\varepsilon=B+\overline B,
\qquad
j_\varepsilon=\frac{B-\overline B}{2i}.
\]
Consequently, $dR_\varepsilon\wedge j_\varepsilon=iB\wedge\overline B$, and
\begin{align*}
T_\varepsilon
&=\kappa+\frac1{2\pi}
\left(\frac{dj_\varepsilon}{R_\varepsilon}
-\frac{dR_\varepsilon\wedge j_\varepsilon}{R_\varepsilon^2}\right)\\
&=\left(1-\frac{Q_\varepsilon}{R_\varepsilon}\right)\kappa+\frac{i}{2\pi R_\varepsilon^2}
\left[R_\varepsilon
\bigl(Df\wedge\overline{Df}+\varepsilon^2Dv\wedge\overline{Dv}\bigr)
-B\wedge\overline B\right].
\end{align*}
Here $1-Q_\varepsilon/R_\varepsilon=\varepsilon^4/R_\varepsilon$.
Using $R_\varepsilon=|f|^2+\varepsilon^2|v|^2+\varepsilon^4$ and expanding gives
\begin{align*}
&R_\varepsilon
\bigl(Df\wedge\overline{Df}+\varepsilon^2Dv\wedge\overline{Dv}\bigr)
-B\wedge\overline B\\
&\qquad=\varepsilon^2(vDf-fDv)\wedge\overline{(vDf-fDv)}
+\varepsilon^4Df\wedge\overline{Df}
+\varepsilon^6Dv\wedge\overline{Dv}.
\end{align*}
For an $L$-valued one-form $\zeta=\xi\otimes e$, we use the convention
$
h(\zeta\wedge\overline\zeta)
:=h(e,e)\,\xi\wedge\overline\xi.
$
Thus $h$ contracts the bundle factors in $L\otimes\overline L$,
leaving an ordinary complex-valued two-form. The same convention applies
to $L^{\otimes2}$ with its induced metric $h^{\otimes2}$.
Set
\[
W:=V\otimes Ds-s\otimes DV\in\mathcal A^1(X,L^{\otimes2}).
\]
Locally, $W=(vDf-fDv)\otimes(e\otimes e)$.
Since $e$ is unitary, we have obtained
\begin{equation}\label{eq:effective-euler-regularization}
\begin{split}
T_\varepsilon={}\frac{\varepsilon^4}{R_\varepsilon}\kappa&+\frac{i}{2\pi R_\varepsilon^2}
\left(\varepsilon^2h^{\otimes2}(W\wedge\overline W)
+\varepsilon^4h(Ds\wedge\overline{Ds})
+\varepsilon^6h(DV\wedge\overline{DV})\right).
\end{split}
\end{equation}
We now estimate $T_\varepsilon\wedge\alpha$, first considering the derivative
terms and then the curvature term. Throughout the estimates, $C>0$ denotes
a constant independent of $\varepsilon$ and may increase from line to line.

Set $d\mu:=\psi\wedge\beta$.
For an $L$-valued one-form $\zeta=\xi\otimes e$ in a unitary frame,
\[
i\,h(\zeta\wedge\overline\zeta)\wedge\beta
=i\,\xi\wedge\overline\xi\wedge\beta.
\]
At a point $x\in X$, write $\beta=i\chi\wedge\overline\chi$ for a
$(1,0)$-form $\chi$; this is possible because $\beta$ is semipositive of rank one.
Since $\psi\wedge\beta>0$, the forms $\theta,\chi$ are linearly independent, and
\[
T_x^*X\otimes_{\mathbb R}\mathbb C
=\operatorname{span}_{\mathbb C}
\{\theta_x,\chi_x,\overline\theta_x,\overline\chi_x\}.
\]
Write $\xi=a\theta+b\chi+c\overline\theta+d\overline\chi$ at $x$.
Wedging with $\beta$ eliminates all terms containing $\chi$ or $\overline\chi$,
so
\[
i\,\xi\wedge\overline\xi\wedge\beta
=i\bigl(|a|^2-|c|^2\bigr)\theta\wedge\overline\theta\wedge\beta.
\]
Since $\theta(Z)=1$ and $\beta(Z,\cdot)=0$, we have $\chi(Z)=0$,
$\xi(Z)=a$, and $\xi(\overline Z)=c$. Thus
\begin{equation}\label{eq:elliptic-formula}
i\,h(\zeta\wedge\overline\zeta)\wedge\beta
=\bigl(|\zeta(Z)|_h^2-|\zeta(\overline Z)|_h^2\bigr)\,d\mu.
\end{equation}
The same identity holds for $L^{\otimes2}$-valued one-forms with $h^{\otimes2}$.
Since $V$ is holomorphic, $D_{\overline Z}V=0$, and therefore
\[
W(\overline Z)=V\otimes D_{\overline Z}s,
\qquad
|W(\overline Z)|_{h^{\otimes2}}^2=p\,|D_{\overline Z}s|_h^2.
\]
Applying \eqref{eq:elliptic-formula} to \eqref{eq:effective-euler-regularization} gives
\begin{align*}
T_\varepsilon\wedge\beta
=\frac{\varepsilon^4}{R_\varepsilon}\kappa\wedge\beta+\frac1{2\pi R_\varepsilon^2}
\Bigl[\varepsilon^2|W(Z)|_{h^{\otimes2}}^2
+\varepsilon^4|D_Zs|_h^2+\varepsilon^6|D_ZV|_h^2-(\varepsilon^2p+\varepsilon^4)|D_{\overline Z}s|_h^2\Bigr]\,d\mu.
\end{align*}
By hypothesis, $|D_{\overline Z}s|_h^2
=|(\bar\partial_Ls)(\overline Z)|_h^2\leq Aq$.
The derivative contribution to $T_\varepsilon\wedge\beta$ is consequently
bounded below by $-CF_\varepsilon\,d\mu$, where
\[
F_\varepsilon:=\frac{q(\varepsilon^2p+\varepsilon^4)}{R_\varepsilon^2}.
\]
The inequality $4ab\leq(a+b)^2$ for $a,b\geq0$, applied with
$a=q$ and $b=\varepsilon^2p+\varepsilon^4$, gives
$0\leq F_\varepsilon\leq1/4$.
Moreover, $F_\varepsilon\to0$ pointwise, since it vanishes wherever $q=0$,
while at every point with $q>0$ its numerator tends to zero and its
denominator tends to $q^2>0$.

Let $g_0$ be the Hermitian metric with fundamental form $\psi+\beta$.
All norms below are pointwise norms computed using $g_0$ and the bundle
metrics induced by $h$. Its volume form is
$(\psi+\beta)^2/2=\psi\wedge\beta=d\mu$.
Since $s,V$ and $D$ are smooth and $X$ is compact, $|Ds|$ and $|DV|$
are bounded. Furthermore,
\[
|W|^2\leq2p|Ds|^2+2q|DV|^2\leq C(p+q).
\]
Denote the derivative component of $T_\varepsilon$ by
\[
E_\varepsilon:=\frac{i}{2\pi R_\varepsilon^2}
\left(\varepsilon^2h^{\otimes2}(W\wedge\overline W)
+\varepsilon^4h(Ds\wedge\overline{Ds})
+\varepsilon^6h(DV\wedge\overline{DV})\right).
\]
The same calculation as in \eqref{eq:elliptic-formula}, with $\beta$ and
$\psi$ interchanged, shows that
\[
\bigl|i\,h(\zeta\wedge\overline\zeta)\wedge\psi\bigr|_{g_0}
\leq C|\zeta|^2,
\]
and similarly for $h^{\otimes2}$. Applying this estimate to each summand gives
\[
|\rho E_\varepsilon\wedge\psi|_{g_0}
\leq C\frac{\rho}{R_\varepsilon^2}
\left[\varepsilon^2|W|^2+\varepsilon^4|Ds|^2+\varepsilon^6|DV|^2\right].
\]
Using $\rho\leq Aq$ and the preceding bounds, we obtain
\[
|\rho E_\varepsilon\wedge\psi|_{g_0}
\leq C\frac{q\bigl[\varepsilon^2(p+q)+\varepsilon^4+\varepsilon^6\bigr]}
{R_\varepsilon^2}
\leq C(F_\varepsilon+\varepsilon^2).
\]
Indeed, since $R_\varepsilon\geq q$ and $R_\varepsilon\geq\varepsilon^4$,
\begin{align*}
\frac{q\bigl[\varepsilon^2(p+q)+\varepsilon^4+\varepsilon^6\bigr]}
{R_\varepsilon^2}
&=F_\varepsilon
+\varepsilon^2\left(\frac q{R_\varepsilon}\right)^2
+\varepsilon^2\frac q{R_\varepsilon}\frac{\varepsilon^4}{R_\varepsilon}\leq F_\varepsilon+2\varepsilon^2.
\end{align*}

The curvature contribution to $T_\varepsilon\wedge\alpha$ is
$G_\varepsilon\kappa\wedge\alpha$, where
\[
G_\varepsilon:=\frac{\varepsilon^4}{R_\varepsilon},\qquad 0<G_\varepsilon\leq1.
\]
Since $\kappa$ and $\alpha$ are fixed smooth forms on $X$,
\[
|G_\varepsilon\kappa\wedge\alpha|_{g_0}\leq CG_\varepsilon.
\]
Set $N:=\{s=0\}\cap\{V=0\}$.
At points where $q>0$, $G_\varepsilon\to0$; where $q=0$ and $p>0$,
\[
G_\varepsilon=\frac{\varepsilon^2}{p+\varepsilon^2}\longrightarrow0.
\]
On $N$, one has $G_\varepsilon=1$ for every $\varepsilon$.
Since $V$ is holomorphic and not identically zero on the connected surface $X$,
its zero set has zero $d\mu$-measure. Thus $N$ also has measure zero, and
$G_\varepsilon\to0$ almost everywhere.
Combining the estimates in
\[
T_\varepsilon\wedge\alpha
=G_\varepsilon\kappa\wedge\alpha
+E_\varepsilon\wedge\beta-\rho E_\varepsilon\wedge\psi
\]
gives the following inequality of real top-degree forms:
\[
T_\varepsilon\wedge\alpha
\geq-C(G_\varepsilon+F_\varepsilon+\varepsilon^2)\,d\mu.
\]
Because $\alpha$ is closed, Stokes' theorem and $[\kappa]=c_1(L)$ imply
\[
c_1(L)\cdot[\alpha]=\int_XT_\varepsilon\wedge\alpha
\geq-C\int_X(G_\varepsilon+F_\varepsilon+\varepsilon^2)\,d\mu.
\]
The integrand on the right tends to zero almost everywhere and is bounded for $0<\varepsilon\leq1$. Dominated convergence therefore gives
$c_1(L)\cdot[\alpha]\geq0$.
\end{proof}

We recall some facts about the Chern classes of $X$ using its 
structure of $\mathbb{CP}^1$-bundle. Write $X=\mathbb P_E(\mathcal E)$, where $\mathcal E$ is a rank-two
holomorphic vector bundle over $E$. Let
$$
L:=T^{1,0}_{X/E}=\ker(d\pi:T^{1,0}X\longrightarrow\pi^*T^{1,0}E),
\qquad \nu:=c_1(L).
$$
Thus $L$ is the holomorphic line bundle tangent to the fibers.
A point $y\in X$ above $p\in E$ corresponds to a line
$\ell_y\subset\mathcal E_p$. These lines define the tautological bundle of $X$, i.e., the holomorphic
line subbundle $S\subset\pi^*\mathcal E$, with
$$
S_y=\ell_y\subset(\pi^*\mathcal E)_y=\mathcal E_p.
$$
Let $Q:=\pi^*\mathcal E/S$. This is a holomorphic line bundle
over $X$, with fiber $Q_y=\mathcal E_p/\ell_y$, and we have
the exact sequence
$$
0\longrightarrow S\longrightarrow\pi^*\mathcal E
\longrightarrow Q\longrightarrow0.
$$
Set $x:=c_1(Q)$.

Keeping $p$ fixed, varying $y$ amounts to varying the line
$\ell_y$ inside $\mathcal E_p$. After choosing a complementary
line, nearby lines are graphs of linear maps from $\ell_y$
to that complement. Thus the tangent space to the fiber at
$y$ is naturally
$$
L_y\simeq\operatorname{Hom}(\ell_y,\mathcal E_p/\ell_y)
=\operatorname{Hom}(S_y,Q_y).
$$
Globally, this identifies $L\simeq S^*\otimes Q$.
The exact sequence above gives
$c_1(S)=\pi^*c_1(\mathcal E)-x$, and hence
\begin{equation}\label{eq:projective-vertical-chern}
\nu=2x-\pi^*c_1(\mathcal E).
\end{equation}
Since $E$ is a curve, $c_2(\mathcal E)=0$. The Whitney formula
yields $0=\pi^*c_2(\mathcal E)=c_1(S)x$, so
\begin{equation}\label{eq:projective-relation}
x^2=x\,\pi^*c_1(\mathcal E).
\end{equation}
\begin{prop}
\label{prop:elliptic-vertical-rigidity}
Let $\pi:X\to E$ be an elliptic ruled surface such that
$T^{1,0}_{X/E}$ admits a nonzero holomorphic section.
Every geometrically Dolbeault formal Hermitian metric on $X$
is K\"ahler.
\end{prop}
\begin{proof}
Let $g$ be such a metric and use
Lemma~\ref{lem:ruled-harmonic-normalization}. Put
\[
\psi=i\theta\wedge\bar\theta,\qquad
\Phi=i\theta\wedge\eta,\qquad
\Sigma=\ast\Phi=a\omega-\bar\Phi,
\]
where $\ast$ is the $\C$-antilinear Hodge operator.
Both $\Phi$ and $\Sigma$ are Dolbeault harmonic and satisfy
$\Phi^2=\Sigma^2=0$. Thus
\[
\tau=\Phi-\Sigma=\Phi+\bar\Phi-a\omega,\qquad
\alpha=\psi-\tau
\]
are real closed $(1,1)$-forms: $\tau$ is real and
$\bar\partial$-closed, hence $d$-closed.

For each fiber $F\simeq\mathbb{CP}^1$, the restriction of the vertical holomorphic tangent space $L|_F$ coincides
with $T^{1,0}F$. Thus $\int_F\nu=\deg T^{1,0}F=2$.
The identities \eqref{eq:projective-vertical-chern} and
\eqref{eq:projective-relation} give $\nu^2=0$. Since $\theta$ is pulled
back from the base, $[\psi]^2=0$ and $\nu\cdot[\psi]>0$.
Consequently, $\nu$ and $[\psi]$ form a basis of $H^{1,1}(X;\R)$.

We now determine the class of the real closed form $c\psi-\tau$.
Since $\Sigma=\Phi-\tau$ and $[\Phi]_{\bar\partial}=c[\psi]_{\bar\partial}$,
we have $[\Sigma]_{\bar\partial}=[c\psi-\tau]_{\bar\partial}$.
Here $c$ is the positive constant supplied by
Lemma~\ref{lem:ruled-harmonic-normalization}.
Identifying Dolbeault and de Rham classes of type $(1,1)$ on the
K\"ahler surface $X$, the identity $\Sigma^2=0$ gives $[c\psi-\tau]^2=0$.
Write
$$
[c\psi-\tau]=u[\psi]+\lambda\nu,\qquad u,\lambda\in\R.
$$
Since $\theta$ is pulled back from the base, $\psi|_F=\Phi|_F=\bar\Phi|_F=0$.
Therefore $(c\psi-\tau)|_F=\Sigma|_F=a\omega|_F$, and
$$
2\lambda=\int_F(c\psi-\tau)=\int_F a\omega>0.
$$
On the other hand, $0=[c\psi-\tau]^2=2u\lambda\,\nu\cdot[\psi]$, so $u=0$.
Thus $[c\psi-\tau]=\lambda\nu$, and the definition $\alpha=\psi-\tau$ gives
\begin{equation}\label{eq:elliptic-rigidity-intersection}
[\alpha]=(1-c)[\psi]+\lambda\nu,\qquad
\nu\cdot[\alpha]=(1-c)\nu\cdot[\psi].
\end{equation}

We construct the data for Lemma~\ref{lem:effective-euler-intersection}
directly from $\theta$ and $\eta$. Set
$$
t:=\frac1a\geq1,\qquad \beta:=\alpha+(t-1)\psi.
$$
Thus $\beta$ is a globally defined real $(1,1)$-form. To check its
positivity, choose a local unitary coframe $(\theta,\chi)$ and write
$\eta=a\bar\theta+b\bar\chi$, with $a^2+|b|^2=1$. Using
$\alpha=\psi-\Phi-\bar\Phi+a\omega$ and completing the square gives
$$
\beta=ia\left(\chi-\frac ba\theta\right)
\wedge\overline{\left(\chi-\frac ba\theta\right)}.
$$
Since $a>0$, this form is semipositive of rank one, and
$\psi\wedge\beta=a\omega^2/2>0$.
Choose local holomorphic coordinates $(z,w)$ adapted to the ruling,
with $\theta=dz$. The properties of $\beta$ give a unique expression
$$
\beta=iH(dw+A\,dz)\wedge(d\bar w+\bar A\,d\bar z),\qquad H>0.
$$
Define $\xi:=dw+A\,dz$ and $Z:=\partial_z-A\partial_w$.
The two expressions for $\beta$ show that there is a smooth
nowhere-vanishing function $v$ such that
$$
\chi-\frac ba\theta=v\xi,\qquad H=a|v|^2.
$$
Define $r:=b\bar v$. Since $\eta=a\bar\theta+b\bar\chi$ and
$a^2+|b|^2=1$, we obtain
\begin{equation}\label{eq:elliptic-harmonic-identities}
\eta=t\bar\theta+r\bar\xi,\qquad \frac{|r|^2}{H}=t-\frac1t.
\end{equation}
The metric defined by $\beta$ on the vertical tangent line
raises $\eta|_{\overline{L}}$ to a global smooth section
\[
s=\frac rH\partial_w\in \Gamma(X,L),\qquad
|s|_\beta^2=t-t^{-1}\geq t-1.
\]
Thus $\rho=t-1\leq C_h|s|_h^2$ for every fixed smooth metric
$h$ on $L$. 
The equation $\bar\partial\eta=0$ yields
\[
\bar Zr=t_{\bar w}+(\bar A)_{\bar w}r.
\]
Differentiating the second identity in
\eqref{eq:elliptic-harmonic-identities} gives
\[
(1+t^{-2})t_{\bar w}
=\frac{r_{\bar w}\bar r+r\,\overline{r_w}}H
-\frac{|r|^2H_{\bar w}}{H^2}.
\]
Choose finitely many coordinate neighborhoods $V_i$ and relatively
compact open subsets $U_i\Subset V_i$ such that the $U_i$ cover $X$.
On each $\overline U_i$, the smooth local coefficients and their first
derivatives are bounded, and $H$ has a positive lower bound. Since there
are finitely many charts, the constants can be chosen uniformly. Hence \(t_{\bar w}=O(|r|)\) and
\(\bar Z(r/H)=O(|r|)\). Since \(\partial_w\) is a holomorphic
frame of \(L\), this proves
\[
|(\bar\partial_Ls)(\bar Z)|_h^2\leq C_h|s|_h^2.
\]
Lemma~\ref{lem:effective-euler-intersection} and
\eqref{eq:elliptic-rigidity-intersection} now imply
\[
0\leq\nu\cdot[\alpha]=(1-c)\nu\cdot[\psi].
\]
Thus $c\leq1$, whereas Lemma \ref{lem:ruled-harmonic-normalization} gives $c\geq1$.
It follows that $c=1$ and $\int_X(1-a)\,dV_g=0$, so
$a\equiv1$ and $\eta=\bar\theta$. The closed form $\alpha$
then equals $\omega-\psi$. Since $d\psi=0$, we obtain
$d\omega=0$.
\end{proof}

\begin{prop}
\label{thm:elliptic-ruled-metric-rigidity}
Every geometrically Dolbeault formal Hermitian metric on an
elliptic ruled surface is K\"ahler.
\end{prop}
\begin{proof}
Let $X=\mathbb{P}_E(\mathcal{E})$ be a $\mathbb{CP}^1$-bundle over an elliptic curve $E$. By the classification of ruled surfaces over an elliptic curve,
$X$ is either decomposable or isomorphic, as a projective bundle over $E$,
to one of the two indecomposable ruled surfaces $S_0=\mathbb{P}_E(\mathcal{E}_0)$ and $S_1=\mathbb P_E(\mathcal{E}_1)$, see
\cite[Section~2.2]{Diaw21} and \cite[Section~3]{LorayMarin09}.

If $\mathcal E=A\oplus B$ is split, scaling the second summand
defines a holomorphic one-parameter group of automorphisms of $X$ by
$$
\Phi_t(p,[u:v]):=(p,[u:e^t v]),\qquad t\in\mathbb C,
$$
where $u\in A_p$ and $v\in B_p$. Since $\pi\circ\Phi_t=\pi$,
the infinitesimal generator of this action $V$ is a holomorphic
vertical vector field. In a local affine coordinate $w=v/u$, we can write  $V=w\,\partial/\partial w$, so it is not identically zero.

For $X=S_0$, consider the nonsplit exact
sequence
$$
0\longrightarrow\mathcal O_E
\xrightarrow{\iota}\mathcal E_0
\xrightarrow{q}\mathcal O_E\longrightarrow0.
$$
Set $N:=\iota\circ q$. Since $q$ is surjective and $\iota$ is
injective, $N$ is nonzero. Moreover, $q\circ\iota=0$ gives
$N^2=0$. Thus $\mathrm{id}+tN$ is a holomorphic bundle
automorphism with inverse $\mathrm{id}-tN$, and these maps
induce a holomorphic one-parameter group
$
\Phi_t(p,[v])=(p,[v+tNv])$, $t\in\mathbb C,
$
on $X$. Its infinitesimal generator is a holomorphic vertical
vector field, which is not identically zero because $N$ is
nonzero and nilpotent, hence not a scalar endomorphism.

In both of the above cases, the elliptic ruled surface admits a non-zero holomorphic vector field. Hence 
Proposition~\ref{prop:elliptic-vertical-rigidity} applies and every geometrically Dolbeault formal metric is K\"ahler.

The remaining case, i.e., the stable indecomposable type $S_1$, has projective monodromy
generated by $w\mapsto-w$ and $w\mapsto1/w$
\cite[Theorem~2.21]{Diaw21}. This is a Klein four group, so pulling back $X$
via the multiplication-by-two map $[2]:E\to E$ yields a trivial projective bundle $Y$. In particular, we obtain a
degree-four holomorphic unramified covering
\[
p:Y=E\times\mathbb {CP}^1\longrightarrow X.
\]
Let $g$ be a geometrically Dolbeault formal Hermitian metric on
$X$. Since $p$ is a local isometry for $p^*g$, pullback
commutes with the Dolbeault Laplacian and is injective on harmonic
forms. The projective bundle formula gives, for every elliptic
ruled surface,
\[
h^{p,q}(X)=h^{p,q}(E)+h^{p-1,q-1}(E).
\]
Thus $X$ and $Y$ have equal
Hodge numbers, and every $p^*g$-harmonic form is pulled back from
a $g$-harmonic form. Their wedge products are therefore Dolbeault harmonic,
so $p^*g$ is geometrically Dolbeault formal on $Y$.
The trivial ruling has a nonzero holomorphic vertical vector field, given by $V=w\partial /\partial w$, in local holomorphic coordinates $(z,w)$. Therefore,  
Proposition~\ref{prop:elliptic-vertical-rigidity} makes $p^*g$
K\"ahler. Since $p^*(d\omega)=0$, the metric $g$ is K\"ahler
as well on $X$.
\end{proof}

\begin{cor}\label{cor:elliptic-ruled-classification}
Let $X$ be an elliptic ruled surface. If $X$ admits a geometrically Dolbeault formal
Hermitian metric, then its ruling is projectively unitary.
\end{cor}
\begin{proof}
Let $X$ be an elliptic ruled surface. If $X$ admits a geometrically Dolbeault formal Hermitian metric $g$, then, by
Proposition \ref{thm:elliptic-ruled-metric-rigidity}, the Hermitian metric $g$ is K\"ahler. We prove that $X$ is projectively unitary. Indeed, for a K\"ahler metric, geometric Dolbeault formality is equivalent
to Riemannian geometric formality. If $\theta$ is the pullback of
a nonzero holomorphic one-form on the base, then
\cite[Theorem~1.3]{GN} gives $\nabla^g\theta=0$.
Its kernel and orthogonal complement are parallel and
$J$-invariant. Since the fibers are simply connected, the complex de Rham theorem yields the holomorphic and isometric  
splitting of the universal cover
\[
(\tilde X,\tilde g)
\cong(\mathbb C,g_{\mathrm{flat}})\times(\mathbb {CP}^1,h).
\]
The deck group acts by translations on the first factor and
holomorphic isometries on the second. The group
$\operatorname{Isom}_{\mathrm{hol}}(\mathbb P^1,h)$ is compact,
hence conjugate in $PGL(2,\mathbb C)$ to a subgroup of $PU(2)$.
\end{proof}

\subsubsection{$k(X)=0$}
The surfaces in this case are precisely the hyperelliptic surfaces
\cite[Chapter~V]{BPV}. Making use of their description as K\"ahler solvmanifolds, it was shown in \cite{ST24} that they admit geometrically formal K\"ahler metrics.
See also
Theorem~\ref{thm:classification-kahler} for an argument exploiting their curvature properties.

\subsubsection{$k(X)=1$}

We determine the structure of the geometrically Dolbeault formal surfaces with $b_1(M)=2$, $b_2(M)=2$, and $k(X)=1$.
\begin{prop}\label{prop:kappa-one-dolbeault}
Let $X$ be a compact K\"ahler surface satisfying
\[
b_1(X)=2,\qquad h_{\dbar}^{1,1}(X)=2,\qquad k(X)=1.
\]
If $X$ admits a geometrically Dolbeault formal Hermitian metric, then
\[
X\simeq (E\times C)/G,
\]
where:
\begin{enumerate}
    \item $E$ is an elliptic curve and $C$ is a compact curve of genus
    $g(C)\geq2$;
    \item $G$ is a finite group acting freely and diagonally on
    $E\times C$;
    \item the action of $G$ on $E$ is by translations;
    \item $C/G\simeq\mathbb{P}^1$.
\end{enumerate}
\end{prop}

\begin{proof}
Assume that $X$ admits a geometrically Dolbeault formal Hermitian
metric. By Remark \ref{rmk:proj}, $X$ is
projective. Since $h^{1,0}(X)=1$,
its Albanese variety $A$ is an elliptic curve. A nonzero holomorphic
one-form on $X$ has constant positive norm by
Theorem~\ref{thm:st25}. Since this form is pulled back from $A$,
the Albanese map is a submersion. Its Stein
factorization \cite[Chapter III, Corollary 11.5]{Har} gives a smooth
fibration $p\colon X\to B$ with connected fibers and a finite
unramified map $B\to A$. In particular, $B$ is elliptic.

Let $C$ be a fiber of $p$. If $g(C)=0$, then $X$ is ruled, contrary
to $k(X)=1$. If $g(C)=1$, the canonical bundle formula \cite[Chapter~V]{BPV} gives
$K_X\simeq p^*L$ with $\deg L=\chi(X)=0$.
Since the fibers are connected,
$H^0(X,K_X^{\otimes m})=H^0(B,L^{\otimes m})$ has dimension
at most one for every $m$, contradicting $k(X)=1$. Thus $g(C)\geq2$. By \cite[Theorem~0.1(ii)]{VZ}, a smooth
projective family of curves of genus at least two over an elliptic
curve is isotrivial. The Fischer--Grauert theorem \cite{FG} makes
$p$ a holomorphic fiber bundle. Since $\operatorname{Aut}(C)$ is
finite, the finite unramified Galois cover $E\to B$ associated with
the kernel of its monodromy trivializes the bundle. Consequently,
\[
X\simeq(E\times C)/G,
\]
where $G$ is the deck group, acting diagonally and by translations
on $E$. Its action on $E$, and hence on the product, is free.
Moreover,
\[
H_{\dbar}^{1,0}(X)
 \simeq H_{\dbar}^{1,0}(E)\oplus H_{\dbar}^{1,0}(C)^G.
\]
Since $h^{1,0}(X)=1$, we obtain
$H_{\dbar}^{1,0}(C)^G\simeq H_{\dbar}^{1,0}(C/G)=0$,
so $C/G\simeq\mathbb P^1$.

\end{proof}

\subsection{Classification}\label{sec:kahler-classification}

We can now combine the preceding cases to obtain the following theorem.

\begin{thm}\label{thm:classification-kahler}
Let $X=(M,J)$ be a compact K\"ahler surface. Then $X$ is geometrically Dolbeault formal if and only if $X$ is either
\begin{enumerate}
\item a complex two-dimensional torus,
\item $\mathbb{CP}^2$ or $\mathbb{CP}^1\times\mathbb{CP}^1$,
\item a fake projective plane,
\item a fake quadric $\Gamma\backslash(\mathbb D\times\mathbb D)$, where
$\Gamma\subset\operatorname{Aut}(\mathbb D)\times
\operatorname{Aut}(\mathbb D)$ is a torsion-free cocompact lattice,
\item a projectively unitary $\mathbb{CP}^1$-bundle over an elliptic curve,
\item a hyperelliptic surface,
\item a quotient $(E\times C)/G$, where $E$ is an elliptic curve,
$g(C)\geq2$, and a finite group $G$ acts freely and diagonally,
by translations on $E$, with $C/G\simeq\mathbb{CP}^1$.
\end{enumerate}
In particular, the following are equivalent:
\begin{enumerate}
\item $X$ admits a geometrically Dolbeault formal Hermitian metric,
\item $X$ admits an ABC-geometrically formal Hermitian metric,
\item $X$ admits a geometrically formal K\"ahler metric.
\end{enumerate}
\end{thm}

\begin{proof}

Suppose $X$ admits a geometrically Dolbeault formal Hermitian metric.
Corollary~\ref{cor:dolb-cohom-ring-kahler-1} shows that $X$ is minimal
and that $b_1(X)\in\{0,2,4\}$.

If $b_1(X)=4$,
Lemma~\ref{lemma:geom-dolb-2} proves that $X$ is a complex torus.

If $b_1(X)=0$, Proposition~\ref{prop:geom-dolb-kahler-b_1=0} gives
$k(X)\in\{-\infty,2\}$. In the rational case, minimality leaves
$\mathbb{CP}^2$ and the Hirzebruch surfaces, and
Proposition~\ref{prop:no-hirzebruch} excludes all the latter except
$\mathbb{CP}^1\times\mathbb{CP}^1$. In the general-type case,
$b_2(X)\in\{1,2\}$, we have established  in \ref{subsection:k=2b1=0} and
Proposition~\ref{prop:fake-quadric-uniformization} that $X$ is, respectively, a fake
projective plane or a factor-preserving bidisk fake quadric.

Finally, if $b_1(X)=2$, the three possible Kodaira dimensions give
the projectively unitary elliptic ruled surfaces of
Corollary~\ref{cor:elliptic-ruled-classification}, hyperelliptic
surfaces, and the quotients of
Proposition~\ref{prop:kappa-one-dolbeault}. This proves that $X$ belongs
to the stated list.

Conversely, let us construct a geometrically formal K\"ahler
metric (thus a geometrically Dolbeault formal Hermitian metric) on each surface in the list.

We first treat the cases admitting K\"ahler metrics with
nonnegative curvature operator, i.e., cases (1), (2), (5), and (6).

First, observe that on complex tori, we can consider flat K\"ahler metrics, whereas on $\mathbb{CP}^2$ and
$\mathbb{CP}^1\times\mathbb{CP}^1$, we can take, respectively, the
Fubini--Study metric and the product of Fubini--Study metrics. These are Hermitian symmetric metrics with nonnegative curvature operator.

If $X$ is a hyperelliptic surface $(E\times F)/G$, with $E$ and $F$ elliptic curves and $G$ a finite group acting holomorphically and diagonally, we can choose $G$-invariant flat metrics on $E$ and $F$.
Their product descends to a flat K\"ahler metric on $X$.

If $X$ is a projectively unitary elliptic ruled surface, write
$$
X=(\mathbb C\times\mathbb{CP}^1)/\Lambda,
\qquad
\lambda\cdot(z,p)=(z+\lambda,\rho(\lambda)p),
\qquad
\rho:\Lambda\to\mathrm{PU}(2).
$$
The product of the Euclidean and Fubini--Study metrics is
invariant under the monodromy action and descends to a K\"ahler metric
on $X$. Its curvature operator is nonnegative, since the first
factor is flat, the second has positive Gaussian curvature,
and all mixed curvature terms vanish.

Thus, in all these cases, by
\cite[Theorem~10]{STprelim} the harmonic forms are parallel and
the K\"ahler metric is geometrically formal.

(3) We have proved in Proposition \ref{prop:fake-pr-dolb} that every K\"ahler metric on a fake projective plane is geometrically formal.

(4) If $X$ is a fake quadric $
X=\Gamma\backslash(\mathbb D\times\mathbb D)$ and $
\Gamma\subset\operatorname{Aut}(\mathbb D)\times
\operatorname{Aut}(\mathbb D)
$ as above, the product hyperbolic metric descends to a K\"ahler metric $g$,
and the fundamental forms of the hyperbolic metric on each factor descend to linearly
independent parallel forms $\Omega_1,\Omega_2$.
Since $b_1=b_3=0$ and $b_2=2$,
$
\mathcal H^\bullet_{\mathrm{dR}}(X;\mathbb C)
=\mathbb C\langle
1,\Omega_1,\Omega_2,\Omega_1\wedge\Omega_2
\rangle.
$
Hence all harmonic forms are parallel and $g$ is geometrically
formal.

(7) Let $X$ be a properly elliptic surface
$X=(E\times C)/G$ as above. Choose a flat metric on $E$
and the hyperbolic metric of curvature $-1$ on $C$.
Translations preserve the first metric, and holomorphic
automorphisms preserve the second. Their product therefore
descends to a K\"ahler metric $g$ on $X$.
A nonzero parallel holomorphic one-form on $E$ is invariant
under translations and descends to a parallel form $\theta$.
The fundamental forms of the fixed metrics on each factor also descend to parallel forms
$\omega_E,\omega_C$.
Since $b_1(X)=b_2(X)=2$, the linearly independent forms
$\theta,\bar\theta$ span the harmonic one-forms, and
$\omega_E,\omega_C$ span the harmonic two-forms.
Hodge duality then shows that the harmonic three-forms are
parallel as well. Thus all harmonic forms are parallel,
and $g$ is geometrically formal.

Thus,  every element of the list admits a geometrically formal K\"ahler metric and in particular, if $X$ admits a geometrically Dolbeault formal Hermitian metric, then it also admits a geometrically formal K\"ahler one, i.e.,  $(1)\Rightarrow (3)$. Conversely, a geometrically formal K\"ahler metric is also geometrically Dolbeault formal and ABC-geometrically formal, i.e., $(3)\Rightarrow (1),(2)$. 
 Finally, by Lemma \ref{lemma:abc-harmonic-comparison}, every ABC-geometrically formal metric is geometrically Aeppli formal  on a compact K\"ahler manifold, and in particular geometrically Dolbeault formal. Hence,  $(2)\Rightarrow (3)$.

\end{proof}

\section{Metric rigidity}\label{sec:rigidity}
The classification result (Theorem \ref{thm:classification-kahler}) shows that the existence of a geometrically Dolbeault or ABC-geometrically formal Hermitian metric is equivalent to the existence of a geometrically formal K\"ahler metric on compact K\"ahler surfaces.
For ABC-geometrically formal metrics, the equivalence holds at the level of each metric.
\begin{prop}\label{lemma:aeppli-forces-kahler}
Let $X$ be a compact K\"ahler surface. A Hermitian metric  $g$
on $X$ is ABC-geometrically formal if and only if it
is K\"ahler and geometrically formal.
\end{prop}
\begin{proof}
Suppose $g$ is ABC-geometrically formal. Since $X$ satisfies the
$\del\dbar$-lemma,
Lemma~\ref{lemma:abc-harmonic-comparison} makes $g$ geometrically
Aeppli formal, and in particular, geometrically Dolbeault formal. If $b_1(X)=4$, Theorem~\ref{thm:geom-BC-A}
shows that $g$ is flat and K\"ahler.

Otherwise, Lemmas~\ref{lemma:geom-dolb-5}
and~\ref{lemma:geom-dolb-7} give $b^+=1$ and
$b^-\in\{0,1\}$. Choose a nonzero real self-dual de Rham harmonic
two-form $\Phi$. Since $h^{2,0}(X)=h^{0,2}(X)=0$,
Lemma~\ref{lemma:abc-harmonic-comparison} gives
\[
\mathcal H^2_{dR}(X;\mathbb R)
=\mathcal H^{1,1}_{BC}(X)\cap\mathcal A^2_{\mathbb R}(X).
\]
Thus $\Phi$ is a real self-dual $(1,1)$-form, so
$\Phi=\lambda\omega$ for a real smooth function $\lambda$.
Its constant positive norm forces $|\lambda|$ to be a positive
constant. Since $X$ is connected, we have $\lambda\equiv c\ne0$.
Since $d\Phi=0$, it follows that $d\omega=0$.
For a K\"ahler metric, ABC-geometric formality and Riemannian
geometric formality agree. This also proves the converse.
\end{proof}

\subsection{Rigidity of geometrically Dolbeault formal metrics}
Aside from the case of projectively unitary bundles, 
the spaces of geometrically Dolbeault formal and geometrically formal K\"ahler metrics may differ on compact K\"ahler surfaces. Nonetheless, we have the following result.
\begin{prop}\label{prop:dolbeault-lck-rigidity}
Let $X$ be a compact K\"ahler surface. A geometrically Dolbeault
formal Hermitian metric on $X$ which is locally conformally
K\"ahler is K\"ahler.
\end{prop}
\begin{proof}
By \cite[Theorem~2.1]{Vaisman80}, such a metric $g$ is globally
conformally K\"ahler, because $X$ admits a K\"ahler metric.
Write $\widetilde g=u g$ for a K\"ahler representative. Its
K\"ahler form $\widetilde\omega$ is Dolbeault harmonic for $g$
by Lemma~\ref{lemma:conform-inv}. From $\widetilde\omega=u\omega$, we have $
|\widetilde\omega|_g^2=2u^2$, which
is constant since $g$ is geometrically Dolbeault formal. Thus $u$ is constant and $g$ is K\"ahler.
\end{proof}

We can also compare the notion of geometric Dolbeault formality with that of Gauduchon metrics.

\begin{prop}\label{prop:gauduchon-b1-zero}
Let $X$ be a compact K\"ahler surface with $b_1(X)=0$.
Every geometrically Dolbeault formal Gauduchon Hermitian metric on
$X$ is K\"ahler.
\end{prop}
\begin{proof}
Let $g$ be such a metric. Proposition \ref{prop:geom-dolb-kahler-b_1=0}  gives
$h^{2,0}(X)=0$, $b^+=1$, and $b^-\in\{0,1\}$. Let $V:=\mathcal{H}_g^-\otimes \C$. By
Lemma~\ref{lemma:asd-dolb}, $V\subset\mathcal H_{\bar\partial}^{1,1}(X)$ and, since $h^{1,1}(X)=1+b^-$, its $L^2$-orthogonal complement $V^{\perp}$ is a complex line in $\mathcal H_{\bar\partial}^{1,1}(X)$.
Now, $V^{\perp}$ is preserved by the complex-antilinear Hodge involution $\ast$, so we can choose a nonzero harmonic form $\eta\in V^{\perp}$ with $\ast\eta=\eta$.

By Lemma~\ref{lemma:dolbeault-middle-trace},
$
\eta=c\omega+i\phi,
$
where $c\in \R$ and $\phi$ is a real primitive, hence, anti-self-dual $(1,1)$-form.
Geometric Dolbeault formality implies constant pointwise norms and inner
products of Dolbeault harmonic forms, so the orthogonality to $V$ is pointwise and
$$
|\phi|_g^2=|\eta|_g^2-2c^2
$$
is constant. If $c=0$, then $\bar\partial\phi=0$. Since $\phi$ is real,
it is closed, and anti-self-duality makes it de Rham harmonic.
Thus $\eta\in V$, a contradiction. Hence $c\ne0$.

Suppose $\phi\ne0$. Its constant norm makes it nowhere-vanishing.
If $b^-=0$, then $b_2=1$, $c_2(X)=3$, and $c_1(X)^2=9$.
Applying Lemma~\ref{lem:primitive-complement-chern} to $\phi$ gives
complex line bundles with Chern classes $x,y$ satisfying
$$
\langle(x-y)^2,[X]\rangle=9-4\cdot3=-3.
$$
This contradicts positive definiteness of the intersection form.
If $b^-=1$, choose a nonzero real anti-self-dual harmonic form $\beta$.
Its norm is a positive constant and it is pointwise orthogonal to $\phi$.
Lemma~\ref{lem:primitive-complement-chern} would give $c_1(X)^2=4c_2(X)$.
But $c_1(X)^2=8$ and $c_2(X)=4$, another contradiction.

Therefore $\phi=0$. Since $c$ is a nonzero constant,
$\bar\partial\eta=0$ implies $\bar\partial\omega=0$. The form $\omega$
is real, so $d\omega=0$.
\end{proof}

The following proposition treats geometrically Dolbeault formal Gauduchon metrics on both hyperelliptic surfaces and the quotients
of Kodaira dimension one. Recall that a hyperelliptic surface has the
same quotient description with $C$ an elliptic curve \cite[Chapter~V]{BPV}.

\begin{prop}\label{lem:gauduchon-product-quotient}
Let $X=(E\times C)/G$, where $E$ is an elliptic curve, $C$ is a compact
curve of genus at least one, and $G$ a finite group which acts freely and diagonally, by
translations on $E$, with $C/G\simeq\mathbb P^1$.
Every geometrically Dolbeault formal Gauduchon Hermitian metric on $X$
is K\"ahler.
\end{prop}
\begin{proof}
Let $q:E\times C\to X$ be the quotient map. Translations act trivially
on the cohomology of $E$, whereas
$H^{1,0}(C)^G=H^{1,0}(C/G)=0$. The K\"unneth formula therefore gives
$h^{1,0}(X)=1$ and $h^{1,1}(X)=2$.
Let $g$ be a geometrically Dolbeault formal Gauduchon metric on $X$,
and choose $\theta,\eta,a$ as in
Lemma~\ref{lem:ruled-harmonic-normalization}. Since $g$ is Gauduchon, the trace of the harmonic form
$\Phi:=i\theta\wedge\eta$ is constant by
Lemma~\ref{lemma:dolbeault-middle-trace}. Thus $a$ and $t:=a^{-1}\geq1$
are constants.

Set
$$
\psi:=i\theta\wedge\bar\theta,\qquad
\Sigma:=\ast\Phi=a\omega-\bar\Phi,\qquad
\tau:=\Phi-\Sigma=\Phi+\bar\Phi-a\omega.
$$
The form $\tau$ is real and $\bar\partial$-closed, hence closed.
Also $d\psi=0$, since $\theta$ is closed. Therefore
$$
\beta:=t\psi-\tau
$$
is a globally defined closed real $(1,1)$-form.
In a local unitary coframe $(\theta,\chi)$, write
$\eta=a\bar\theta+b\bar\chi$, with $a^2+|b|^2=1$. Then
$$
\beta=ia\left(\chi-\frac ba\theta\right)
\wedge\overline{\left(\chi-\frac ba\theta\right)}.
$$
Thus $\beta$ is semipositive of rank one and $\psi\wedge\beta>0$.

Suppose $t>1$ and put $\zeta:=\eta-t\bar\theta$. Directly from the
preceding expression,
$$
\beta=\frac{i}{t-t^{-1}}\bar\zeta\wedge\zeta,\qquad
|\zeta|_g^2=t^2-1>0,\qquad \bar\partial\zeta=0.
$$
Let $\sigma$ be the form on $X$ induced by a $G$-invariant positive
area form on $C$. The classes $[\psi],[\sigma]$ form a basis of
$H^{1,1}(X;\R)$, and
$$
[\psi]^2=[\sigma]^2=0,\qquad [\psi]\cdot[\sigma]>0.
$$
Writing $[\beta]=u[\psi]+v[\sigma]$, the inequalities
$[\psi]\cdot[\beta]>0$ and $[\beta]^2=0$ give $v>0$ and $u=0$.
Consequently, for every $p\in C$,
$$
\int_{E\times\{p\}}q^*\beta=0.
$$
Semipositivity makes the restriction of $q^*\beta$ to this elliptic
curve vanish pointwise. The formula in terms of $\zeta$ then shows
that $q^*\zeta$ annihilates the tangent directions of $E$.

For each fixed $p\in C$, the form $q^*\zeta$ can therefore be viewed as
a map from $E$ to the complex line $(T_p^{0,1}C)^*$. The equation
$\bar\partial(q^*\zeta)=0$ makes this map holomorphic. Compactness of
$E$ makes it constant. Hence
$$
q^*\zeta=\operatorname{pr}_C^*\zeta_C
$$
for a smooth nowhere-vanishing $(0,1)$-form $\zeta_C$ on $C$.
Deck invariance implies that $\zeta_C$ is $G$-invariant.

Since $g(C)\geq1$, the quotient map $C\to C/G\simeq\mathbb P^1$
is branched. Thus some element of the effective action on $C$ fixes
a point $p$ and has derivative $\lambda\ne1$ there. Its action on
$(T_p^{0,1}C)^*$ is multiplication by $\bar\lambda$, so invariance forces
$\zeta_C(p)=\bar\lambda\zeta_C(p)$, and hence $\zeta_C(p)=0$.
This contradicts the fact that $\zeta_C$ never vanishes.

Therefore $t=1$, so $a=1$ and equality in Cauchy--Schwarz gives
$\eta=\bar\theta$. Then $\Phi=\psi$ and
$\Sigma=\omega-\psi$. Since $\Sigma$ is real and $\bar\partial$-closed,
it is closed. Together with $d\psi=0$, this gives $d\omega=0$.
\end{proof}
Combining the classification theorem (Theorem \ref{thm:classification-kahler}) with the preceding rigidity results yields a uniform rigidity
statement. 

\begin{thm}\label{thm:kahler-gdf-metric-rigidity}
Let $X$ be a compact K\"ahler surface and let $g$ be a geometrically
Dolbeault formal Hermitian metric on $X$. The following are equivalent:
\begin{enumerate}
\item $g$ is locally conformally K\"ahler,
\item $g$ is Gauduchon,
\item $g$ is K\"ahler.
\end{enumerate}
\end{thm}

\begin{proof}
The implication $(1)\Rightarrow(3)$ is
Proposition~\ref{prop:dolbeault-lck-rigidity}. A K\"ahler metric is
Gauduchon and locally conformally K\"ahler, so it remains to prove
$(2)\Rightarrow(3)$.

Assume $g$ is Gauduchon. If $b_1(X)=4$, Lemma~\ref{lemma:geom-dolb-2} shows that $g$ is flat
and K\"ahler. We may therefore assume that $b_1(X)\in\{0,2\}$. By Theorem~\ref{thm:classification-kahler},
there are two cases. If $b_1(X)=0$, apply
Proposition~\ref{prop:gauduchon-b1-zero}. If $b_1=2$, then $X$ is  either elliptic ruled, a hyperelliptic surface, or a properly elliptic surface of Proposition~\ref{prop:kappa-one-dolbeault}. For the first case, every geometrically Dolbeault formal Hermitian metric is K\"ahler by Proposition~\ref{thm:elliptic-ruled-metric-rigidity}, hence it holds in particular for geometrically Dolbeault formal Gauduchon metrics. For the remaining surfaces, we can apply
Proposition~\ref{lem:gauduchon-product-quotient}. Thus $g$ is K\"ahler in
every case.
\end{proof}

\subsection{Conformal classes of geometrically Dolbeault formal metrics}
We finish the metric rigidity discussion with a study of conformal classes of geometrically Dolbeault formal metrics on $\mathbb C\mathbb P^2$ and fake projective planes, showing that a geometrically Dolbeault formal Hermitian metric need not be K\"ahler, in contrast with the very special situation on projectively unitary elliptic ruled surfaces.

First, in analogy with  \cite[Proposition 2.4]{OP}, we prove the following general fact.
\begin{prop}\label{prop:conformal-gdf-uniqueness}
Let $X$ be a compact K\"ahler surface. Each Hermitian conformal class
contains at most one geometrically Dolbeault formal metric up to
constant scaling.
\end{prop}
\begin{proof}
Suppose $g$ and $\widehat g=ug$, with $u>0$, are both geometrically
Dolbeault formal. Since $X$ is K\"ahler, $h^{1,1}(X)>0$. Choose
$0\ne\gamma\in\mathcal H_{\bar\partial}^{1,1}(X,g)$.
Lemma~\ref{lemma:conform-inv} makes $\gamma$ harmonic also for
$\widehat g$. Its squared norms for both metrics are positive constants,
and
$
|\gamma|_{\widehat g}^2=u^{-2}|\gamma|_g^2.
$
Thus $u$ is constant.
\end{proof}

\begin{prop}\label{thm:conformal-gdf-b2-one}
Let $X$ be $\mathbb{CP}^2$ or a fake projective plane.
Every Hermitian conformal class on $X$ contains a unique geometrically
Dolbeault formal metric up to constant scaling.
\end{prop}
\begin{proof}
The only nonzero Hodge numbers of $X$ are
$$
h^{0,0}(X)=h^{1,1}(X)=h^{2,2}(X)=1,
$$
and $b^-(X)=0$.
Let us fix a Hermitian conformal class and choose its Gauduchon
representative $h$, with fundamental form $\omega_h$.
The complex-antilinear Hodge operator $\ast_h$ preserves the complex
line $\mathcal H_{\bar\partial}^{1,1}(X,h)$ and squares to the identity
there. We may therefore choose a nonzero generator $\eta$ satisfying
$\ast_h\eta=\eta$.
By Lemma~\ref{lemma:dolbeault-middle-trace},
$
\eta=c\omega_h+i\phi,
$
where $c\in\mathbb R$ is constant and $\phi$ is a real primitive
$(1,1)$-form. If $c=0$, then $\bar\partial\eta=0$ implies
$\bar\partial\phi=0$. Since $\phi$ is real and anti-self-dual, it is
closed and de Rham harmonic. This contradicts $b^-(X)=0$ and
$\eta\ne0$. Thus $c\ne0$, and rescaling $\eta$ gives
$
\eta=\omega_h+i\phi.
$
Since $\phi$ is primitive and anti-self-dual, we have that
$
\omega_h\wedge\phi=0$ and $
\phi\wedge\phi=-1/2({|\phi|_h^2})\omega_h^2.
$
Consequently,
$$
\eta\wedge\eta
=\left(1+\frac{|\phi|_h^2}{2}\right)\omega_h^2.
$$
Define
$$
\lambda:=\left(1+\frac{|\phi|_h^2}{2}\right)^{1/2},
\qquad g:=\lambda h,
\qquad \omega_g=\lambda\omega_h.
$$
The function $\lambda$ is smooth and positive.
By Lemma~\ref{lemma:conform-inv}, $\eta$ is Dolbeault harmonic
for $g$, and the preceding identity becomes
$$
\eta\wedge\eta=\omega_g^2.
$$
Since $\omega_g^2$ is Dolbeault harmonic, the Hodge numbers give
$
\mathcal H_{\bar\partial}^{\bullet,\bullet}(X)
=\operatorname{span}_{\mathbb C}\{1,\eta,\omega_g^2\}.
$
This space is closed under wedge product, so $g$ is geometrically
Dolbeault formal. Uniqueness up to constant scaling follows from
Proposition~\ref{prop:conformal-gdf-uniqueness}.
\end{proof}

\begin{cor}\label{cor:nonkahler-gdf-b2-one}
On $\mathbb{CP}^2$ and every fake projective plane, there exist
geometrically Dolbeault formal Hermitian metrics which are neither
Gauduchon nor locally conformally K\"ahler, and hence are not K\"ahler.
\end{cor}
\begin{proof}
Let $\omega_0$ be a K\"ahler form on $X$. Choose a smooth $(0,1)$-form
$\xi$ such that $d(\partial\xi+\bar\partial\bar\xi)\ne0$.
For example, in a holomorphic coordinate chart, one may take
$\xi:=\chi|z_2|^2d\bar z_1$, where $\chi$ is a smooth cutoff supported
in the chart and equal to one on a smaller open set, and extend
$\xi$ by zero.
For sufficiently small nonzero real $\varepsilon$, the form
$
\omega_\varepsilon
:=\omega_0+\varepsilon(\partial\xi+\bar\partial\bar\xi)
$
is positive and satisfies
$\partial\bar\partial\omega_\varepsilon=0$ but
$d\omega_\varepsilon\ne0$.
Thus its associated Hermitian metric $h_\varepsilon$ is Gauduchon
and not K\"ahler.
By uniqueness of the Gauduchon representative up to constant scaling,
the conformal class of $h_\varepsilon$ contains no K\"ahler metric.
Proposition~\ref{thm:conformal-gdf-b2-one} supplies a geometrically
Dolbeault formal metric in this class. It is not K\"ahler, so
Theorem~\ref{thm:kahler-gdf-metric-rigidity} shows that it is neither
Gauduchon nor locally conformally K\"ahler.
\end{proof}
\subsection*{Declaration on the use of generative AI}
ChatGPT (OpenAI) was used to assist with exposition and the review of the proofs.
The author is responsible for the final text and the correctness and originality of the
mathematical content.

\end{document}